\documentclass[12pt]{amsart}
\usepackage{amssymb}
\ifx\pdfoutput\undefined
\usepackage{graphicx}  \DeclareGraphicsExtensions{.ps}

\else
\usepackage{amsfonts}
\usepackage[pdftex]{graphicx}  \DeclareGraphicsExtensions{.pdf,.mps}

\fi

\newcommand{\Op}{{\operatorname{Op}^{{w}}_h}}
\newcommand{\Oph}{{\operatorname{Op}_h}}

\newcommand{\CC}{{\mathbb C}}
\newcommand{\CI}{{\mathcal C}^\infty }
\newcommand{\CIc}{{\mathcal C}^\infty_{\rm{c}} }

\newcommand{\ZZ}{{\mathbb Z}}

\newcommand{\RR}{{\mathbb R}}

\newcommand{\NN}{{\mathbb N}}

\newcommand{\vol}{\operatorname{vol}}
\newcommand{\rank}{\operatorname{rank}}

\newcommand{\supp}{\operatorname{supp}}

\newcommand{\Spec}{\operatorname{Spec}}
\newcommand{\Res}{\operatorname{Res}}

\newcommand{\rest}{\!\!\restriction}

\renewcommand{\Re}{\mathop{\rm Re}\nolimits}
\renewcommand{\Im}{\mathop{\rm Im}\nolimits}

\newcommand{\esssupp}{\operatorname{ess\hspace{0.1cm}supp}}
\newcommand{\ad}{\operatorname{ad}}

\newcommand{\sigs}{\sigma_{\Sigma, h, \tilde h }}
\newcommand{\pxx}{\partial_x^\alpha \partial_\xi^\beta }

\newcommand{\Pkh}{\Psi_{\Sigma,\frac12}^{{\mathfrak p }}}

\newcommand{\Ssgh}{S_{\Sigma, \frac12}^{{\mathfrak p}'}}

\newcommand{\Ssh}{S_{\Sigma,\frac12}^{{\mathfrak p }}}
\newcommand{\Ssd}{S_{\Sigma,\delta}^{{\mathfrak p }}}

\newcommand{\WFh}{\operatorname{WF}_{ h }}
\newcommand{\Ops}{\operatorname{Op}_{\Sigma, h, \tilde h }}

\newcommand{\Pkgh}{\Psi_{\Sigma, \frac12}^{{\mathfrak p}' }}

\newcommand{\Opt}{\widetilde{\operatorname{Op}}_{h, \tilde h }}

\newcommand{\dn}{\|}
\newcommand{\hph}{{\widehat \varphi}}

\theoremstyle{plain}

\newtheorem{thm}{Theorem}

\newtheorem{prop}{Proposition}[section]

\newtheorem{lem}[prop]{Lemma}

\theoremstyle{definition}

\numberwithin{equation}{section}

\def\squarebox#1{\hbox to #1{\hfill\vbox to #1{\vfill}}} 
\newcommand{\stopthm}{\hfill\hfill\vbox{\hrule\hbox{\vrule\squarebox 
                 {.667em}\vrule}\hrule}\smallskip}

\title
[Fractal upper bounds on the density of resonances]
{Fractal upper bounds on the density of semiclassical resonances}
\author[J. Sj{\"o}strand]{Johannes Sj{\"o}strand}
\address{Centre de Math{\'e}matiques, {\'E}cole Polytechnique \\
UMR 7460, CNRS \\
F-91128 Palaiseau}
\email{johannes@math.polytechnique.fr}
\author[M. Zworski]{Maciej Zworski}
\address{Mathematics Department, University of California \\
Evans Hall, Berkeley, CA 94720}
\email{zworski@math.berkeley.edu}

\begin{document}    
   
\maketitle   
   
\section{Introduction}   
\label{in}

Let 
$ P = -h^2 \Delta_g + V ( x) $
 be a self-adjoint 
Schr\"odinger operator on 
a compact Riemannian $n$-manifold, $ ( X, g ) $, $ V \in \CI ( X ; \RR) $.
The spectral asymptotics as $ h \rightarrow 0 $ are given by the
celebrated {\em Weyl law} -- see \cite{DiSj} and \cite{Ivr} for
recent advances and numerous references. If we assume that 
the zero energy surface is nondegenerate,
\[ p \stackrel{\rm{def}}{=} |\xi|_g^2 + V( x ) = 0 \ \Longrightarrow \ dp \neq 0 \,, \]
then 
\begin{equation}
\label{eq:1}
 | {\rm Spec ( P ) } \cap D( 0 , C h ) | = {\mathcal O} ( h^{-n+1} ) \,,
\end{equation}
where $ D ( 0 , r ) = \{ z \in \CC \; : \; |z| < r  \}$, though of 
course in this case the eigenvalues are all real -- 
see \S \ref{pe} for yet another proof of this
well known result. 

Let $ H_p $ be the  Hamilton vector field of $ p $ on $ T^* X $, 
locally given by 
$$ H_p = \sum_{ j=1}^n \frac{\partial p}{\partial \xi_j} 
\partial_{ x_j } - 
\frac{\partial p}{\partial x_j} 
{\partial_{ \xi_j} } \,, \ \ ( x , \xi ) \in T^* \RR^n \,. $$ 
When the flow, $ \exp t H_p \; :\; p^{-1} ( E ) \rightarrow p^{-1} ( E ) $,
has the property that the set of its closed orbits has Liouville measure
zero on $ p = 0 $, then we have the infinitesimal version of the Weyl law:
\begin{equation}
\label{eq:1'}
 | {\rm Spec ( P ) } \cap D( 0 , C  h ) | = \frac{2 C h } 
{ ( 2 \pi h )^n } \int_{ p( x, \xi ) = 0 } d {\mathcal L} ( x , \xi ) + 
o ( h^{-n + 1 } ) \,,\end{equation}
where $ d {\mathcal L} $ 
is the Liouville measure on $ p = 0 $, that is 
$ d {\mathcal L} dp = dx d\xi  $. This result is the mathematical 
starting point of many recent investigations, mostly in physical 
literature,  of the finer structure
of the spectrum and its relation to classical dynamics -- see \cite{Bog} and
references given there.

When the manifold is non-compact the situation is dramatically different.
The simplest case is that of a manifold which is Euclidean outside of a compact
set and $ V + 1 \in \CIc ( X ; \RR ) $. The discrete eigenvalues of 
$ P $ are replaced by {\em quantum resonances} which are defined as
the poles of the meromorphic continuation of 
\[  ( P -  z)^{-1} \; : \; \CIc ( X ) \ \longrightarrow \ \CI ( X ) \,, \ \ 
\Im z > 0 \,,\]
and we denote the set of resonances by $ \Res ( P ( h ) ) $.
The basic physical interpretation
 is that a resonance at $ z = E_0 - i \Gamma/2 $
corresponds to a state with time evolution given by 
$ \exp ( - i t E_0 / h - \Gamma t /2 h ) $, and to a Breit-Wigner
peak in energy density, $ \Gamma / ( ( E - E_0)^2 + \Gamma^2 / 4 ) $.
This intuitive picture becomes however 
complicated when many resonances are present 
which is natural for $ h$ is small 
-- see \cite{NaStZw} and \cite{BoSj},\cite{BrPe},\cite{PZ2}
for recent 
results and references. 

\begin{figure}[htbp]
\begin{center}
\includegraphics[width=4.0in]{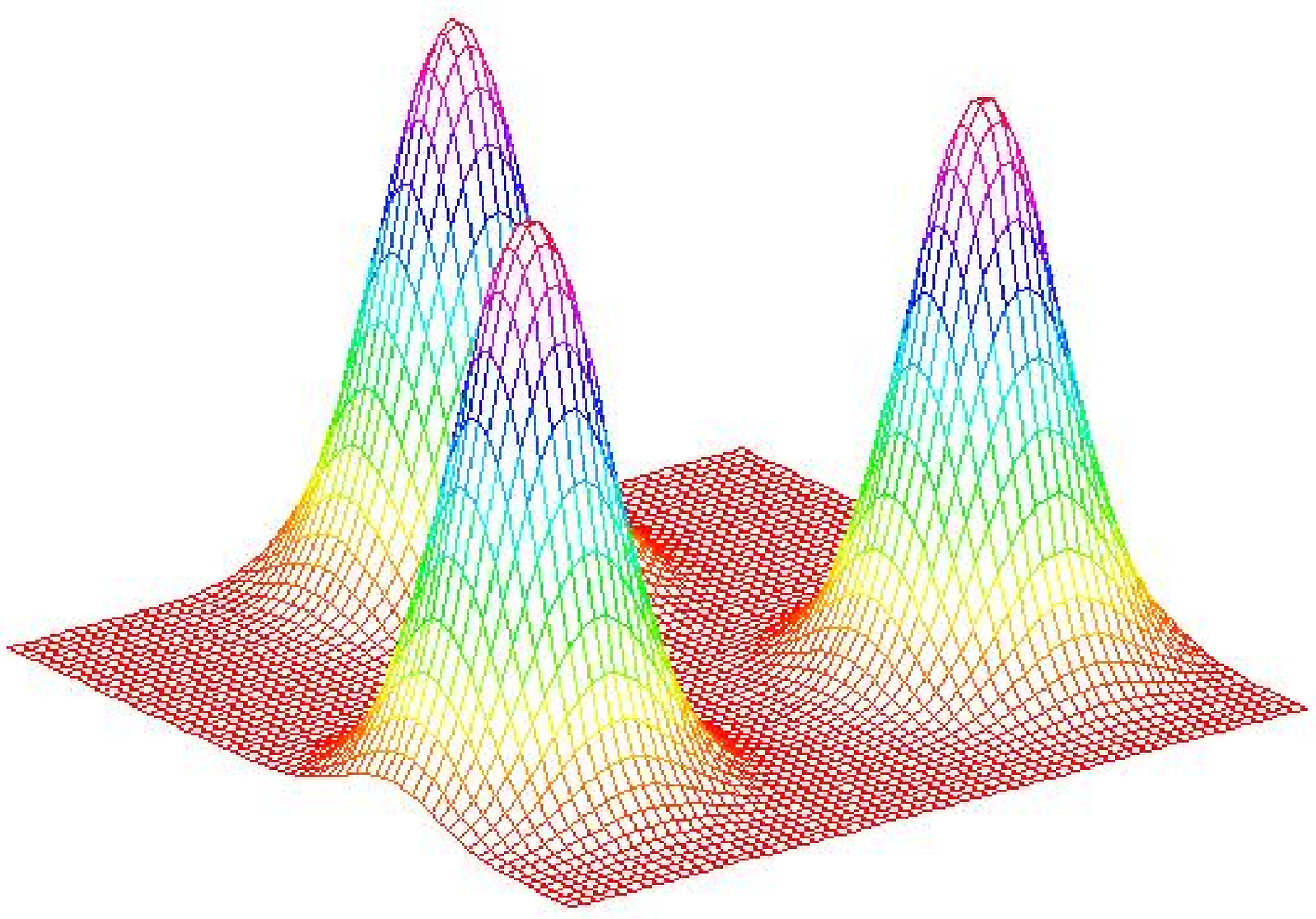}
\end{center}
\caption
{\label{fig:1} A three bump potential exhibiting hyperbolic dynamics
at an interval of energies.}
\end{figure}

Here we provide upper bounds for the number of resonances
of $ P $ in $ D( 0 ,C h ) $. The main result (Theorem 3) states that
for classical Hamiltonians $ p $ with hyperbolic flow on $ p = 0 $,
\begin{equation}
\label{eq:2}
  | \Res P ( h) \cap D ( 0 , C h ) | = {\mathcal O} ( h^{-\nu } )\,, 
\end{equation}
where $ 2 \nu + 1 $ is essentially the dimension of the trapped
(non-wandering) set in $ p^{-1} ( 0 ) $, 
\[ K \stackrel{\rm{def}}{=} \{ ( x , \xi ) \in T^* X \; : \; 
p ( x , \xi ) = 0 \,, \ \exp ( t H_p ) ( x , \xi ) \not \rightarrow 
\infty \,, \ \ t \rightarrow \pm \infty \}\,.\]
In the case of a compact manifold $ \nu = n -1 $ 
so that \eqref{eq:2} reduces to \eqref{eq:1}. By dimension we always
mean the Minkowski dimension 
\[ m_0 = 2 n - 1  - \sup\{ d \; : \; \limsup_{ \epsilon \rightarrow 0 }
\epsilon^{-d} \vol ( \{ \rho \in p^{-1} ( 0 ) \; : \; d ( \rho , K ) < 
\epsilon \} ) < \infty \} \,. \] 
A simple example is provided by a three bump potential shown in 
Fig.\ref{fig:1}.
\bigskip

The first estimate of this form was proved by the first author 
in \cite[Theorems 4.6, 5.5, and  5.7]{SjDuke}: there exists
constants $ C_0, C_1 > 0 $, such that for $ \delta_0 >0 $ fixed and
small enough
\begin{gather}
\label{eq:sjduke}
\begin{gathered}
| \Res ( P ( h ) ) \cap \{ z \; : \; \ |z | < \delta \,, \ \  
\Im z > - \mu \} |
\leq  C_1  \delta 
\left( \frac{h}{\mu} \right)^{-n} \mu^{- \frac12 \widetilde
m } \,, \\
C_0h \leq \mu \leq 1/C_0 \,, \ \ 
C_0 h^{\frac12}  \leq \delta  \leq \delta_0 \,, \ \ 0 < h < 1/C_0 \,,
\end{gathered}
\end{gather}
where now $ \widetilde
m $ is any number greater than the dimension of the trapped set in
$$ p^{-1} ( [ - 2 \delta_0  ,  2 \delta_0 ] ) \,. $$
 In homogeneous situations, such as
for instance obstacle scattering, $ \widetilde m = m + 1 $.
When $ \mu = C_0 h $, the improvement in Theorem \ref{t:3} 
lies in allowing  $ 
\delta \simeq h $, which is the natural limit for this type of 
spectral estimates. 

Earlier, non-geometric, bounds
on the number of resonances (scattering poles) were obtained
by Melrose \cite{M1},\cite{M2} and the second author \cite{Z1},\cite{Z2}.
In the case of convex co-compact Schottky quotients (and any convex co-compact 
quotients in dimension two) the analogue of 
\eqref{eq:2} was proved in \cite{GLZ}
using zeta function techniques, improving earlier estimates of \cite{ZwIn} 
the proof of which 
was largely based on \cite{SjDuke}. These technique gave similar results 
for the zeros
of zeta functions of rational maps \cite{Ch},\cite{StZw}, in which case the
dimension of the trapped set becomes essentially the dimension of the 
Julia set.

\begin{figure}[ht]
$$\includegraphics[width=7.5cm]{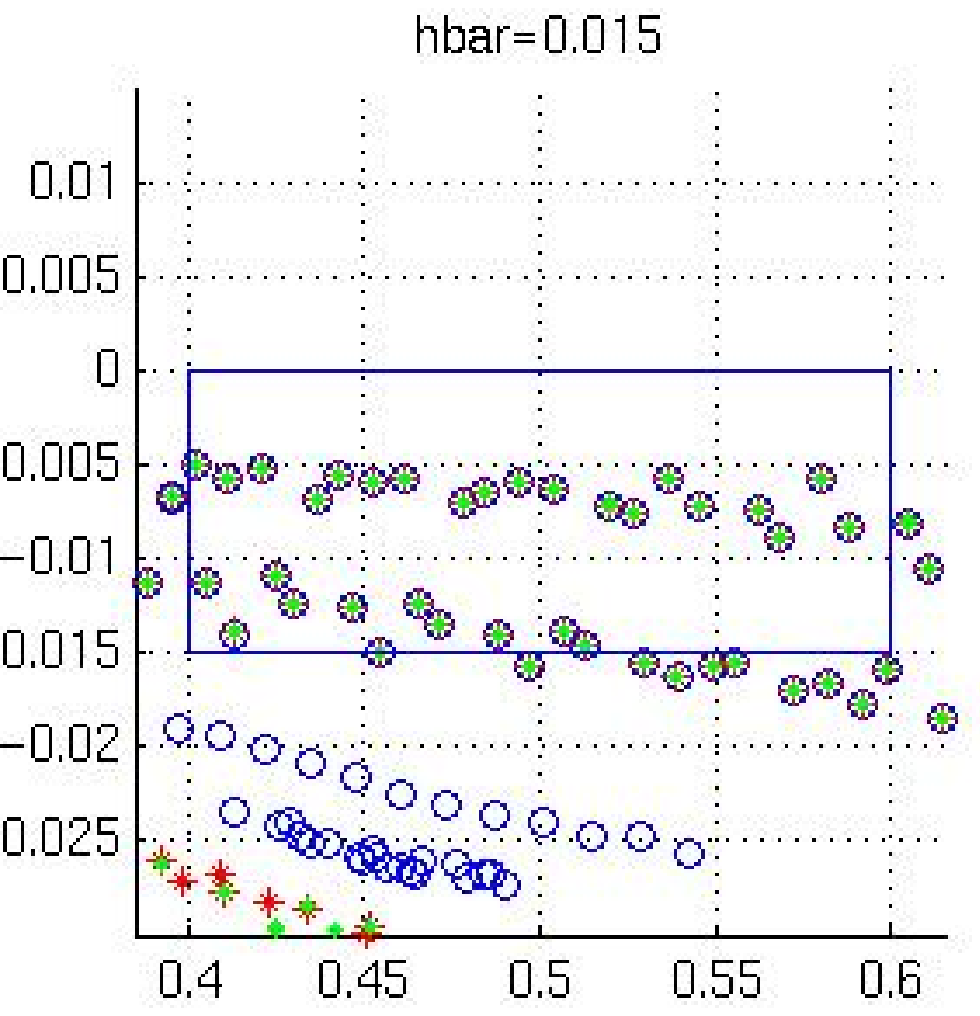}\qquad
\includegraphics[height=7.5cm]{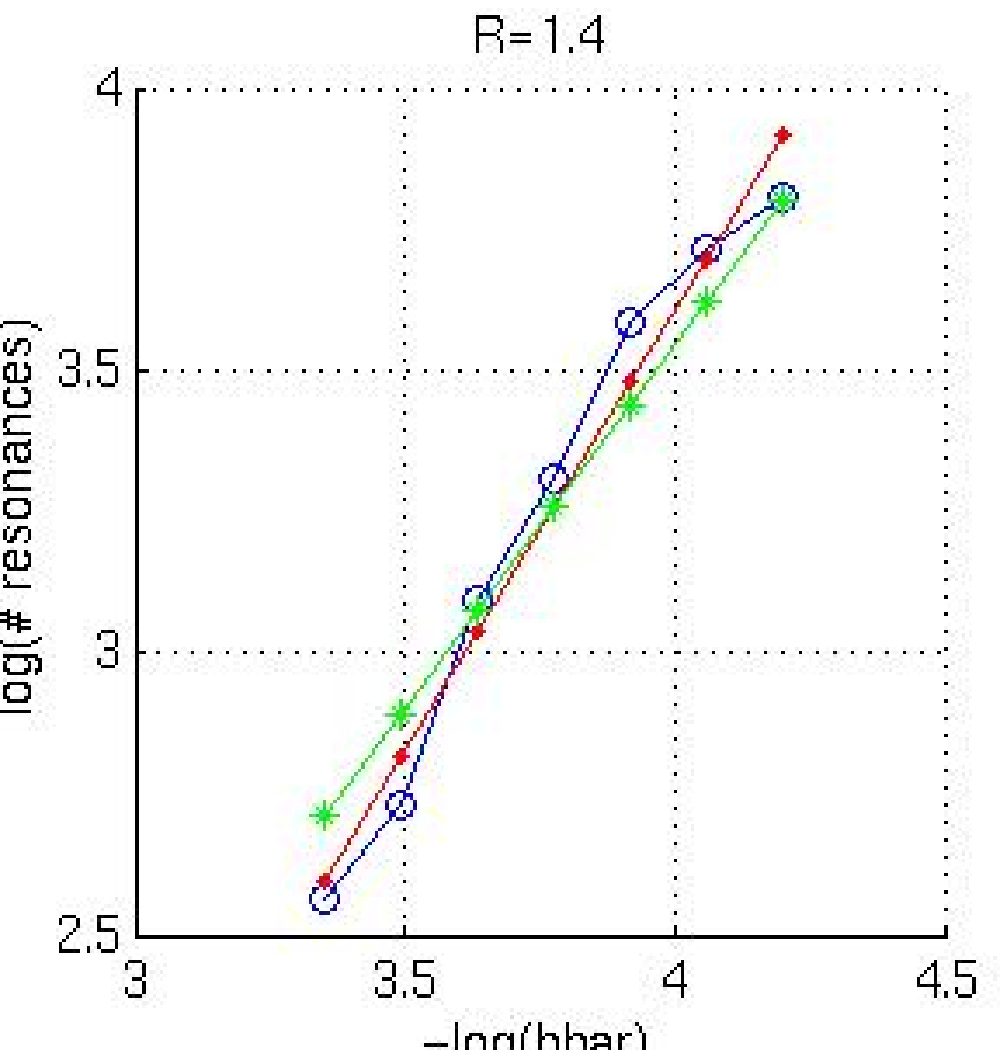}$$
\caption{A sample of results of \cite{L}: the plot on the left shows
resonances for $ h = 0.015 $, and the plot on the right is the
log-log plot of the number of resonances vs. $ h$ with $ \circ $
denoting numerical data, $\ast$ the density predicted by the
fractal Weyl law, and $ \bullet $, the least square
interpolants. }
\label{f:2}
\end{figure}

\begin{figure}[htbp]
$$\includegraphics[width=7.5cm]{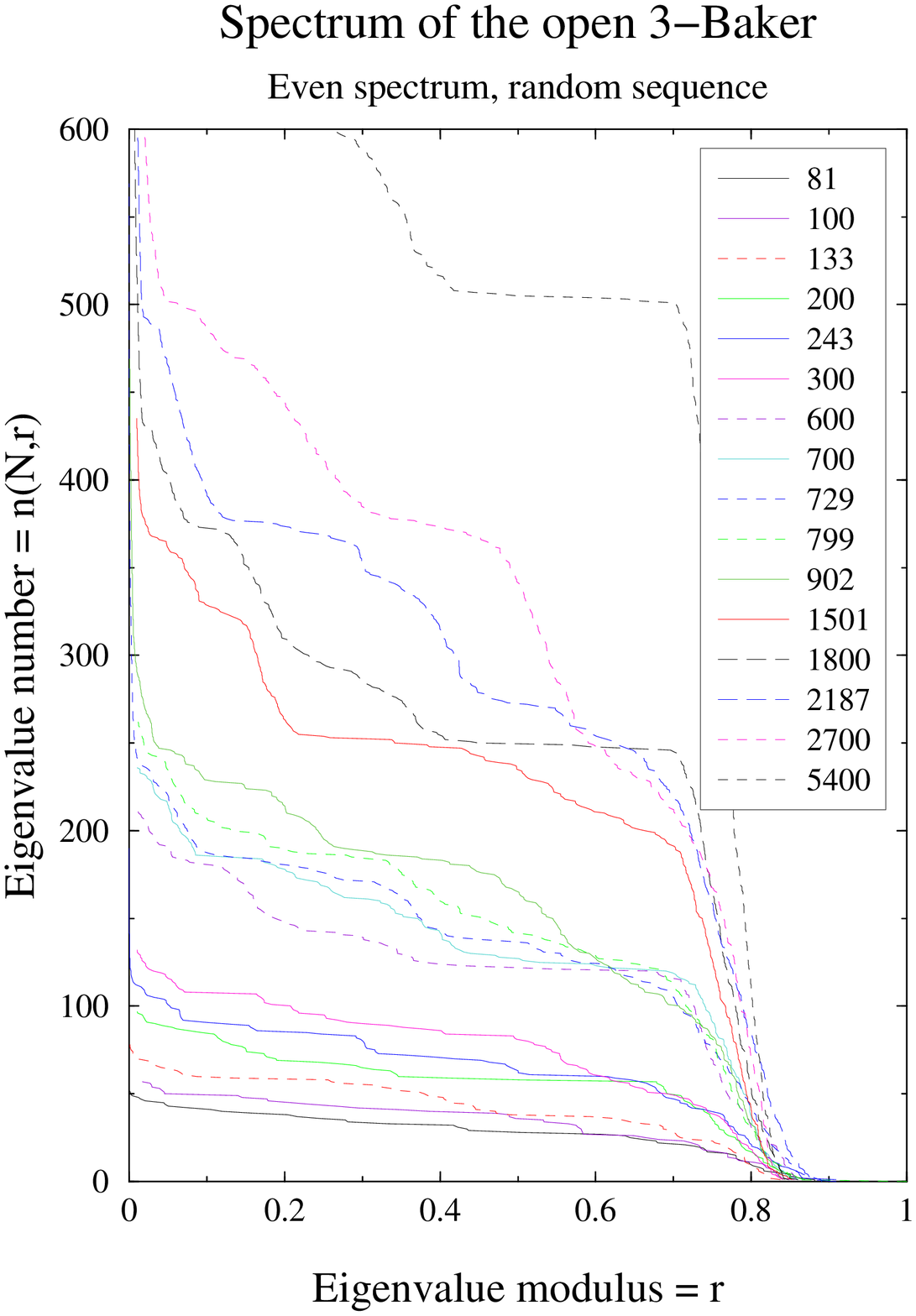}\qquad
\includegraphics[width=7.9cm]{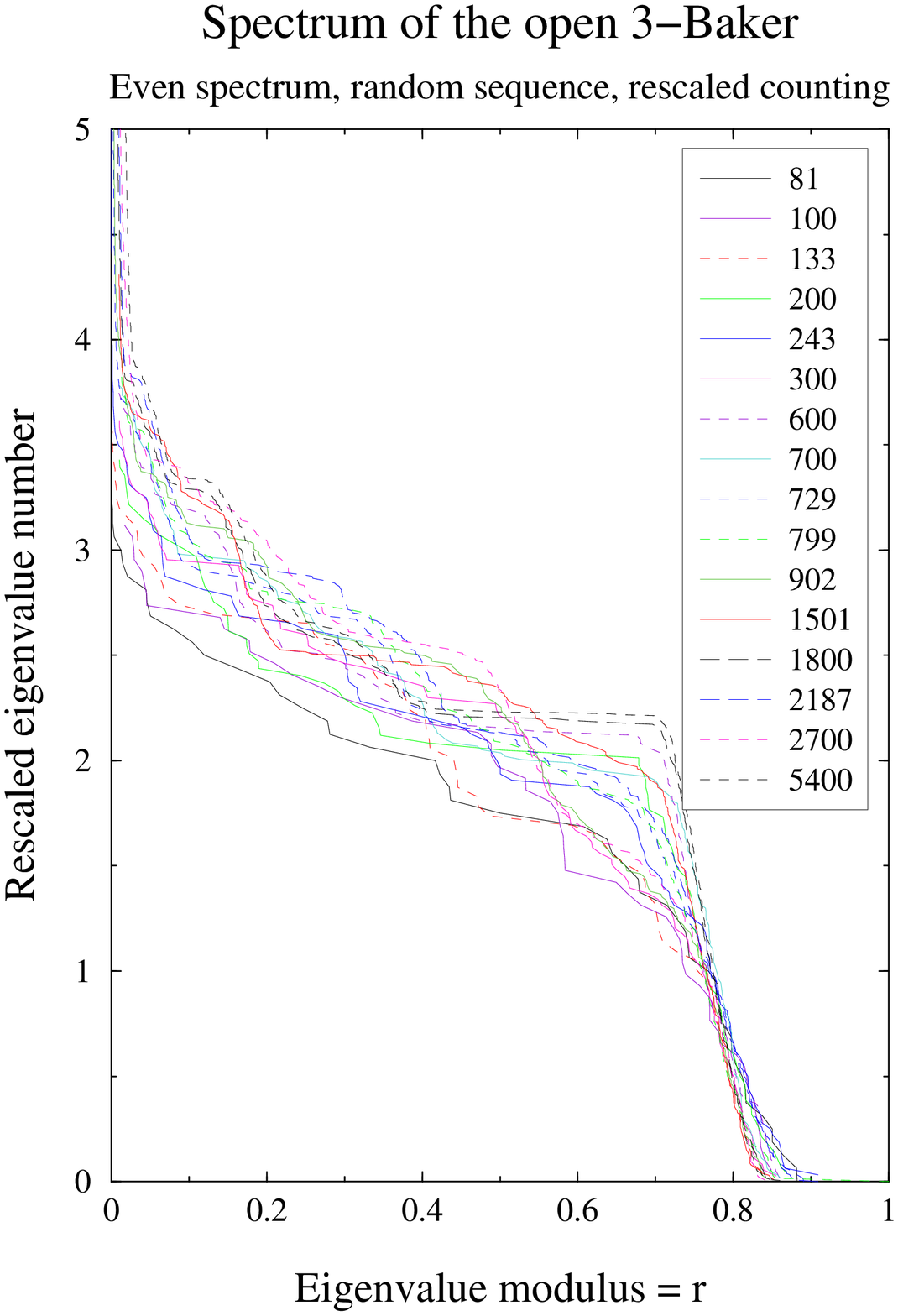}$$
\caption{\label{f:1} An example of numerical results of \cite{NoZ}. 
The eigenvalues of a $ 3 N \times 3 N $ matrix $ A_{3N} $ model the 
resonances  
in a disc of
size $ h \sim 1/N$. Dashed lines corresponding to 
 counting the resonances (eigenvalues of $ A_{3N} $)
with $ \exp ( \Im z /h ) \sim r  $. 
Full lines give the 
 counting function rescaled  using the dimension of the trapped 
set in this problem. The coalescence of the graphs confirms the 
power law based on the dimension of the trapped set. A similar 
correction based on \cite{LZ},\cite{LSZ},\cite{SjDuke} 
was given for the random matrix models of open system in
\cite{SchTw}.}
\end{figure}

Numerical investigations in different settings of semiclassical three bump
potentials \cite{L},\cite{LZ} (see Fig.\ref{f:2}), Schottky quotients \cite{GLZ}, three disc scattering
\cite{LSZ}, and Cantor-like Julia sets for $ z \mapsto z^2 + c $, $ c < - 2$ \cite{StZw},
suggest that for $ \mu \simeq C h $ and $ \delta \simeq 1 $ the estimate
\eqref{eq:sjduke} is optimal. A different model was recently considered
in \cite{NoZ}: quantum resonances were defined using an open quantum 
map with a classical ``trapped set'' corresponding to $ K $ intersected
with a hypersurface transversal to the flow -- see Fig.\ref{fig:3}. 
The numerical results and a simple linear algebraic toy model suggest
that the fine estimate \eqref{eq:2} is optimal -- see Fig.\ref{f:1}. 
A similar model
was also used in \cite{SchTw} where the fractal Weyl law gave 
corrections to the applications of random matrix theory to open 
quantum systems.

We now state the general assumptions on the operator $ P $.
We reiterate that the simplest case to keep in mind is 
\[ P = - h^2 \Delta + V( x ) - 1  \,, \ \ V \in \CIc ( \RR^n ) \,.\]
In general we consider
\[ P ( h ) \in \Psi^{0,2} ( X ) \,, \ \ P ( h ) = P ( h)^* \,, \]
where the calculus of semiclassical pseudodifferential operators 
is reviewed in \S \ref{rsa},
\begin{gather}
\label{eq:gac}
\begin{gathered}
 P(h) = p^w ( x, hD) + h p_1^w ( x , h D; h ) \,, \ \ 
p_1 \in S^{0,2} ( T^* X) 
\,, \\
 | \xi | \geq C \ \Longrightarrow \ 
p ( x , \xi ) \geq \langle \xi \rangle^2 / C \,, \ \ \
 p = 0 \ \Longrightarrow dp \neq 0 \,,\\
\exists \; R,  \ \ 
u \in \CI ( X \setminus B ( 0 , R) )  \ \Longrightarrow \ P ( h ) u ( x ) = 
Q ( h ) u ( x ) \,, \ \
Q(h) u =\sum_{|\alpha|\leq2} a_{\alpha}(x;h){(hD_x)}^{\alpha} u   \,,
\end{gathered}
\end{gather}
where $a_\alpha(x;h)=a_\alpha(x)$ is independent of $h$ for $|\alpha|=2$,
$a_{\alpha}(x;h)\in C_b^\infty(\RR^n)$ are uniformly bounded with respect to
$h$, here $C_b^\infty(\RR^n)$
denotes the space of $C^\infty$ functions on $\RR^n$ with bounded derivatives 
of all orders,
\begin{gather}
\label{eq:valid}
\begin{gathered}
  \sum_{|\alpha|=2} a_\alpha(x) {\xi}^\alpha\geq {(1/c)}{|\xi|}^2, \;\; 
{\forall\xi\in\RR^n}\,, \text{ for some constant $c>0$, } \\
\text{ $\sum_{|\alpha|\leq2}a_\alpha(x;h){\xi}^\alpha \longrightarrow 
{\xi}^2 -1 $ 
uniformly with respect to $h$ as $|x|\to \infty.$}
\end{gathered}
\end{gather}
We also need the following analyticity
assumption in a neighbourhood of infinity:
there exist $\theta\in {[0,\pi)},$ $\epsilon>0$ and $R\geq{R_0}$ such that the
coefficients $a_\alpha(x;h)$ of $Q(h)$ extend holomorphically in $x$ to
$$\{ r\omega: \omega\in {{\mathbb C}^n}\,, \ \ 
 \text{dist}(\omega,{\bf{S}}^n)<\epsilon\,, \ 
r\in {{\mathbb C}}\,, \   |r|>R\,, \ 
 \text{arg}\,r\in [-\epsilon,\theta_0+\epsilon) \}\,, $$
with \eqref{eq:valid} valid also in this larger set of $x$'s.
We remark that in \cite{HeSj},\cite{SjDuke} the operators were required to be 
globally real analytic but the conditions at infinity were much more
general. A particularly nice feature of the theory developed in \cite{HeSj}
is allowing arbitrary homogeneous polynomials as potentials (see
\cite[(c.31)-(c.33)]{SjDuke}).

\renewcommand\thefootnote{\dag}%

The first theorem we present fits naturally in the methodology of this paper.
It is a slight generalization of a result of
Martinez 
\cite{Mart} which in turn was a $ \CI $ version of a similar
result in \cite{SjDuke}, already implicit in \cite{HeSj}. 
As discussed at the end of \S \ref{rfr} it is essentially
optimal. 
\begin{thm}
\label{t:1}
Suppose that $ X$ is non compact and euclidean outside of a compact set,
and that the operator $ P \in \Psi^{0,2}(  X) $
satisfies the assumptions below. Suppose in addition that no
orbit of $ H_p $ on $ p^{-1} ( 0 ) $ is trapped:
\begin{equation}
\label{eq:th3a}
\forall \; K \Subset p^{-1} ( 0
 ) \,, \ \exists \; T_K  \,, \ \
( x , \xi ) \in K \ \Longrightarrow \ \
\exp ( t H_p ( x , \xi ) ) \notin
K \,, \ t > T_K \,. \end{equation}
Then, for any $ M > 0 $ there exists $ h( M ) $ such that for
$ 0 < h < h(M ) $,
\begin{equation}
\label{eq:th3}
| \Res ( P ) \cap D( 0 , M h \log(1/h) )   | = \emptyset \,.
\end{equation}
Here $ \Res ( P ) $ denotes the set of resonances of $ P $ defined 
in some $ h $-independent neighbourhood of $ 0 $.
\end{thm}
Before stating the main result which requires hyperbolicity 
of the flow on the energy surface we first give the general 
upper bound which generalizes
slightly\footnote{Here we consider an arbitrary manifold in the
compact part. It is clear that the methods of \cite{Bo} easily
allow this type of generalization and the point is that our
method is different and more robust.} the results of Bony \cite{Bo}
which in turn generalized earlier results of Petkov and the second author
\cite{PZ1},\cite{PZ2}.

\begin{thm}
\label{t:2}
Suppose that $ X$ is Euclidean outside a compact set,
and that the operator $ P \in \Psi^{2,0}(  X) $
satisfies the general assumptions below. Then, as $ h \rightarrow 0 $,
\begin{equation}
\label{eq:th2}
| \Res ( P ) \cap D( 0 , Ch )  | =
{\mathcal O} ( h^{-n+1} ) \,.
\end{equation}
\end{thm}

The estimate \eqref{eq:th2} is also optimal in the same way that the
analogous estimate for eigenvalues of a self-adjoint operator with a compact
and smooth energy surface. That follows for instance from applying 
\cite[Corollary, \S 5]{NaStZw}. 

The basic hyperbolicity assumption at an energy $ E $ can be
stated as follows: for $ \rho \in p^{-1} ( E) $ lying in a 
neighbourhoood of the trapped set $ K_E $ we have,
\begin{equation}
\label{eq:hyp}
\begin{split}
& T_\rho (p^{-1} (E)) = \RR H_p (\rho ) \oplus E_+ ( \rho) \oplus E_- ( \rho ) \
\,, \ \ \dim E_\pm ( \rho ) = n -1
\,, \\
& p^{-1} ( E )  \ni \rho \longmapsto E_{\pm} ( \rho ) \subset
T_\rho (p^{-1} (E)) \
\text{ is continuous,} \\
&  d (\exp t H_p )_\rho ( E_\pm ( \rho ) ) = E_\pm ( \exp t H_p ( \rho ) )
\,, \\
& \exists \; \lambda > 0  \ \ \ \| d (\exp t H_p )_\rho ( X ) \| \leq C e^{ \pm \lambda t } \| X \| \,,
\ \text{ for all $ X \in E_\pm ( \rho ) $, $ \mp t \geq 0 $.}
\end{split}
\end{equation}
An example of a potential satisfying this assumption at a range of
non-zero energies is given in Fig.\ref{fig:1} -- see \cite{Mo} and 
\cite[Appendix c]{SjDuke}.
Following \cite{SjDuke} we will formulate a weaker dynamical hypothesis
in \S \ref{efh}.

The main result of this paper is 
\begin{thm}
\label{t:3}
Suppose that $ P( h ) $ satisfies our general assumptions 
\eqref{eq:gac}, 
the flow of $ H_p $ near zero energy is hyperbolic in the sense of
\eqref{eq:hyp}, or the weaker sense given in \S \ref{esf}, and that 
the trapped set at zero energy has 
Minkowski dimension $ m_0 = 2 \nu_0 + 1$. Then for any $ \nu > \nu_0 $,
and $ C_0 > 0 $
there exists $ C_1 $ such that
\begin{equation}
\label{eq:t1}
| \Res ( P ( h ) ) \cap D ( 0 , C_0 h ) | 
\leq  C_1 h^{- \nu } \,. 
\end{equation}
When the trapped is set is of {\em pure dimension}, $ \nu$ can be replaced by $ \nu_0 $.
\end{thm}
A sharp rigorous lower bound is known when $ K $ is an isolated
hyperbolic trajectory. A very precise asymptotic 
description of resonances in that case is given in \cite{GeSj} and it
implies \eqref{eq:t1} with $ \nu = 0 $. In spite of the convincing 
numerical evidence cited above 
no rigorous examples with non-integral values $ \nu_0 $ are
known. A recent indication of the delicate nature of lower bounds for
resonances was given in \cite{TCh} where a class of {\em complex} 
compactly supported potentials in $ \RR^3 $ with {\em no} resonances at all.
We have no reasonable hope of obtaining any analogue of \eqref{eq:1'}
at the present moment.

The methods of this paper apply also to a simpler problem of 
operators with complex absorbing barriers. Let $ V \in \CIc ( B ( 0 ,R_0 ) 
 ; \RR) $ be 
a potential for which $ H_p $, $ p = \xi^2 + V ( x ) - 1 $,  
has hyperbolic flow on $ p = 0 $, for instance a ``three bump'' potential 
\cite{Mo}. Now let $ W \in \CI ( \RR^n ) $ satisfy
\[ W ( x ) \geq 0 \,, \ \ \supp W \cap B ( 0 , R_0 ) = \emptyset \,, \ \
W ( x ) \geq 2 C_0 > 0  \ \text{ for $ | x| > 2 R_0 $.} \]
Consider then (see \cite{Stef} and references given there)
$$ \widetilde P( h )  = - h^2 \Delta + V ( x ) - i  W ( x ) \,. $$ 
The spectrum of this non-selfadjoint operator lies in $ \overline \CC_- $
and we have the exact analogue of \eqref{eq:2}:
\[ | \Spec \widetilde 
P ( h ) \cap D( 0 , C h ) | = {\mathcal O} ( h^{-\nu } ) \,.\]

\medskip

\noindent
{\sc Acknowledgements.} 
The second author would like to thank Jean-Yves Chemin and Jean-Michel 
Bony for useful discussions, the National Science Foundation 
for partial support under the grant DMS-0200732, and \'Ecole Polytechnique
for its generous hospitality in Fall 2004.

\section{Outline of the proof}
\label{oof}

To prove the main result on fractal upper bounds 
(Theorem 3) we first develop methods
for proving the natural results on the absence of resonances 
(Theorem 1) and on general upper bounds at non-degenerate 
energies (Theorem 2). In this section we present the 
general ideas. All of them have origins in other works and
pointers to the literature will be given in corresponding sections. 

The absence of resonances for operators with $ \CI $ coefficients
in domains of size $  h \log ( 1/ h ) $ around an 
energy level hold under a nontrapping condition:
\[ 
\exists \; \epsilon_0 > 0 \,, \
\forall \; K \Subset p^{-1} (0
 )\,,  \ \exists \; T_K  \,, \ \
( x , \xi ) \in K \ \Longrightarrow \ \
\exp ( t H_p ( x , \xi ) ) \notin
K \,, \ t > T_K \,. \]
This implies the existence of an escape function in a neighbourhood of
$ p^{-1} ( 0 ) $:
\[ \exists \; G_1 \in \CI ( T^* X )\,, 
 \  \ H_p G_1 ( x , \xi ) \geq c_0 > 0 \,, 
 \ \text{for $ | p ( x , \xi ) | < \epsilon_0 $.} \]
The resonances of $ P $ are given by the eigenvalues of the deformed
operator $ P_\theta $. In the case of $ P = -h^2 \Delta + V ( x )  $ with 
$ V $ analytic in a conic neighbourhood of $ \RR^n $,  
$$ V ( x ) + 1 \longrightarrow 0, \ \ |x| \longrightarrow \infty \,,  $$ 
the scaled operator is simply 
$$ P_\theta = -h^2 e^{- 2 i \theta} \Delta + V ( e^{ i \theta} x ) \,,  $$
and  it behaves as $ - h^2 e^{ -2 i \theta } \Delta - 1 $ near
infinity. For $ \theta > 0 $ that last operator is clearly invertible. 

We can introduce a modified 
$ G = \chi G_1 $, $ \chi \in \CIc ( X ) $, so that for 
$$ \theta \sim \epsilon \sim M h \log ( 1 / h ) \,, $$
we have 
\[ |\Re p_\theta | < \delta \ \Longrightarrow \ 
- \Im p_\theta +  \epsilon H_{p} G \geq c_0 \epsilon \,, \ \ 
p_\theta = \sigma ( P_\theta ) \,. \]
The operators $ \exp ( \pm \epsilon G^w ( x , hD ) / h ) $ are 
now pseudodifferential operators in a mildly exotic class and we
consider
\[ P_{\theta , \epsilon } = e^{-\epsilon G^w / h } P_\theta 
e^{ \epsilon G^w / h } \,.\]
The spectrum of $ P_\theta $ in $ D ( 0 , M h \log ( 1/h ) ) $ is 
the same as that of $ P_{\theta, \epsilon } $ but the properties
of $ G $ imply that 
$$ \| P_{\theta, \epsilon }^{-1} \| \leq C/ \epsilon $$
showing that in fact there is no spectrum in $ D ( 0 , M' h \log ( 1/h )) $.
This approach allows us to obtain the absence of resonances 
very directly.

In Theorem 2 we show that if $ 0 $ is a non-critical energy level then
\[ | \Res P \cap D ( 0 , C h ) | = {\mathcal O} ( h^{-n + 1 } ) \,.\]
The proof follows from a ``robust'' 
proof of the same estimate for an operator
with a compact resolvent (for instance, an elliptic operator on a compact
manifold). Let  $P $ be such an operator, say, $ P = -h^2 \Delta_g - 1 $, 
on a compact Riemannian manifold. 
We would like to consider a modified operator
\[ \widetilde P ( h ) 
 \stackrel{\rm{def}}{=} 
 P ( h ) - i M h \psi ( M P( h ) / h ) \,, \ \
\psi \in \CIc ( \RR ) \,, \]
whose ``symbol'', $ p - i M h \psi ( M p / h ) $, has 
the absolute value bounded from below by $ M h / 2 $ everywhere. That does not
make sense at first since 
$$ \psi ( M P ( h ) / h ) $$ 
is {\em not}
an $h$-pseudodifferential operator. To remedy this we construct
a second microlocal calculus with a new Planck constant $ \tilde h 
\sim 1/ M $. The new operator $ \widetilde P ( h )$ becomes elliptic
in this calculus and for $ \tilde h $ small enough, independent of
$ h $, it is invertible. We then have 
\[ ( P ( h ) - z )^{-1} =( 
I +  K ( z )  )^{-1}  ( \widetilde P ( h ) - z )^{-1}  \,, \ \ 
K ( z )  \stackrel{\rm{def}}{=} i ( \widetilde P ( h ) - z )^{-1}  M h \psi ( 
 P ( h ) / h ) \]
and the eigenvalues of $ P ( h ) $ near $ 0 $ coincide with the
zeros of $ \det ( I + K ( z ) ) $. The zeros of this determinant are
the same as the zeros of a determinant $ \det (  I + R ( z ) ) $ 
where $ R ( z) $ is a finite rank operator with the rank 
proportional to $ h^{-n} $ times the volume of the support of 
$ \psi ( M p / h ) $. That gives estimates on the deteminant which 
imply \eqref{eq:1}. A slight modification of this argument is
needed to obtain Theorem 2.

\begin{figure}[htbp]
\begin{center}
\includegraphics[width=4.5in]{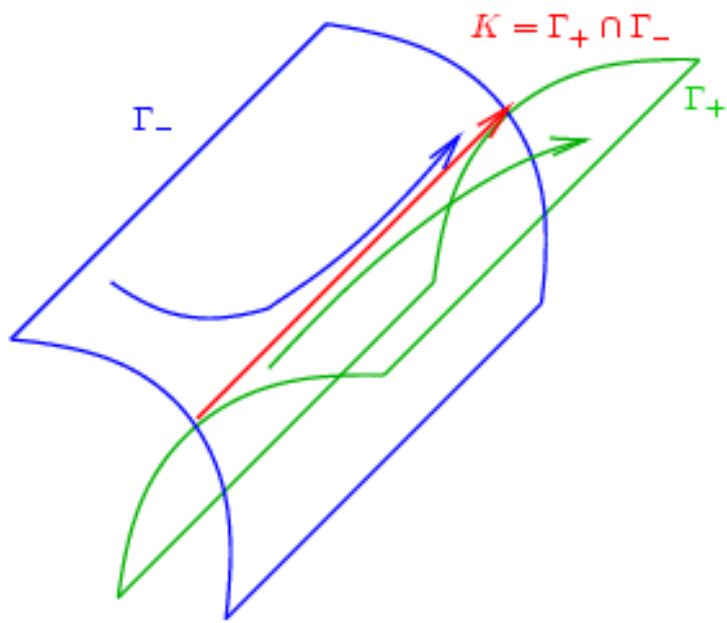}
\end{center}
\caption
{\label{fig:2} Outgoing and incoming sets in the case of one
orbit in a three dimensional energy hypersurface.}
\end{figure}

We now assume that the flow of $ H_p $ is hyperbolic and 
introduce the sets
\begin{equation}
\label{eq:gpm}
  \Gamma_\pm  \stackrel{\rm{def}}{=} \{  ( x , \xi ) \in T^* X \; : 
\; 
 p ( x , \xi ) =0 \,, \ \ \exp ( t H_p ) ( x , \xi ) 
\not \rightarrow \infty \,, \ t \rightarrow \mp \infty \} \,,\end{equation}
depicted in a simple case in Fig.\ref{fig:2}. The trapped set at 
zero energy is 
\begin{equation}
\label{eq:ktr}
 K = \Gamma_+ \cap \Gamma_- \,.\end{equation}
If we assume that $ K \subset \overline {\Gamma_\pm \setminus K } $,
that is $ K $ has no component isolated from infinity, then $ K $ 
is a set of Liouville measure $ 0 $. 

To prove an upper bound involving the dimension of $ K $ we combine
the methods used to prove Theorems 1 and 2. There exist 
functions $ \varphi_\pm \in {\mathcal C}^{1,1} ( T^* X ) $ such that, 
uniformly on compact sets,
\[   H_p \varphi_\pm \sim \mp \varphi_\pm \,, \ \ \varphi_\pm \sim d ( \Gamma_\pm , 
\bullet )^2  \,, \ \ \varphi_+ + \varphi_- \sim d ( K , \bullet )^2 \,. \] 
A local model for the simplest case of one trajectory is given by 
$ p = \xi_1 + x_2 \xi_2 $, $ ( x, \xi ) \in T^* \RR^2 $, 
so that 
\begin{equation}
\label{eq:simo}  H_p = \partial_{x_1} + 
x_2 \partial_{x_2} - \xi_2 \partial_{\xi_2} \,, \ \ 
\varphi_+ = \xi_2^2 \,, \ \ 
\varphi_- = x_2^2 \,, \ \ K = \{ ( t , 0 ; 0 , 0 ) \; : \; t \in \RR \} \,.
\end{equation}
A new escape function is given by 
\begin{equation}
\label{eq:1es} G \stackrel{\rm{def}}{=} 
\left( 
\log (  C \epsilon + \hph_- ) - \log (  C \epsilon + \hph_+ ) \right)  \,,
\ \ \epsilon \sim M h \,,  \  \ M \gg 1 \,,  \end{equation}
where $ \hph_\pm $ are suitable $ h$-dependent regularizations of 
$ \varphi_\pm $.

The logarithmic flattening of the more straightforward escape function 
$ \varphi_- - \varphi_+ $ is forced by the requirement that $ G = {\mathcal O} 
( \log ( 1/ h ) ) $ so that the conjugation used in the proof of
Theorem 1 can be applied. However, even for uniformly smooth 
$ \widehat \varphi_\pm $ the regularization of $ G $ is essentially in the 
symbolic class $ S_{\frac12} $ and the situation becomes more 
complicated in general. Nevertheless we obtain the following estimates:
\begin{gather*}
 \partial_{( x, \xi ) }^\alpha H_p^k G = 
 {\mathcal O} ( \epsilon ^{-\frac{|\alpha|}2 } )  \,, 
\ \ \text{for $ | \alpha |+ k \geq 1 $, uniformly on compact sets,}
\end{gather*}
and 
\begin{gather*}
 d ( (x ,\xi ) , K)^2 
\geq C \epsilon   \ 
\Longrightarrow \ H_p G \geq 1 / C  \,. 
\end{gather*}
As in the proof of Theorem 1 (but with very different parameters and 
escape functions) we introduce a conjugated operator,
\[   P_{\theta , t} ( h )   \stackrel{\rm{def}}{=} 
e^{- t G^w } P_\theta ( h ) e^{ tG^w } \,,
\]
which now is in an exotic $\frac12$-class, with the second Planck constant
$ \tilde h \sim 1/M $ playing the r\^ole of the asymptotic parameter.
The escape function used here, $ G$, has compact support.

\begin{figure}[htbp]
\begin{center}
\includegraphics[width=4.5in]{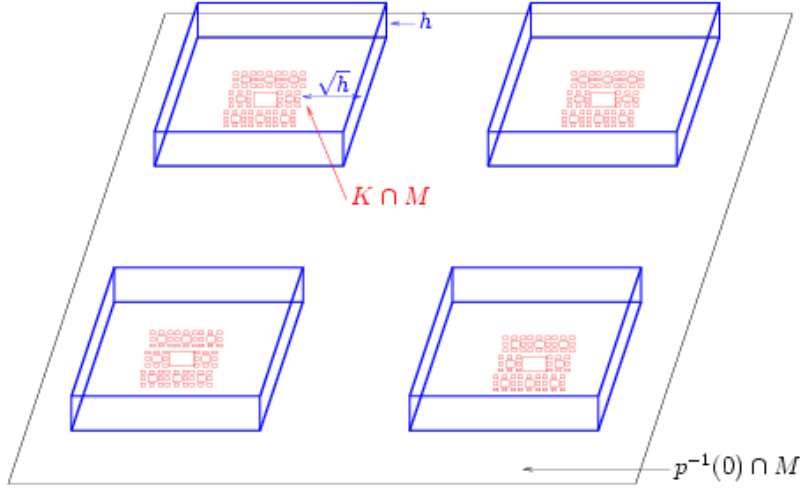}
\end{center}
\caption
{\label{fig:3} The trapped set $ K \subset p^{-1} (0 ) $ intersected with a
hypersurface $ M \subset T^* X $ transversal to the flow and 
surrounded by an $h$-dependent neighbourhood admissible in the 
second microlocal calculus. For $ K $ of  dimension less than 
 $ 2 \nu + 1 $,  
the volume of this neighbourhood is bounded by  $ h^{n - \nu } $, 
$ n = \dim X $.}
\end{figure}

We now build a second microlocal calculus which combines this 
exotic class with the one used in the proof of Theorem 2. The first
allows us the use of the irregular escape function and the second
allows a localization to an $ h $-neighbourhood of the energy 
surface. In the new calculus the operator
\[ \widetilde  P_{\theta,t} =  P_{\theta, t}
 - i M h \widetilde {\Op} ( a ) 
\,, \ \ 
 a ( x , \xi ) \stackrel{\rm{def}}{=}
 \chi \left( \frac{p ( x , \xi ) }{ K_1 h } \right) 
\chi \left(\frac{ C_0 H_p G ( x, \xi ) }{ K_1 h}  \right)  \,, \]
is globally elliptic (here $ \widetilde {\Op } $ describes a
second microlocal  quantization operator). As in the proof of 
Theorem 2 the number of eigenvalues of $ P_{\theta, t}  $, and
hence $ P_\theta$, near $ 0 $, is estimated by $ h^{-n}$ times the 
volume of the support of  $a $. A cross-section of that support 
with a hypersuface transversal to the flow is illustrated in 
Fig.\ref{fig:3}. That volume is bounded by $ h^{1 + ( 2n - 2 - 2 \nu )/2 } $,
where $ \nu > \nu_0 $, $ 2 \nu_0 + 1 $ is the dimension of $ K $. 
That gives \eqref{eq:2}.

\section{Preliminaries}

\subsection{Review of semiclassical pseudodifferential calculus}
\label{rsa}
Let $ X $ be a $\CI$ manifold which is agrees with $ \RR^n $ outside a
compact set, or more generally
\[  X = X_0 \sqcup (\RR^n \setminus B( 0 , R_0 ) ) \sqcup \cdots \sqcup
(\RR^n \setminus B( 0 , R_0 ) )  \,, \ \ X_0 \Subset X \,. \]

We introduce the usual class of semiclassical 
symbols on $ X $:
\[ S^{m,k} ( T^* X ) = \{ a \in \CI( T^* X \times (0, 1]  ) :
|\partial_x ^{ \alpha } \partial _\xi^\beta a ( x, \xi ;h ) | \leq
C_{ \alpha, \beta} h^{ -m } \langle \xi \rangle^{ k - |\beta| } 
\} \,, \]
where outside $ X_0 $ we take the usual $ \RR^n $ coordinates in this 
definition.
The corresponding  class of pseudodifferential operators is denoted by 
$ \Psi_{h}^{m,k} ( X ) $, and we have the quantization and symbol maps:
\[ 
\begin{split}
  & \Op \; : \; S^{m , k } ( T^* X ) \ \longrightarrow 
  \Psi^{m,k}_h ( X) \\
  & \sigma_h \; : \; \Psi_h^{m,k} ( X ) \ \longrightarrow 
  S^{ m , k } ( T^* X ) / S^{ m-1, k-1} ( T^* X ) \,, \end{split}
\]
with both maps  surjective,  and the usual properties 
\begin{gather}
\label{eq:psym}
\begin{gathered}
 \sigma_h ( A \circ B ) = \sigma_h ( A)\sigma _h ( B ) \,,\\
 0 \rightarrow \Psi^{m-1, k-1} ( X) \hookrightarrow \Psi^{m, k} ( X)
\stackrel{\sigma_h}{\rightarrow} S^{ m , k } ( T^* X ) 
/ S^{ m-1, k-1} ( T^* X ) \rightarrow 0 \,, 
\end{gathered}
\end{gather}
a short exact sequence, and
\[ \sigma_h \circ \Op : S^{m,k} ( T^* X ) 
\ \longrightarrow  S^{ m , k } ( T^* X ) / S^{ m-1, k-1} ( T^* X ) \,,
\]
the natural projection map.
The class of operators and the quantization map are defined locally using the
definition on $ \RR^n $:
\begin{equation}
  \label{eq:weyl}
  \Op (a) u ( x) = \frac{1}{ ( 2 \pi h )^n } 
  \int \int  a \left( \frac{x + y }{2}  , \xi \right
  ) e^{ i \langle x -  y, \xi \rangle / h } u ( y ) 
dy d \xi \,, \end{equation}
and we refer to \cite{DiSj} for a detailed discussion.
  We remark only that when we consider the operators
  acting on half-densities we can define 
  the symbol map, $ \sigma_h $, onto
  $$  S^{ m , k } ( T^* X ) / S^{ m-2, k-2} ( T^* X ) \,, $$ 
   see
  \cite[Appendix]{SjZw02}. We keep this in mind but 
  for notational simplicity we suppress the half-density notation. 

  For $ a \in S^{m,k} ( T^* X ) $ we define 
  \[ {\text{ess-supp}}_h\; a \subset T^* X \sqcup S^* X \,, \ \ 
  S^* X \stackrel{{\rm{def}}}{=}  ( T^* X  \setminus 0 ) / {\mathbb R}_+ \,, 
  \]
  where the usual $ \RR_+ $ action is given by multiplication on 
  the fibers: $ ( x, \xi ) \mapsto  ( x , t \xi ) $, as 
   \[  \begin{split}  {\text{ess-supp}}_h\; a = &   
    \;  \complement \{ ( x, \xi ) \in T^* X \; : \; 
    \exists \; \epsilon > 0 \,, \ \partial_x ^\alpha \partial_\xi ^\beta 
    a ( x', \xi' ) = {\mathcal O} ( h^\infty ) 
    \,, \  d( x, x' ) + | \xi - \xi' |  < \epsilon \}  \\
    & \   \sqcup \;    \complement \{ 
    ( x, \xi ) \in T^* X \setminus 0  :  
    \exists \; \epsilon > 0\,,  \ \partial_x ^\alpha \partial_\xi ^\beta 
    a ( x', \xi' ) = {\mathcal O} ( h^\infty  \langle \xi'
     \rangle^{  -\infty}) \,, \\
&    \ \ \ \ \   
    d( x, x' ) + 1 / |\xi'| + | \xi/ |\xi|  - \xi'/ |\xi'| |  
    < \epsilon \} / {\mathbb R}_+ \,,
  \end{split} \]
where the second complement is in $ S^*X $.
  For $ A \in \Psi^{m,k} _h ( X) $, then define
  \[ \WFh ( A) =  {\text{ess-supp}}_h\; a \,, \ \ A = \Op ( a ) \,, 
  \]
  noting that, as usual, the definition does not depend on the choice of
  $ \Op $, and 
\[ {\rm{Char}} ( A ) = \bigcup \left\{ WF_h ( B ) \; : \; 
B \in \Psi^{m,k}_h ( X ) \,, \ \sigma_h ( B ) = \sigma_h ( A ) \right\} \,,\]
where $ \sigma_h $ is the principal symbol map in \eqref{eq:psym}.
 For 
  \[  u \in \CI ( ( 0 , 1]_h ; \CI ( X) ) \,, \ \  \forall \; 
K \Subset X \,, N \in \NN 
\ \exists \; P\,, h_0 \,, \   \   \| u \|_{ C^N ( K ) } \leq h^{-P} \,, \ 
h < h_0 \,. \]
we 
define
\[
\WFh ( u ) =  \left( \bigcup \left\{ {\rm{Char}} ( A ) \; : \; 
 A \in \Psi^{0,0} ( X ) \; : \; 
\ A u \in h^\infty  \CI ( ( 0 , 1]_h ; \CI ( X) ) \right\}  
\right)^\complement
 \,,\]
where the complement is taken in $ T^* X \sqcup S^* X $.
Here we will be concerned with a purely semiclassical 
theory and deal only with {\em compact} subsets of $ T^* X $.

To illustrate the $h$-pseudodifferential calculus at work 
we prove two simple lemmas which will be used later. We say that
$ A \in \Psi^{m,k} ( X ) $ is elliptic on $ K \Subset T^*X $ if
$ |\sigma ( A ) \rest_K | > h^{-m}/C $. This is equivalent to 
saying 
\begin{lem}
  \label{l:r1}
  Suppose $ Q \in \Psi^{0,m} ( X) $ is elliptic at $ ( x_0 , \xi_0 ) $, 
$ \| u \|_{L^2 } = 1 $, and $ \WFh ( u ) $ 
  is contained in a sufficiently small neighbourhood of $ ( x_0 , \xi_0 )$.
  Then for $ h $ small enough,
  \[ \| Q u \|_{L^2} \geq 1/C  \,.\]
\end{lem}

\begin{lem}
\label{l:r2}
Suppose that $ \psi_j \in \CI_{\rm b} ( T^* X) $, $ \psi_1^2 + \psi^2_2 = 1 $,
$ \supp \psi_1 \subset \{ ( x, \xi ) \; : \; |\xi | \leq C \} $. 
Then, there exist $ \Psi_1 \in \Psi^{0, - \infty} ( X )  $ and 
$ \Psi_2  \in \Psi^{0,0} ( X ) $,
 with principal symbols $ \psi_1 $ and $ \psi_2 $ respectively, such that 
\[ \Psi_1^2 + \Psi_1^2 = I + R \,, \ \  R \in \Psi^{-\infty, - \infty } ( X ) 
\,, \ \ \Psi_j^* = \Psi_j 
\,.\]
\end{lem}
\begin{proof}
Functional calculus 
gives
$$ (\psi_1^w) ^2 + ( \psi_2^w )^2 = I + r_1^w \,, \ \ r_1 \in 
S^{-1,-\infty} ( T^* X ) \,, $$
in particular
$ r = {\mathcal O}  ( h ) : H^{-M} ( X )  \rightarrow H^M ( X ) $. 
If $ h $ is small enough we put 
$$ \Psi_j^{1} = ( 1+ r_1^w)^{-\frac14} \psi_j^w ( 1+ r_1^w)^{-\frac14} \,, $$
so that 
\[  (\Psi_1^1)^2 + (\Psi_2^1)
^2 = I + r_2^w \,, \ \  r_2 \in S^{-2,-\infty} ( T ^* X ) \,, \ \ 
(\Psi_j^1)^* = \Psi_j^1\,.  \]
and we can then proceed by iteration.
\end{proof}

The semiclassical Sobolev spaces, $ H_h^s ( X ) $  are defined by 
choosing a globally elliptic, self-adjoint operator, 
$ A \in \Psi^{0,1} ( X ) $ (that is an operator satisfying 
$ \sigma ( A) \geq \langle \xi \rangle / C $ everywhere) and
putting 
\[  \| u \|_{ H_h^s } = \| A^{s} u \|_{L^2 ( X ) } \,.\]
When $  X = \RR^n $, 
\[  \| u \|^2_{ H_h^s } \sim  \int_{\RR^n} 
\langle h \xi \rangle^{2s} |{\mathcal F} u ( \xi ) |^2 d \xi \,, \ \ 
{\mathcal F} u ( \xi ) \stackrel{\rm{def}}{=} 
\int_{\RR^n } u ( x ) e^{ - i \langle x , \xi \rangle } d x \,. \]

The following lemma will also be useful:
\begin{lem}
\label{l:r3}
Suppose that $ P_t $, $ t \in (0, \infty ) $, 
is a family of operators such that 
\begin{gather*}
 P_t \; : \; H^s_h ( X ) \; \longrightarrow \; H^{s-m} _h ( X) \,, \\  
 \forall \; A \in \Psi^{0,-\infty} ( X ) \,,  
 \ \  \ad_{P_t} A = {\mathcal O} ( h ) \; : \; L^2 ( X) \; \longrightarrow 
\; L^2 ( X ) \,, \  \ 0 < h < h_0 ( t ) ,.\end{gather*}
Let $ \Psi_j $ be as in Lemma \ref{l:r2} and suppose that
\[ \| P_t \Psi_j u \| \geq t h \|  \Psi_j u \| - {\mathcal O} ( h / t  ) 
\| u \|  \,, \ \ j = 1,2\,,  \ \ u \in 
\CIc ( X) \,. \]
Here the constants in $ {\mathcal O}$ are independent of $ h $ and $ t $.
Then for $ t > t_0 \gg 1 $ and $ 0 < h < h_0 ( t) $,
\[  \| P_t u \| \geq t h \| u \| /2 \,. \]
\end{lem}
\begin{proof}
We recall from Lemma \ref{l:r2} that
\begin{equation}
\label{eq:pss}
 \| \Psi_1 v \|^2 + \|\Psi_2 v \|^2  = \|v\|^2  + 
\langle  R v ,  v \rangle = \| v\|^2 + {\mathcal O} ( h^\infty ) 
\| v \|_{ H_h^{-N} }  \,,\end{equation}
and hence with $ v = P_t u $,
\vspace{0.1in}
\begin{equation*}
\begin{split}
 \| P_t u \|^2 &  = \| \Psi_1 P_t u \|^2 + \| \Psi_2 P_t u \|^2 - {\mathcal O}(h^\infty)
\| u \|^2 \\
&  \geq  \| P_t \Psi_1  u \|^2 + \| P_t \Psi_2  u \|^2 
-  \| [ \Psi_2, P_t ] u \|^2  -  \| [ \Psi_2, P_t ] u \|^2 \\
& \hspace{0.2in}  - \; 2  \left( \| \Psi_1 P_t u \|\| [ \Psi_1, P_t ] u \|^2 
+ \| \Psi_2 P_t u \|
\| [ \Psi_2, P_t ] u \|^2 
\right)^{\frac12}
- {\mathcal O}(h^\infty) \| u \|^2 \\
&  \geq  \| P_t \Psi_1  u \|^2 + \| P_t \Psi_2  u \|^2 \\
& \hspace{0.2in}  - 2 C ( \| [ \Psi_1, P_t ] u \|^2  +
 \| [ \Psi_2, P_t ] u \|^2 )
-  \| P_{t} u \|^2/C - {\mathcal O}(h^\infty) \| u \|^2  \\
&  \geq  \| P_{t} \Psi_1  u \|^2 + \| P_{ t} \Psi_2  u \|^2 
-   C' h^2 \|  u \|^2 
-  \| P_t  u \|^2/C \,.
\end{split} \end{equation*}
We now use the hypothesis of the lemma and \eqref{eq:pss} with $ v = u $ 
to obtain
\[ \begin{split}
\| P_t u \|^2 & \geq  t^2 h^2 (\| \Psi_1  u \|^2 + \|  \Psi_2  u \|^2)
-   C' h^2 \|  u \|^2
-  \| P_t  u \|^2/C \\
& \geq  t^2 h^2 \|  u \|^2 -   C' h^2 \|  u \|^2
-  \| P_t  u \|^2/C  \end{split} \]
and the lemma follows.
\end{proof}

\subsection{Semiclassical Fourier integral operators.}
\label{fio}

We now follow \cite{SjZw02} and review some 
aspects of the theory of semiclassical
Fourier Integral Operators. We take a point of view which will
be used in showing invariance of the second microlocal calculus
developed below in \S \ref{smc}.

Thus let $ A ( t ) $ be a smooth 
family of pseudodifferential
operators, 
$$ A ( t ) = \Op ( a ( t ) ) \,, 
\ \  a (t)  \in \CI ([-1  , 1]_t ; S^{ 0, -\infty }
( T^* X ))  \,, $$ 
 such
that for all $ t $, $ \WFh ( A ( t ) ) \Subset T^* X $. We then define
a family of operators
\begin{gather}
\label{eq:2.2}
\begin{gathered}
  U( t ) \; : \; L^2 ( X) \rightarrow L^2 ( X) \,, \\
h D_t U ( t) +  U ( t) A(t) = 0 \,, \ \ U( 0 ) = U_0 \in \Psi_h^{0,0} ( X) \,.
\end{gathered}
\end{gather}
This is an example of a family of $h$-Fourier Integral Operators, $ U ( t)$,
associated to canonical transformations $ \kappa(t)$, 
generated by the Hamilton 
vector fields $ H_{a_0(t)} $, where the real valued
$ a_0 (t) $ is the $h$-principal 
symbol of $ A ( t ) $,
\[ \frac{d}{dt} \kappa( t) ( x, \xi ) =  ( \kappa ( t ))_* ( H_{a_0 ( t)} 
( x, \xi )) \,, \ \ \kappa( 0 ) ( x, \xi ) = ( x, \xi )\,,   \ \ 
( x, \xi ) \in T^* X \,.\]
We will often need the {\em Egorov theorem} which 
can be proved directly from this definition: 
when $ U_0 $ in \eqref{eq:2.2} is elliptic (that is 
$ | \sigma (U_0 )| > c > 0 $) on $ T^* X $, then 
for $ B \in \Psi_h^{m,k} ( X ) $
\begin{gather}
\label{eq:egorov}
\begin{gathered}
\sigma ( V(t) B U ( t) ) = ( \kappa ( t) )^* \sigma ( B ) \,, 
\end{gathered}
\end{gather}
where $ V ( t) $ is an approximate inverse to $ U ( t) $, 
$$  V(t) U ( t ) - I \,, \,U ( t) V ( t) - I 
\in \Psi_h^{- \infty , - \infty } ( T^* X)  \,. $$
The approximate inverse is constructed by taking
\[ h D_t V ( t) - A ( t) V(t) = 0 \,,  \ \ V ( 0 ) = V_0 \,, \ \ 
V_0 U_0 - I \,, \, U_0 V_0 - I \in \Psi_h^{- \infty , - \infty } ( T^* X) 
\,, \]
the existence of $ V_0 $ being guaranteed by the ellipticity of $ U_0 $.
The proof of \eqref{eq:egorov} follows from writing $ B(t) = 
V ( t) B U ( t) $, so that, in view of the properties of $ V (t)$, 
\[ h D_t B( t) \equiv [ A( t) , B( t ) ] \, \mod \Psi_h^{-\infty, -\infty}
\,, \ \ B( 0 ) = B_0 \,. \]
Since the symbol of the commutator is given by $ (h/i) H_{ a_0 ( t) }
\sigma ( B ( t) )$, \eqref{eq:egorov} follows directly from the definition
of $ \kappa( t)$.

If $ U = U(1) $, say, and the graph of $ \kappa ( 1) $ is denoted by $C$,
we conform to the usual notation and write
\[  U \in I^0_h ( X \times X ; C' ) \,, \ \ 
C' = \{ ( x, \xi;  y , - \eta) \; : \; ( x, \xi) = \kappa ( y , \eta )   \}
\,, \]
which means that $ U $ is an $h$-Fourier Integral Operator associated to 
the canonical graph $ C$. Locally all $h$-Fourier Integral Operators
associated to canonical graphs are of the form $ U ( 1) $ 
since each local canonical transformation with a fixed point can be 
deformed to identity, see \cite[Lemma 3.2]{SjZw02} and the 
proof of Lemma \ref{l:3.2s} below. 

Our definitions of pseudo-differential operators and of (the special class of)
$h$-Fourier Integral Operators were global. It is useful and natural to 
consider the operators and their properties microlocally.
We consider classes of {\em tempered} operators:
\[ T  \; : \; \CI ( X ) \rightarrow \CI ( X) \,, \]
and for any semi-norm $ \| \bullet\|_1 $ 
on $ \CI ( X) $ there exist a seminorm 
$ \| \bullet\|_2 $ on $ \CI ( X) $ and a constant $ M_0$ such that
\[  \| T u \|_1  = {\mathcal O } ( h ^{ - M_0 } ) 
\| u \|_2 \,. \]
We remark that since we deal with compact subsets of $ T^* X $ here,
we could consider operators $ T \; : \; {\mathcal D}' ( X ) 
\rightarrow \CI ( X ) $ in which case we can ask for existence 
of $ M_0 $ for any {\em two} seminorms $ \| \bullet \|_j $, $ j = 1, 2$.

For open sets, $ V \subset T^* X$, $ U \subset T^* X $, the operators 
{\em defined microlocally} near  $ V \times U $ 
are given by equivalence classes of tempered operators
given by the relation
\[ T \sim T' \ \Longleftrightarrow \ A ( T - T' ) B = {\mathcal O} ( h^\infty )
\; : \; {\mathcal D}' ( X) \ \longrightarrow \CI ( X ) \,, \]
for any 
$ A, B \in \Psi_h ^{0,0 } ( X ) $ such that 
\begin{gather}
\label{eq:2.7}
\begin{gathered}
\WFh ( A ) \subset \widetilde V \,, \ \ \WFh ( B ) \subset \widetilde U\,,
\\
\overline  V \Subset  \widetilde V \Subset T^* X \,, \ \ 
\overline U 
\Subset  \widetilde U \Subset T^* X \,, \ \ \widetilde U \,,\,
\widetilde V \ \text{ \  open}\,.
\end{gathered}
\end{gather}
The equivalence class $ T $, 
 $h$-Fourier Integral Operator associated to a local canonical graph
$ C $ if, again for any $ A $ and $ B$ above 
\[ A T B 
\in I^0 ( X \times X ; \widetilde C' ) \,,\]
where $ C$ needs to be defined only near $ U \times V $. 

We say that $ P = Q $ {\em microlocally} near $ U \times V $ if 
$  A P B - A Q B = {\mathcal O}_{L^2 \rightarrow L^2 }
 ( h ^ \infty ) $, where because of the assumed pre-compactness of 
$ U $ and $ V$ the $ L^2 $ norms can be replaced by any other norms.
For operator identities this will 
be the meaning of equality of operators in this paper, with $ U , V $
specified (or clear from the context). Similarly, we say that $ 
B = T^{-1} $ microlocally near $ U \times V $, if $ B T = I$ microlocally 
near $ U \times U $, and $ T B = I $ microlocally near $ V \times V $.
More generally, we could say that $ P = Q $ microlocally on $ 
W \subset T^* X \times T^* X $ (or, say, $ P $ is 
microlocally defined there), if for any $ U, V$, $ U \times V \subset
W $, $ P = Q $ microlocally in $ U \times V$. 

If the open sets 
$ U $ or $ V $ in \eqref{eq:2.7}
are small enough, so that they can be identified with 
neighbourhoods of points in $ T^* \RR^n $, we can use that identification to 
state that $ T $ is microlocally defined near, say, $ ( m , ( 0 , 0 ) )$, 
$ m \in T^* X $, $ ( 0, 0 ) \in T^* \RR^n $.

To give a useful example of this formalism we state
a semiclassical version of Egorov's theorem.
\begin{prop}
\label{p:egorov}
Suppose that $ F $ is an $h$-Fourier integral operator, microlocally 
defined near $ U \times V \subset T^* X \times T^* X$, and associated
to a locally defined canonical transformation $ \kappa : V \rightarrow U $,
and elliptic near $ ( \kappa( \rho ) , \rho ) \in U \times V $. 
Let $ F^{-1} $ be the microlocal inverse of $ F $ near $ ( \kappa( \rho) ,
\rho )$. 
Then for any $ A \in \Psi^{m,k} ( X ) $,
\begin{equation}
\label{eq:egorov'}
F^{-1} \circ A \circ F = B \in \Psi^{m,k} ( X ) \,, \ \ \sigma_h (B ) 
= \sigma_h ( A ) \circ \kappa \,, \end{equation}
microlocally near $ \rho \in T^* X  $.
\end{prop}
The proof is a localized adaptation of the argument giving 
\eqref{eq:egorov} and the observation recalled above that any canonical 
transformation with a fixed point (as we can assume that 
$ \kappa ( m) = m $) can be deformed to the identity.

\subsection{$S_{\frac12}$ spaces with two parameters}
\label{sma}
We define the following symbol class:
\begin{equation}
\label{eq:sthf}
a \in S^{m , {\widetilde m}, k }_{\frac12} ( T^* \RR^n ) 
\ \Longleftrightarrow \ |\partial_x^\alpha \partial^\beta_\xi 
a (  x, \xi ) | \leq C_{\alpha \beta} h^{-m } \tilde h^{-{\widetilde m}} 
\left( \frac{\tilde h }{ h } 
\right)^{\frac12 ( |\alpha | + |\beta| ) } \langle \xi \rangle^{k-|\beta|}
 \,,\end{equation}
where in the notation we suppress the dependence of $ a $ on $ h 
$ and $ \tilde h $.
We define the Weyl quantization of $ a $ in the usual way
\[ a^w ( x , h D_ x ) u = \frac{1}{ (2 \pi h)^n } 
\int a \left( \frac{x + y }2 , \xi 
\right) e^{ \frac{i}h \langle x -y , \xi \rangle
} u ( y ) dy d \xi \,,\]
and the standard results (see \cite{DiSj}) show that if  
$ a \in   S^{m , {\widetilde m},k }_{\frac12} ( T^* \RR^n ) $ and
$ b \in   S^{m', \widetilde m', k'}_{\frac12} ( T^* \RR^n ) $ then 
\[ a ( x , h D_ x ) \circ b ( x , h D_ x ) = c ( x , h D_ x  ) 
\ \text{ with } \ 
c \in   S^{m +m', {\widetilde m}+ {\widetilde m}' , 
k+ k' }_{\frac12} ( T^* \RR^n ) \,.
\]
The presence of the additional parameter $ \tilde h $ allows us to 
conclude that 
\[ c \equiv \sum_{ |\alpha | <  M } \frac{1}{\alpha !} 
\partial_\xi^\alpha a D_x^\alpha b \ \mod S^{m  + m ' , {\widetilde m} + {\widetilde m}' 
- M , k + k' - M }_{\frac12} ( T^* \RR^n ) \,, \]
that is, we have a symbolic expansion in powers of $ \tilde h $. 
We could also consider an expansion in the Weyl quantization -- see 
\eqref{eq:weylc}.

We denote our class of operators by $ \Psi_{\frac12}^{m , {\widetilde m}, k } 
(T^* \RR^n ) $. 
For simplicity we will only state the characterization \`a la Beals
for a simpler class of symbols:
\begin{lem}
\label{l:1hf}
Suppose that $ A \; : \; {\mathcal S} ( \RR^n ) \rightarrow 
{\mathcal S}' ( \RR^n ) $.
Then $ A = \Op ( a ) $ with 
\begin{equation}
\label{eq:symbh}
 \partial_x ^\alpha \partial_\xi^\beta a = 
{\mathcal O}(h^{-m} \tilde h^{-\widetilde m } )  \left( \frac{\tilde h }{ h } 
\right)^{\frac12 ( |\alpha | + |\beta| ) } \,, \end{equation}
 if and only if 
for 
any sequence $ \{ \ell_j \}_{j=1}^N $ of linear functions on 
$ T^* \RR^n $ we have 
\[  
\| \ad_{\Op (\ell_1)} \circ \cdots \ad_{\Op (\ell_
N ) } A u \|_{L^2 ( \RR^n ) }
\leq C h^{-m  + N/2} \tilde h^{-{\widetilde m} + N/2} \| u \|_{ L^2 ( \RR^n )} \,, \]
for any $ u \in {\mathcal S} ( \RR^n ) $.
\end{lem}
\begin{proof} 
We can assume that $ m  = {\widetilde m} = 0 $. 
The statement follows from the proof in \cite[Chapter 8]{DiSj} and a
rescaling:
\[  ( \tilde x , \tilde \xi ) = ( \tilde h / h )^{\frac12} ( x , \xi ) \,.\]
In fact, we define the following unitary operator on $ L^2 ( \RR^n ) $:
\[ U_{ h , \tilde h } u ( \tilde x ) = ( \tilde h / h )^{\frac{n}{4}}
u ( ( h / \tilde h )^{\frac12}  \tilde x ) \,,\]
for which we can check that
\[ a ( x , h D_ x ) = U_{ h , \tilde h }^{-1} a_{ h , \tilde h } 
( \tilde x , \tilde h D_{\tilde x } ) U_{ h , \tilde h} \,, \ \ 
a_{ h , \tilde h } ( \tilde x , \tilde \xi ) = a ( 
( h / \tilde h)^{\frac12} ( \tilde x , \tilde \xi ) ) \,. \]
Clearly $ a $ satisfies \eqref{eq:symbh} if and 
only if
$ a_{ h, \tilde h } \in \CI_{\rm{b}} ( T^* \RR^n ) $. 
The Beals condition for $ \tilde h $-pseudodifferential operators
is 
\[ \| \ad_{ \tilde 
\ell_1 ( \tilde x , \tilde h D_{ \tilde x } )} 
\circ \cdots \circ 
 \ad_{  
\tilde \ell_N ( \tilde x , \tilde h D_{ \tilde x } )} a_{ h , \tilde h 
} ( \tilde x , \tilde h D_ {\tilde x } ) u \|_{ L^2 } \leq 
C \tilde h ^N \| u \|_{ L^2 } \,. \]
But this is  the condition in the lemma since we should take
\[  \tilde \ell_j = ( \ell_j )_{ h , \tilde h } = ( \tilde h / h )^{\frac12} 
\ell_j \,, \]
and this completes the proof.
\end{proof}
We remark that the proof given in  \cite{HeSj1},\cite[Chapter 8]{DiSj} is 
recalled, in a more complicated setting, 
in the proof of Proposition \ref{p:Beals} below.
It is based on showing that 
\[
 \forall \alpha \in \NN^{2n}\,, \ \| (\partial^\alpha a ) ( x, D ) u \|_{ L^2}
\leq C_\alpha \| u \|_{L^2 } 
\; \Longleftrightarrow \; \forall \beta \in \NN^{2n} \,, \ 
\sup |\partial^\beta a  | \leq C'_\beta \,,
\]
which in the setting presented in Lemma \ref{l:1hf} becomes
\begin{gather}
\label{eq:bealss}
\begin{gathered}
  \forall \alpha \in \NN^{2n}\,, \ \| (\partial^\alpha a ) ( x, h 
 D ) u \|_{ L^2} \leq C_\alpha ( \tilde h / h )^{|\alpha|/2}\| u \|_{L^2 } 
\\ \; \Longleftrightarrow \;\\
 \forall \beta \in \NN^{2n} \,, \ 
\sup |\partial^\beta a  | \leq C'_\beta 
 ( \tilde h / h )^{|\alpha|/2} 
\,.
\end{gathered}
\end{gather}
We will also need the following application of the semi-classical
calculus:
\begin{lem}
\label{l:Sjnew}
Suppose that 
$ \partial^\alpha a $, $ \partial^\alpha b = {\mathcal O}_{\alpha}
( ( \tilde h/ h )^{|\alpha|/2} ) \,, $
and that $ c^w ( x , hD ) = a^w ( x, hD) \circ b^w ( x, hD) $. 
Then 
\begin{equation}
\label{eq:weylc}  c ( x, \xi) = \sum_{k=0}^N \frac{1}{k!} \left( 
\frac{i h}{2} \sigma ( D_x , D_\xi; D_y , D_\eta) \right)^k a ( x , \xi) 
b( y , \eta) \rest_{ x = y , \xi = \eta} + e_N ( x, \xi ) \,,
\end{equation}
where for some $ M $
\begin{equation}
\label{eq:Sjnew1}
\begin{split}
& | \partial^{\alpha} e_N | \leq C_N h^{N+1}
\times \\
& \ \ 
\sum_{ \alpha_1 + \alpha_2 = \alpha} 
 \sup_{ 
{{( x, \xi) \in T^* \RR^n }
\atop
{ ( y , \eta) \in T^* \RR^n }}} \sup_{
|\beta | \leq M  \,, \beta\in \NN^{2n} }
\left|
( h^{\frac12} 
\partial_{( x, \xi; y , \eta) } )^{\beta } 
(i  \sigma ( D) / 2) ^{N+1} \partial^{\alpha_1} a ( x , \xi)  
\partial^{\alpha_2} b ( y, \eta ) 
\right| \,,
\end{split} 
\end{equation}
where $ \sigma ( D) = 
 \sigma ( D_x , D_\xi; D_y, D_\eta ) \,. $ 
\end{lem}
\begin{proof}
This follows from from the standard estimates of symbolic 
calculus (see \cite[Proposition 7.6]{DiSj}): suppose that
$ A ( D ) $ is a non-degenerate real quadratic form. Then
there exists $ M $ such that 
\[ | \partial^\alpha \exp( i A ( D ) ) a ( x , \xi ) | 
\leq C \sum_{ |\beta | \leq M } \sup_{( x , \xi) \in T^* \RR^n } 
|\partial^{\alpha + \beta } a ( x , \xi ) | \,. \]
We observe that a rescaling $ \tilde x = x / \sqrt s $, $ 
s > 0 $,  shows that
\[  | \partial^\alpha \exp( i s A ( D ) ) a ( x , \xi ) | 
\leq C \sum_{ |\beta | \leq M } \sup_{( x , \xi) \in T^* \RR^n } 
|\partial^{\alpha } ( \sqrt s \partial) ^\beta  a ( x , \xi ) | \,.
\] 
To obtain an expansion we 
use the Taylor expansion: 
\[ \exp ( i h A ( D) ) = \sum_{k=0}^N \frac{ ( i h A ( D) )^k}{k!} 
+ \frac{1}{ N!} \int_0^1 ( 1 - t)^N \exp ( i t h A ( D) ) 
( i h A( D) )^{N+1} dt \,. \]
In the notation of the lemma and with 
 $ A ( D ) = \sigma ( D_x, D_\xi; D_y , D_\eta ) /2 $, 
\[ c ( x , \xi ) = \exp ( i A ( D ) ) a ( x , \xi ) b ( y , \eta ) 
\rest_{ x = y ,  \eta = \xi } \,, \]
and the lemma follows.
\end{proof}
As a particular consequence we notice that if $ a \in 
S^{0,0,-\infty}_{\frac12} ( T^* \RR^n ) $ and $ b \in S^{0,-\infty} ( 
T^* \RR^n ) $ then 
\begin{gather*}
a^w ( x , h D) \circ b^w ( x , h D) = c^w ( x, hD ) \,, \ \ 
  c ( x, \xi) = \\
 \sum_{k=0}^N \frac{1}{k!} \left( 
i h \sigma ( D_x , D_\xi; D_y , D_\eta) \right)^k a ( x , \xi) 
b( y , \eta) \rest_{ x = y , \xi = \eta} + {\mathcal O} 
( h^{\frac{N+1}2 } \tilde h^{\frac{N+1}2} ) 
\,, 
\end{gather*}
and the usual pseudodifferential calculus allows a remainder 
improvement to  
$$  {\mathcal O} 
( h^{\frac{N+1}2 } \tilde h^{\frac{N+1}2}
 \langle \xi \rangle^{-\infty} ) \,.$$

\subsection{One parameter groups of elliptic operators.}
We recall a special case of a result of
Bony and Chemin \cite[Th\'eoreme
6.4]{BoCh}. Let $ m ( x , \xi ) $ be an order function in 
the sense of \cite{DiSj}:
\begin{equation}
\label{eq:orderm}
  m ( x , \xi ) \leq C m ( y , \eta ) \langle  ( x-y , \xi - \eta ) 
\rangle ^N \,. 
\end{equation}
The class of symbols, $ S ( m ) $,  corresponding to $ m $ is defined as
\[  a \in S ( m ) \ \Longleftrightarrow \ |\partial^\alpha_x \partial_\xi^\beta
a ( x , \xi ) | \leq C_{\alpha \beta} m ( x , \xi ) \,. \]
If $ m_1$ and $ m_2 $ are order functions in the sense of \eqref{eq:orderm},
and $ a_j \in S (m_j ) $ then (we put $ h=1 $ here),
\[ a_1^w ( x, D ) a_2^w ( x , D ) = b^w ( x , D) \,, \ \ b \in S ( m_1 m_2 ) \,,\]
with $ b $ given by the usual formula,
\begin{equation}
\label{eq:usual} \begin{split} b ( x , \xi ) & = a_1 \; \# \; a_2 ( x, \xi ) \\
& \stackrel{\rm{def}}{=}
\exp ( i \sigma ( D_{x^1}, D_{\xi^1} ; D_{x^2 } , D_{\xi^2} )/2  ) 
a_1 ( x^1 , \xi^1 ) a_2 ( x^2 , \xi^2 ) \rest_{ x^1 = x^2 = x , 
\xi^1 = \xi^2 = \xi } \,. \end{split} \end{equation}
A special case of \cite[Th\'eoreme 6.4]{BoCh} gives 
\begin{prop}
\label{p:bbc}
Let $m $ be an order function in the sense of \eqref{eq:orderm} and
suppose that $ G \in \CIc ( T^* \RR^n; \RR  )  $ satisfies 
\begin{equation}
\label{eq:bc1}
 G ( x , \xi ) - \log m ( x , \xi )  = {\mathcal O} ( 1 ) \,, \ \
\partial_x^\alpha \partial_\xi^\beta G ( x, \xi ) = {\mathcal O} ( 1 ) \,,
 \ \ |\alpha | + |\beta| \geq 1 \,.
\end{equation}
Then 
\begin{equation}
\label{eq:bc2} 
\exp ( t G^w ( x , D ) ) = B_t^w ( x , D ) \,, \ \ B_t \in S ( m^t ) \,.
\end{equation}
Here $ \exp ( t G^w ( x , D) ) $ is constructed 
using spectral theory of bounded
self-adjoint operators. The estimates on $ B_t \in S( m^t ) $ depend 
{\em only} on the constants in \eqref{eq:bc1} and in \eqref{eq:orderm}.
In particular they are independent of the support of $ G $.
\end{prop}
In Appendix at the end of the paper we give a simple direct proof 
of this proposition. We should stress that the main difficulties in
\cite{BoCh} came from considering general Weyl calculi of pseudodifferential
opearators. Here we need only the case of the simplest metric $
g = dx^2 + d\xi^2 $.

\subsection{Review of complex scaling}
\label{rcs}

We very briefly
recall the procedure described in \cite{Sj}. 
It follows the long tradition of the complex scaling method -- 
see \cite{SjZw1} for the presentation for compactly supported
perturbations and references to earlier work.

Let $ \Gamma_\theta
\subset \CC^n $ be a totally real contour with the following properties:
\begin{gather}
\label{eq:gpr}
\begin{gathered}
\Gamma_\theta \cap B_{\CC^n } ( 0 , R_0 ) = B_{\RR^n } ( 0 , R_0 ) \,, \\
\Gamma_\theta \cap \CC^n \setminus B_{\CC^n } ( 0 , 2 R_0 ) =
e^{ i \theta } \RR^n \cap \CC^n \setminus B_{\CC^n } ( 0 , 2 R_0 ) \,,\\
\Gamma_\theta = \{ x + i f_\theta ( x) \; : \; x \in \RR^n \} \,,  \  \
\partial_x^{\alpha}  f_\theta ( x ) =
{\mathcal O}_{\alpha} ( \theta )  \,.
\end{gathered}
\end{gather}
The contour can be considered as a deformation of the manifold $ X $
as nothing is being done in the compact region. The operator $ P$
defines a dilated operator:
\[ P_\theta \stackrel{\rm{def}}{=} P\rest_{ \Gamma_\theta } \,,
\ \  P_\theta u  = \widetilde P ( \tilde u ) \rest_{\Gamma_\theta } \,,\]
where $ \widetilde P $ is the holomorphic continuation of the
operator $ P $, and $ \tilde u $ is an almost analytic extension of
$ u \in \CIc ( \Gamma_\theta ) $ (here we are only concerned with
$ \Gamma_\theta \cap B_{\CC^n } ( 0 , R_0 ) $).

For $ \theta $ fixed, the scaled operator, $ P_\theta $,
 is uniformly elliptic in $ \Psi^{0,2} ( X ) $ outside a 
compact set (see \eqref{eq:sca} below) and hence the resolvent, 
$ ( P_\theta - z )^{-1} $, is meromorphic for  $ z \in D ( 0 , 1 /C ) $.
We can also take $ \theta $ to be $ h $ dependent 
and the same statement holds for $ z \in D ( 0 , \theta / C ) $. 
The spectrum of $ P_\theta $ in $ z \in D ( 0 , \theta / C ) $ 
is independent of $ \theta $ and 
consists of quantum resonances of $ P$ which are defined as the
poles of the meromorphic continuation of
\[  ( P - z )^{-1} \; : \; \CIc ( X ) \; \longrightarrow \; \CI ( X ) \,.\]
 In fact, that is one of the 
ways of defining resonances, and in this paper we will be estimating 
the number of eigenvalues of $ P_ \theta $.

\section{Resonance free regions under the non-trapping assumption}
\label{rfr}

\subsection{Estimates using weight functions}
\label{cop}
We follow the presentation given in \cite[\S 4.2]{DSZ} and inspired by 
many previous works, including \cite{Mart}. 

Let us suppose that $ P_\theta \in \Psi^{0,2} ( X ) $ (we identify 
$ X $ and $ X_\theta $ here) is a complex scaled operator with
$ \theta = M_1 h \log ( 1/h ) $. 
We choose $ \epsilon $ 
\begin{equation}
\label{4}
 \epsilon \le M_2 h\log \frac{1}{h} \,, 
\end{equation}
where $ M_2 > M_1 $ is a large constant to be fixed later. 

Let $ G \in \CIc ( T^* X ) $ and define
\begin{equation*}
P_{\epsilon, \theta}   \stackrel{\rm{def}}{=} 
e^{-\epsilon G/h}P_\theta   e^{\epsilon G/h}=e^{-\frac{\epsilon}{ h}{\rm
ad}_G}P_\theta \sim \sum_0^\infty  \frac{\epsilon ^k}{ k!}(-\frac{1}{ h}{\rm
ad}_G)^k(P_\theta )\,, \  G = G^w ( x , h D) \,.  \end{equation*}
We note that the assumption on $ \epsilon $ and 
the boundedness of  $ {\rm ad}_G/h $ show that the expansion makes sense.
The operators $ \exp ( \epsilon G/ h ) $ are pseudo-differential in 
an exotic class $ S_\delta^{C_2}   $ for any $ \delta > 0 $ (see 
\cite{DiSj}) but that is not relevant here.

Using the same letters for operators
 and and the corresponding symbols, we see
that
\begin{equation*}
{P_{\epsilon,\theta} 
 =P_\theta - i\epsilon \{ P_\theta 
,G\} +{\mathcal O}(\epsilon ^2)=p_\theta - i\epsilon \{
p_\theta ,G\} +{\mathcal O}(h + \epsilon ^2),}
\end{equation*}
so that
\begin{equation*}
\begin{split}
& {\Re P_{\epsilon, \theta}  =\Re p_\theta 
+ \epsilon \{ \Im p_\theta ,G\} +{\mathcal O}(h+\epsilon ^2)} = 
\Re p_\theta + {\mathcal O} ( h + \theta \epsilon + \epsilon^2 ) \,, \\
& {\Im P_{\epsilon, \theta }  =\Im p_\theta - \epsilon \{ \Re p_\theta 
,G\} +{\mathcal O}(h+\epsilon ^2).}
\end{split}
\end{equation*}

We now make the following assumption: for a fixed $ \delta > 0 $
\begin{equation}
\label{12}
|\Re p_\theta | < \delta \ \Longrightarrow \ 
- \Im p_\theta +  \epsilon H_{p} G \geq c_0 \epsilon \,.
\end{equation}

Now let $ \psi_1, \psi_2 \in \CI_{\rm{b}} ( T^*X ) $ be two functions
satisfying 
\[  \psi_1^2 + \psi^2_2 =1 \,, \ \ \psi_1 \rest_{ |\Re p_\theta | < \delta/2 
} \equiv 1 \,, \ \ \supp \psi_1 \subset \{ |\Re p_\theta | < \delta \} \,.\]
Lemma \ref{l:r2} 
gives two selfadjoint operators $ \Psi_1 $ and 
$ \Psi_2 $ with principal symbols $ \psi_1 $ and $ \psi_2 $ respectively, 
such that 
\[ \Psi_1^2 + \Psi_1^2 = I + R \,, \ \  R = {\mathcal O} ( h^\infty ) \; : 
\; H^{-M} ( X) \rightarrow H^{M} \,.\]
We then write $ P_{\epsilon, \theta } = A_{\epsilon , \theta } + i 
B_{\epsilon , \theta } $, where
\[ A_{\epsilon , \theta }  =  \frac12 ( P_{\epsilon , \theta} + 
P_{\epsilon , \theta}^* ) \,, \ \ 
B_{\epsilon , \theta }  =  \frac1{2i} ( P_{\epsilon , \theta} - 
P_{\epsilon , \theta}^* ) \,. \]
The principal symbol of $ B_{\epsilon, \theta }$ is given by 
$ \Im p_\theta - \epsilon H_p G $ and on the essential support of 
$ \Psi_1 $ it is bounded below by $ c_0 \epsilon \gg h $.
Hence  the sharp G\aa{}rding inequality (see \cite[Theorem 7.12]{DiSj}) implies
that for $ h $ small enough
\begin{equation*}
\begin{split}
\| P_{\epsilon, \theta } \Psi_1 u \| \| \Psi_1 u \| 
& \geq | \langle P_{\epsilon, \theta} \Psi_1 u , \Psi_1 u \rangle | 
\geq | \Im  \langle P_{\epsilon, \theta} \Psi_1 u , \Psi_1 u \rangle |\\
& = - \langle B_{\epsilon , \theta} \Psi_1 u , \Psi_1 u \rangle  \geq 
\frac{\epsilon}{C} \| \Psi_1 u \|^2  \,,
\end{split}
\end{equation*}
and hence
\[ \| P_{\epsilon , \theta } \Psi_1 u \| \geq \frac{\epsilon}{C} \| \Psi_1 u \| 
\,.\]
On the support of $ \psi_2 $ the operator $ A_{\epsilon , \theta } $ is
elliptic and  by Lemma \ref{l:r1},
\[ \| P_{\epsilon , \theta } \Psi_2 u \| \geq \frac{1}{C} \| \Psi_2 u \| 
- {\mathcal O} ( h^\infty ) \| u \| 
\,.\]
We conclude from Lemma \ref{l:r3} that
\[ \| P_{ \epsilon, \theta } u \| \geq \frac{\epsilon }{ C} 
\| u \|\,.\]

This shows that the conjugated operator has no spectrum in   
$ D ( 0 , \epsilon / (2 C) ) $.

\subsection{Construction of an escape function}
\label{cef}

Using the results of \S \ref{cop} all we need to do is to construct
$ G $ so that \eqref{12} holds. For that we modify 
a standard argument with the presentation 
borrowed in part from \cite[Sect.4]{VZ}. 

We recall that \eqref{eq:th3a} implies the same condition with 
$ p^{-1} ( [ - \epsilon_0 , \epsilon_0 ] ) $ for some $ \epsilon_0 > 0 $.
That follows from the compactness of the trapped set in 
$ p^{-1} ([ - \delta , \delta ] ) $ -- see \cite[Appendix]{GeSj} for
a detailed discussion.

Let us now fix $ R $ a large parameter.
We will  define $G_\rho \in \CIc ( T^* X ) $, a local escape function
supported in a neighbourhood
of the bicharacteristic segment 
\[ I_\rho =  \{\exp( t H_p)(\rho):\ t\in[-T,T]\}\,,\]
and which satisfies 
$H_p G_\rho \geq 1 $ on the part of $ I_\rho $ lying over
\begin{equation}
K'=\{\rho'\in T^*X \ : \ |x(\rho') | \leq R \}
\end{equation}
For that, let $\Gamma$ be a hypersurface through $\rho$ which is transversal
to $H_p$. Then there is a neighbourhood $U_\rho$ of $\rho$, such that
\[  V_\rho=\{\exp( t (U_\rho\cap\Gamma)):\ t\in(-T-1,T+1)\} \subset
p^{-1} ( [ - \epsilon_0/2 , \epsilon_0/2 ] ) \,, 
\]
is a neighbourhood of $ I_\rho $. That neighbourhood can be identified with a product, 
$$ V_\rho \simeq (-T-1,T+1)\times(U_\rho\cap\Gamma)\,, $$ 
and, in this identification, we will 
choose $ T $  and $ 0 < \alpha < 1 $ so that 
\[ \left( \left( (-T-1 , - \alpha T ) \cup
 ( \alpha T , T+1) \right) \times(U_\rho\cap\Gamma) 
\right)) \cap K' = \emptyset\,.\]
We now need the following elementary 
\begin{lem}
\label{l:chi}
For any $ 0 < \alpha < 1/2 $ 
and $ T > 0 $ there exist as function $ \chi = \chi_{T, \alpha} \in 
\CI ( \RR; \RR ) $ such that 
\[ \chi ( t)  = \left\{ \begin{array}{cc} 0 &  |t | > T \\ 
                                          t  & \ \ |t | < \alpha T
\end{array} \right. \,, \  \ \chi'( t)  \geq - 2 \alpha  \,. \]
\end{lem}
\begin{proof}
The piecewise linear function 
\[ \chi_\# ( t ) = \left\{ \begin{array}{ccc} 0 & \ &  |t | > T \\ 
                                          t  & \ &  \ \ |t | < \alpha T \\
                                          \pm \alpha(  T  - t) / ( 1 - \alpha)
 & \  & \alpha T \leq \pm t \leq T 
\end{array} \right. \]
satisfies $ {\chi_\sharp}' \geq - \alpha / ( 1 - \alpha ) > -2 \alpha $ 
wherever the derivative is defined. A regularization of this function
gives $ \chi_{ T , \alpha } $.
\end{proof}

Now let $\phi_\rho\in\CIc(U_\rho\cap\Gamma)$
be identically $1$ near $\rho$, and let $\chi_T$ be given by the lemma.
Using the product coordinates, we can think of $\phi_\rho$, $ t $, and 
hence $\chi ( t) $, 
as functions on $T^* X$. The functions $ \phi_\rho $ and $ \chi_T ( t ) $
have compact support in $V_\rho$. Let 
$$ \psi \in \CIc ( ( - \epsilon_0 , \epsilon_0 ) )  \,, \ \
\psi \rest_{ [ - \epsilon_0/2 , 
\epsilon_0/2 ]}  \equiv 1\,, $$
and put
\begin{equation}
G_\rho=\chi_T (t)  \phi_\rho \psi(p)  \,, \ \ G_\rho \in \CIc ( V_\rho ) \,.
\end{equation}
so that 
\begin{equation}
H_p G_\rho= \chi_T' \phi_\rho \psi(p), 
\end{equation}
satisfies 
\[ \text{$ H_p G_\rho = 1 $ on $ V_\rho \cap \{ |x| < R \} $ and 
$ H_p G_\rho \geq - 2 \alpha $ everywhere.}\]

Now let $K\Subset T^*X$ be the compact set
\begin{equation}
K=\{\rho\in p^{-1} ([ -\epsilon_0/3 , \epsilon_0/3] ) \ : \ 
|x(\rho)| \leq R/2\}.
\end{equation}
Since $K$ is compact, applying the previous argument for
every $\rho\in K$ gives a $U_\rho$, and a $U'_\rho\subset U_\rho$ on which
$\phi_\rho=1$. Since 
$\{U'_\rho:\ \rho\in K\}$ covers $K$, the compactness
of $K$ shows that we can pass to a finite subcover, $\{U'_{\rho_j}:
\ j=1,\ldots,N\}$. We let
\begin{equation}
G=\sum_{j=1}^N G_{\rho_j}.
\end{equation}
The construction of $ G_{\rho_j}$'s now shows that by choosing $ \alpha $ small
enough (depending on the maximal number of support overlaps
we obtain
\begin{equation}
\label{eq:hpg}
  H_p G ( \rho ) \geq 1 \,, \ \ \rho \in p^{-1} ( ( - \epsilon_0/2 , 
\epsilon_0/2 ))  \cap \{ |x( \rho ) | < R \}  \ \ \text{and} \ \ 
H_p G ( \rho ) \geq - \delta \,, \ \ \rho \in T^* X \,.\end{equation}

\subsection{Resonance free region}
\label{rfr1}

We now want to choose the scaling so that \eqref{12} holds with
$ G $ satisfying \eqref{eq:hpg}. Once that is done the results of
\S \ref{cop} will give Theorem \ref{t:3}.

For that we choose the complex scaling so that
\begin{gather}
\label{eq:sca}
\begin{gathered}
  - \Im p_\theta ( x , \xi ) \geq \theta  \ \ \text{when $ |p ( x ,
\xi )| \leq \epsilon_0 $ and $
|x| \geq R $,}  \\
 \Im p_\theta  < C_1 \theta \ \
 \ \text{when $ |p ( x , \xi )| \leq \epsilon_0 $,}
\end{gathered}
\end{gather}
where $ R $ is independent of $ \theta $.
With $ \epsilon = M_2 h \log( 1/h ) $ we now choose $ \theta =
M_1 h \log ( 1/h ) $ such that
\[  M_1 < M_2 / C_1 \,, \ \  \delta M_2 < M_1 \,, \]
where $ C_1 $ comes from \eqref{eq:sca} and $ \delta $ comes from
\eqref{eq:hpg}. Since we can choose $ \delta $ as small as we want
this can certainly be arranged leading to \eqref{12}.

For completeness we include a quantitative corollary of 
Theorem \ref{t:1} from \cite{NaStZw}:

\medskip
\noindent
{\bf Theorem \ref{t:1}$'$.}{\em Suppose that the assumptions of
Theorem \ref{t:1} are satisfied and that $ ( P_\theta - z )^{-1} $ 
is the scaled resolvent defined for $ 0 < \theta < 2 M h \log h $, 
$ M \gg 1 $. Then for $ 0 < h < h_0 ( M ) $ we have 
\[ \| ( P_ \theta - z )^{-1} \|_{ L^2 ( \Gamma_\theta 
) \rightarrow L^2 ( \Gamma_\theta ) } = C \frac{\exp ( C |\Im z|/ h ) }{h} 
\,, \ \ z \in D ( 0 , M h \log ( 1/h )) \,. \]}

Theorem \ref{t:1} is essentially optimal as shown by the well known 
one dimensional result going back to Regge (see \cite{Z1} for a
proof and references): if $ V \in {\mathcal C}^N ( [a,b] ) $
is extended by $ 0 $ to a potential on $ \RR $, and 
\[  V( x ) \sim \left\{ \begin{array}{ll} ( x - a )^{p} & x \simeq a+ \\
( x - b)^{q} & x \simeq b- \end{array} \right. \,, \ \ p, q < N \,,\]
then the scattering poles for $ - \Delta + V( x ) $ are given at high energies 
by the sequence
\[  \lambda_{k} = \frac{\pi k}{b-a} - i \alpha \log |k| + {\mathcal O}(1) \,, \ \ k \in \ZZ \,, \ \  \alpha =  \frac{ p + q + 4 }{ 2 ( b - a ) } \,.  \]
The semiclassical resonances, $ z_k ( h ) $,
of $ - h^2 \Delta + h^2 V( x) $ are related to these scattering poles by 
the formula $ z_k ( h ) = h^2 \lambda_k^2 $. Hence 
\[ \Re z_k ( h ) \sim 1  \ \Longrightarrow \  \Im z_k ( h ) \sim h \log(1/h) \,.\]

\section{Second microlocal calculus associated to a hypersurface}
\label{smc}

To obtain Theorem \ref{t:2} we need to localize to 
an $ h$-size neighbourhood  of the energy surface 
\[ \Sigma \stackrel{{\rm{def}}}{=} \{ ( x , \xi ) \in 
T^* X \; : \; p ( x , \xi ) = 0 \}
\,. \] 
That means that we have to work with functions of the form 
$$ a ( x , \xi ; h ) = \psi ( p ( x , \xi ) / h ) \,.$$ 
The usual quantization procedure
(the passage from symbols to pseudodifferential operators) is 
prohibited as 
$$ \partial_{x, \xi}^\alpha a \simeq h^{-|\alpha| } \,. $$

The troublesome symbols have a special form and we can construct
a calculus which includes them by straightening $ \Sigma $ locally 
by means of canonical 
transformations. That means moving $ \Sigma $ 
 to 
$$ \Sigma_0 = \{ \xi_1 = 0 \} \,. $$ 
We may then localize
to rectangles in the $ ( x_1 , \xi_1 ) $-space of length $ \sim 1 $ in 
$ x_1 $ and of length $ \sim h $ in $ \xi_1$. This amounts to
a form of  semiclassical second microlocalization. 
The presentation here is essentially self-contained and we refer to  
\cite{SjZw99} for pointers to the literature.

For Theorem 2 we need even more singular calculus related to 
the $ \Psi_{\frac12} $ calculus described in 
\S \ref{sma}.

We assume that $ \Sigma $ is a compact $ \CI $ hypersurface
in $ T^*X $. 
Since the delicate constructions will only be used in a compact
set and since we are working in  the $ \CI $ category this
creates no restrictions.

\subsection{Basic properties} \label{s.smr}
To construct the calculus, let $ \Sigma \Subset  
T^* X $ be a $ \CI $ compact hypersurface. 
We consider a class of symbols associated to $ \Sigma $, a multiindex
\[ {\mathfrak p} \stackrel{\rm{def}}{=} ( m , {\widetilde m}, k_1, k_2 ) \,,\]
and 
depending on {\em two} small parameters, 
$$ 0 < h < \tilde h \,.$$ 
\begin{equation} 
\label{eq:9.0} 
 a \in \Ssd ( T^* X ) \Longleftrightarrow   
\left\{ 
\begin{array}{l} \text{\underline{near $\Sigma$} } \ \ \
V_1 \cdots V_{l_1} W_1 \cdots W_{l_2 } H_q^k a = \\  {\mathcal O} \left(  
h^{ -m - \delta l_1 -l_2 } \tilde h^{-{\widetilde m}+ \delta l_1 + l_2}
 \langle
( \tilde h/ h) q \rangle^{k_2 + 
\delta k - ( 1- \delta ) l_2 }\right)  \,, \\
\text{where} \ V_1, \cdots, V_{l_1 } \ \text{ are vector fields  
tangent to $ \Sigma $ and } \\
 W_1 , \cdots W_{l_2 } \ \text{ are 
any vectorfields } \\
\ \ \\ 
\text{\underline{away from  $\Sigma$} } \ \ 
\pxx a ( x , \xi ; h ) = \\
{\mathcal O } ( h^{-m  - k_2 - \delta ( |\alpha| 
+|\beta | )} \tilde 
h ^{-{\widetilde m}+k_2 + \delta ( |\alpha| +|\beta | )} 
 \langle \xi \rangle^{k_1 - |\beta | } ) \,. 
\end{array} 
\right.  
\end{equation} 
Here $ q $ is 
any defining function of $ \Sigma $, that is a 
function which vanishes simply on $ \Sigma $. To assure invariance 
we need
\begin{lem}
\label{eq:l.3.inv}
Definition \eqref{eq:9.0} is independent of the choice of $ q $ and
the vector fields applied to $ a $ can be taken in any order.
\end{lem}
\begin{proof}
The independence of the order of vector fields follows from the fact
that $ [V_{\ell_j}, H_q ] $ is a vector field tangent to $ \Sigma $.
To see independence of $ q $ we note that $ H_{uq} = u H_q + q H_u $.
The vector field $ H_u $ can be considered as an arbitrary 
vector field and $ q H_u $ is tangent to $  \Sigma $. 
Hence the application of the second term is estimated by 
$$ | q H_u a |
\leq C \min ( |q| h^{-1} \langle q h^{-1} \rangle^{-( 1 - \delta )} 
, h^{-\delta } ) \leq C' \langle h^{-1} q \rangle^{ \delta } 
\,, $$ 
and this estimate can be iterated. Hence we have the same estimates
for $ q $ replaced by $ u q $, $ u \neq 0 $ near $ \Sigma $.
\end{proof}

The symbol classes in \eqref{eq:9.0} are best understood in the
simple case when $ q = \xi_1 $, say. The condition
for  $ |\xi| \leq C $,  means that, with $ \xi = ( \xi_1 , \xi' ) $, 
\begin{equation}
\label{eq:9.0'}
\partial_{x_1}^k \partial_{x'}^\alpha \partial_{\xi'}^\beta 
( \xi_1 \partial_{\xi_1} )^{\ell_1} \partial_{\xi_1}^{\ell_2} 
a = {\mathcal O} \left(   h^{-m }
\tilde h^{-{\widetilde m}}
 ( \tilde h/ h )^{\delta( |\alpha| + |\beta| ) + \delta \ell_1 + 
 \ell_2 } \langle
( \tilde h / h ) \xi_1\rangle^{k_2 +  \delta k - ( 1- \delta) 
\ell_2) } \right) \,,
\end{equation}
or, if we eliminate the vector fields vanishing at $ \xi_1 =0 $,
\begin{equation}
\label{eq:9.0''}
\partial_{x_1}^{k}
\partial_{x'}^\alpha
 \partial_{\xi'}^\beta ( ( h/ \tilde h) \partial_{\xi_1}) ^p
 a = {\mathcal O} \left(   h^{-m } \tilde h^{-{\widetilde m}}
( \tilde h / h )^{ \delta ( |\alpha| + |\beta|) } 
\langle 
( \tilde h / h ) \xi_1\rangle^{ k_2 - ( 1- \delta) p + \delta k } \right) \,.
\end{equation}
We used the fact  that if $ | \xi_1 | \leq C $ then 
\[ \left( ( \tilde h / h ) \xi_1 \right)^{ \delta -1  } \leq 
 \left( ( \tilde h / h ) \xi_1 \right)^{-1} ( \tilde h / h )^{\delta } \,,\]
to eliminate the need for the $ \xi_1 \partial_{\xi_1 } $.
The advantage of the formulation in \eqref{eq:9.0} (and \eqref{eq:9.0'})
is the geometric invariance. In explicit coordinates \eqref{eq:9.0''} is
however more transparent.

Since it is sufficient for our purposes, and for simplicity of 
presentation we will consider the case of $ \delta = 1/2 $ only.

The second microlocalization associates to the space of 
symbols $ S^{\mathfrak p}_{ \Sigma, \frac12} $ a space
of operators, $ \Psi_{\Sigma, \frac12}^{\mathfrak p} $, defined
in \S \ref{gc} below. The basic properties of these spaces
are described in the following  
\begin{thm} 
\label{th:3.1} 
 Let us define two 
multiindices,
\[ {\mathfrak p } = ( m , {\widetilde m}, k_1, k_2 ) \,, \ \ 
 {\mathfrak p }' = ( m  , {\widetilde m} - 1,  k_1-1, k_2 ) \,.\]
With the definitions of $ S^{\mathfrak p}_{\Sigma,\delta} $ above and
$ \Psi^{\mathfrak p}_{\Sigma,\delta } $ in \S \ref{gc} below, 
there exist maps
\begin{align*} 
& \Ops \; : \; \Ssh ( T^* X ) \longrightarrow \Pkh ( X ) \\ 
& \sigs \; : \; \Pkh ( X ) \longrightarrow \Ssh ( T^* X )  
/ \Ssgh   ( T^*X )    
\end{align*} 
such that  
\begin{equation} 
\label{eq:3.1h} 
\sigs ( A\circ B ) = \sigs ( A) \sigs ( B ) \,  
\end{equation} 
\begin{gather*}
0 \longrightarrow \Pkgh ( X ) \longrightarrow \Pkh ( X)  
\stackrel{\sigs}{\longrightarrow}   \Ssh ( T^* X )  
/  \Ssgh   ( T^* X )  \longrightarrow 0  
\end{gather*}
is a short exact sequence and 
\[ 
\sigs \circ \Ops \; : \;  \Ssh  ( T^* X )  \longrightarrow \Ssh  
( T^ * X )  / \Ssgh    ( T^* X )
\]
is the natural projection map. 
If 
$$ a \in \Ssh ( T^* X ) \,,  \ \ 
 d ( \supp  a  , \Sigma)  \geq 1/C  $$ then  
$$ a  \in  S_{\frac12}^{ m + k_2,  \widetilde m - k_2 , k_1 } ( T^* X) \,, \ \ 
\Ops ( a ) = \Op ( a ) \in  \Psi_{\frac12} ^{ m + k_2 ,
 \widetilde m - k_2,  k_1 } ( X ) \,, $$ 
where $ \Psi_{\frac12} ^{\bullet, \bullet, \bullet }  
( X) $ is the  class of 
pseudodifferential operators defined in \S \ref{sma}.
\end{thm}

\subsection{Calculus in the model case}
\label{cmc}
To define $ \Pkh ( X) $ we proceed locally and put $ \Sigma $ into 
a normal form $ \Sigma_0 = \{ \xi_1 = 0 \} $ (locally). 

The model case is obtained by taking symbols satisfying 
\eqref{eq:9.0''} globally and defining
\[ a = a ( x , \xi, \lambda , h, \tilde h ) \,, \ \ \lambda = \tilde h \xi_1 
/h \,, \]
satisfying
\begin{equation}
\label{eq:newOh}
\partial_{x_1}^k \partial_{x'}^{\alpha} \partial_\xi^\beta \partial_\lambda^p  
a ( x, \xi , \lambda ; h ) = {\mathcal O  } ( h^{- m  }  \tilde 
h^{-\widetilde m} ) (\tilde h / h )^{ (|\alpha| + |\beta|)/2} 
\langle \lambda \rangle ^{k_2 +  k/2 -  p/2 }   \,,\end{equation}
which is the same as
\eqref{eq:9.0}. 
We will write \eqref{eq:newOh} as 
\[ a \in \widetilde{\mathcal O}_{\frac12}
  \left( h^{-m} \tilde h^{- \tilde m } 
\langle \lambda\rangle ^{k_2 } \right) \,. \]
and define an exact quantization in the usual way,
\begin{equation}
\label{eq:3.4''}
\begin{split} 
& \Opt (  a ) u ( x ) = \\
& \ \  \frac{1}{2 \pi \tilde h 
}  \frac{1}{ ( 2 \pi h ) ^{n-1} } \int a ( x, \xi' , (h/ \tilde h) 
 \lambda  , \lambda  ; 
h ) e^{{i} \langle x' - y' , \xi' \rangle/h   + 
{i}  \langle x_1 - y_1, \lambda \rangle/\tilde h } u ( y ) d y  
d \xi' d \lambda  \,, \end{split}
\end{equation}
$ n = \text{dim}\; X $, and 
where 
$$ \lambda = ( \tilde h/ h) \xi_1 \,. $$
For $ a \in \widetilde {\mathcal O} ( \langle \lambda \rangle^{k_2} ) $
and $ b \in \widetilde {\mathcal O} ( \langle \lambda \rangle^{k_2'} ) $
we  have 
 \begin{gather} 
\label{eq:3.5} 
\begin{gathered} 
\Opt (a) \circ
 \Opt (b) = \Opt ( a \sharp_{h , \tilde h } b ) \,, \ \ 
a \; \sharp_{h, \tilde h}  \; b = \widetilde{\mathcal O} 
 \left( \langle \lambda  
\rangle^{ k_2 + k'_2 }  \right) \,,
 \\
a \; \sharp_{h,\tilde h} \; b ( x, \xi , \lambda ; h ) \sim 
\sum_{\alpha \in  \NN ^n } \frac{1}{\alpha!} ( h \partial_{\xi' } )^{\alpha'} 
(   h \partial_{\xi_1} +  \tilde h  \partial_\lambda )^{\alpha_1} 
 a  \,   D_x^\alpha b  \,, 
\end{gathered} 
\end{gather} 
where the asymptotic sum is defined up to terms in 
\[   \widetilde {\mathcal O}_{\frac12}
 (\tilde h^\infty \langle \lambda \rangle^{k_2 +k'_2} 
)   \,. \]
To see this we write $ \Opt ( a ) $ as a quantization of an 
operator valued symbol, $ \Oph ( a ) ( x_1 , \lambda ) $, 
\[ \Oph ( a ) ( x_1, \lambda ) = \frac{1}{( 2 \pi h )^{n-1} } 
 \int a ( x, \xi' , (h/ \tilde h) 
 \lambda  , \lambda  ; 
h ) e^{{i} \langle x' - y' , \xi' \rangle/h  } d \xi' \,, \]
so that 
\begin{equation}
\label{eq:tensq} \Opt ( a ) = \Oph (a ) ( x_1 , \tilde h D_{x_1} ) \,. 
\end{equation}
In view of \eqref{eq:bealss} we have 
\begin{gather*}
 a \in \widetilde {\mathcal O}_{\frac12} ( \langle \lambda\rangle^{k_2} ) 
\; \Longleftrightarrow \; \\
\| \partial_{x_1}^k \partial_{\lambda}
^p \Oph (\partial^\alpha a) ( x_1, \lambda )  \|_{ L^2 ( \RR^{n-1} ) \rightarrow
L^2 ( \RR^{n-1} ) } \leq C_{\alpha, p, k } ( \tilde h / h )^{|\alpha|/2} 
\langle \lambda \rangle^{k_2 + k/2 - p/2 } \,. \end{gather*}
Observing that 
\begin{gather*}
  \Opt ( a ) \circ \Opt ( b ) = \\
\left( \exp ( i \tilde h \langle 
D_\lambda , D_{y_1} \rangle ) 
 \Oph a ( x_1, \lambda ) \circ \Oph b ( y_1, \lambda' ) \rest
_{ x_1 = y_1, \lambda = \lambda' } \right) 
( x_1 , \tilde h D_{x_1} )  
\,, \end{gather*}
proves \eqref{eq:3.5}.

%


For the operator 
$ \Opt ( a ) $, with $ a $ supported in $ |\xi_1 | \leq C $,  we define
its principal symbol 
as the equivalence class of $ a $ in  
$$ \widetilde{\mathcal O }_{\frac12}  ( \langle\lambda\rangle^{k_2} ) / 
 \widetilde{\mathcal O}_{\frac12}
 (  \tilde h  \langle \lambda \rangle^{k_2} ) 
   \,. $$
 In view of \eqref{eq:3.5} this symbol map is 
a homomorphism onto the quotient
of symbol spaces:
\begin{equation}
\label{eq:3.sy}
\Opt ( a ) \mapsto [a] \in
  \widetilde{\mathcal O }_{\frac12}  ( \langle\lambda\rangle^{k_2} ) / 
 \widetilde{\mathcal O}_{\frac12}
 ( \tilde h  \langle \lambda \rangle^{k_2  } ) 
\end{equation}
We note that in the local model 
we are not concerned here with the behaviour as
$ |\xi| \rightarrow \infty $.


It will be useful to have the analogue 
of Beals's characterization of pseudodifferential operators by
stability under taking commutators. It follows from the
proof of the semiclassical analogue of Beals's result  \cite{HeSj1} and its
adaptations in \cite[Chapter 8]{DiSj} and \cite[Lemma 4.2]{SjZw99}.
\begin{prop} 
\label{p:Beals} 
Let $ A = A_h : {\mathcal S}( \RR^n ) \longrightarrow {\mathcal S}' ( \RR^n) 
$ and put $ x' = ( x_2 , \cdots,  x_n ) $. 
For $ p \in \RR $, we define the norms
\[    \|u \|_{(p)} = \| \langle \tilde h D_{x_1} \rangle^p u \|_{L^2} \,.\]
Then
\[ \text{ $A = \Opt ( a ) \  $ for 
$  \ a = \widetilde{\mathcal O }_{\frac12} ( \langle \lambda
\rangle^k )  $, } \ \  \ 
 a = a ( x , \xi' , \lambda ; h , \tilde h ) \,, \]
if and only if for all $ N, p, q \geq 
0 $ and every sequence $ \ell_1 ( x' , \xi' ) , \cdots ,  
\ell_N ( x', \xi' ) $ of linear forms on $ \RR^{2 ( n-1 ) } $ there 
exists $ C> 0 $ for which 
\begin{equation}
\label{eq:beals}
\begin{split}
& \dn \ad_{\ell_1 ( x' , h D_{x'} ) } 
\circ \cdots \circ \ad_{\ell_N ( x' , h D_{x'} ) }  
\circ (\ad_{{\tilde h} D_{x_1} } )^p  \circ ( \ad_{x_1} )^q  
A u \dn_{ ( q/2 - \min( k , 0 ) ) } \\
& \ \ \ \ \ \ \ \ \ \ \leq C h^{N/2} \tilde h^{N/2 + p + q } 
\dn u \dn_{p/2 + (\max(k , 0 ) ) } \,.
\end{split}
\end{equation}
\end{prop} 
\begin{proof}
We first observe that for $ A =  A = \Opt ( a ) $ for $ a = \widetilde{\mathcal O }_{\frac12} ( \langle \lambda 
\rangle^k )  $, \eqref{eq:beals} follows from the calculus.

Let 
$$ V_{h} u ( x_1 , \tilde x' ) \stackrel{\rm{def}}{{=}}
 h^{-\frac{n-1}4} u ( x_1 , h^{\frac12} \tilde x' ) $$ 
so that 
\[ A = V_h ^{-1}  \widetilde A V_h \,,\]
where, as in the proof of Proposition \ref{l:1hf}, we have
\begin{equation}
\label{eq:11.5} 
\begin{split}
& \dn \ad_{\ell_1 ( x' ,  D_{x'} ) } 
\circ \cdots \circ \ad_{\ell_N ( x' ,  D_{x'} ) }  
\circ (\ad_{ D_{x_1} } )^p  \circ ( \ad_{x_1} )^q  
\widetilde A u \dn_{ ( q/2 - \min( k , 0 ) ) } \\
 & \ \ \ \ \ \ \ \ \ \ 
\leq C  \tilde h^{N/2 +  q } 
\dn u \dn_{p/2 + (\max(k , 0 ) ) } \,.
\end{split}
\end{equation}
We can write $ \widetilde A =  a_{h, \tilde h} ( x, D_x )$,
so that, $ A = \Opt(a ) $, where
\[ a ( x , \lambda , \xi' ; h ) = a_{h,\tilde h } ( x_1 , h^{-1/2} 
x' ; \tilde h^{-1} \lambda , h^{-1/2} \xi' )  \,.\]
The required estimate on $ a $ becomes
\begin{equation}
\label{eq:reqa}
  \partial_{x_1}^p \partial^q_{\xi_1}  \partial_{\tilde x'} 
^{\alpha'} \partial_{\tilde \xi'}^{\beta'} a_{ h , h'} = 
{\mathcal O} ( 1) \tilde h^{ (|\alpha'| + |\beta'|)/2 + q } 
 \langle \tilde h \xi_1 \rangle^{k -q/2 + p/2} 
\end{equation}

We have 
\begin{equation} 
\label{eq:a.0} 
 \langle \widetilde A \psi , \phi \rangle = \frac{1}{ ( 2 \pi )^n }  
\int \int e^{ i \langle x , \xi \rangle } a_{h,\tilde h} ( x , \xi )  
\widehat \psi ( \xi ) \overline{\phi ( x ) } d x d \xi \,,  
\end{equation} 
with 
$$ \widehat \psi ( \xi ) = ({\mathcal F} \psi ) ( \xi ) =  
\int e^{ - i \langle x , \xi \rangle } \psi ( x ) d x \,,$$ 
and for $ \phi, \psi \in {\mathcal S} ( \RR^n ) $. 
Let us fix $ ( x_0, \xi_0) \,, (y_0, \eta_0) \in T^* \RR^n $, 
and $ \lambda \gg 1 $. With
$ \phi , \psi \in {\mathcal S} ( \RR^n ) $ we put
\[ 
\begin{split}
&  \psi_{   x_0, \xi_0 }( x)  = \lambda^{\frac14} 
\psi \left( \lambda^{\frac12}( x_1 - x_{0,1} ) , 
x' - x'_0 
\right) e^{ i \langle 
 x , \xi_0 \rangle }  
\,, \\ 
& \phi_{   y_0, \eta_0 }( x ) = \lambda^{\frac14} 
\phi \left( \lambda^{\frac12} ( x_1 - y_{0,1} ) , x' - y'_0 
\right) 
e^{ i \langle  x - y_0 , \eta_0 \rangle }  \,.\end{split} \]
We see that 
\[ \begin{split}
& \widehat \psi_{   x_0, \xi_0 }( \xi )  = 
\lambda^{-\frac14} 
\hat \psi \left( \lambda^{-\frac12}( \xi_1 - \xi_{0,1} )
 , \xi' - \xi'_0 \right) e^{ - i \langle 
\xi - \xi_0 , x_0 \rangle }  \,, \\
& \widehat \phi_{   y_0, \eta_0 }( \xi )  = \lambda^{-\frac14} 
\hat \phi \left( \lambda^{-\frac12}( \xi_1 - \eta_{0,1} )
, \xi' - \eta_0' \right) e^{ - i \langle \xi , y_0 \rangle }  \,. 
\end{split} \]
We have   
\begin{gather*} 
B \stackrel{\text{def}}{=} \text{ad}_{x'}^{\alpha'}  
\text{ad}_{D_{x'}}^{\beta'}  \text{ad}_{x_1}^q \text{ad}^p_{ 
D_{x_1}} \widetilde A = \left(-i\right) ^{ |\alpha'| +  
|\beta'| + q + p } 
b_{h,\tilde h} ( x , D) \,, \\ 
b_{h,\tilde h} ( x, \xi ) = ( - \partial _{\xi'} )^{\alpha' }  
 \partial_{x'} ^{\beta'} ( - \partial_{\xi_1} )^q \partial_{x_1}^p 
a_{h,\tilde h} ( x, \xi ) \,.  
\end{gather*} 
Let us now assume we have the commutator estimate in the lemma with  
$ k \geq 0 $. Since 
$$ \dn u \dn_{( - q /2 ) } = \| \langle  
\tilde h D_{x_1} \rangle^{- q/2 } u \|_{L^2 } $$
is the norm dual  
to $ \dn \bullet \dn_{( q/2 ) } $, we get from \eqref{eq:11.5}
\begin{equation} 
\label{eq:a.1} 
| \langle B  \psi_{x_0, \xi_0}  ,  
 \phi_{y_0, \eta_0}  \rangle | \leq 
\tilde h ^{  (|\alpha'| + |\beta'|)/2 + q  } 
\dn  \psi_{x_0, \xi_0} \dn_{( p /2 + k ) }  
\dn \phi_{y_0, \eta_0 }\dn_{ ( - q/2 ) }  
\,.  
\end{equation} 
Let $ \bullet_{0,1} $ denote the first component of 
$ \bullet_0 \in \RR^n $. 
For  fixed $ \psi, \phi \in {\mathcal S}( \RR^n ) $, we have 
\begin{gather}
\label{eq:fixed} \begin{gathered} 
 \| \psi_{x_0, \xi_0 } \|_{(p/2 + k)}^2  \leq C_N I_{p/2+k}^N 
( \lambda, \xi_{0,1} ) \,, 
\ \ \ 
 \| \phi_{y_0 ,\eta_0 } \|_{(-q/2)}^2
 \leq  
 C_N I_{-q/2}^N ( \lambda, \eta_{0,1} ) 
\\
I_{r}^N ( \lambda, \tau)  \stackrel{\rm{def}}{=} 
\lambda^{-\frac12} \int_{\RR} \langle \tilde h \rho 
\rangle^{2r} \langle \lambda^{-\frac12} 
( \rho - \tau  ) \rangle^{-N} d \rho 
\end{gathered}
\end{gather}
Using \eqref{eq:a.0} we rewrite the left hand side of \eqref{eq:a.1} as  
\[ \frac{1}{  ( 2 \pi )^n}
\left|
\int \int e^{ i \langle x , \xi \rangle }  
b_{h,\tilde h} ( x, \xi ) \widehat \psi_{x_0, \xi_0} ( \xi) 
\bar \phi_{y_0, \eta}  ( x  )  dx d \xi \right|   \,. \] 
Decomposing the first 
exponent in the integral as  
\[ \langle x , \xi \rangle = \langle y_0 , \xi_0 \rangle +  
\langle x - y _0 , \xi_0 \rangle + \langle \xi - \xi_0 , y_0 \rangle 
+ \langle x - y_0 , \xi - \xi_0 \rangle \] 
and using the formul{\ae} for $ \widehat \psi_{x_0, \xi_0} $, 
$ \phi_{y_0, \eta_0} $,
we rewrite it further as  
\begin{gather*} 
\frac{1}
{  ( 2 \pi )^n} \left|  
\int \int b_{h,\tilde h} ( x , \xi ) \exp({ i \langle x - y_0 , 
\xi - \xi_0 }) \widehat \psi ( \lambda^{-\frac12} ( \xi_1 - \xi _{0,1} ) 
, \xi' - \xi'_0  ) 
 \right.
\\
\left.
\bar\phi ( \lambda^{\frac12} (x_1- y_{0,1} ), x' - y'_0 ) 
\exp({ i ( \langle \xi - \xi_0 , y_0 - x_0 \rangle +  
\langle x- y_0 , \xi_0 - \eta_0 \rangle )}) dx d \xi \right|  
\end{gather*} 
Summing up, we get  
\begin{gather}
\label{eq:a.2} 
\begin{gathered}
 \left| 
{\mathcal F}{\left( \chi b_{h, \tilde h } \right)} ( \eta_0 - \xi_0 , x_0 - y_0 ) \right| = 
{\mathcal O} ( 1)  \tilde h^{ (|\alpha'|
+ |\beta '|)/2 +    q} 
\| \psi_{x_0,\xi_0} \|_{p/2 + k } \| \phi_{ y_0, \eta_0 } \|_{ - q/2} \,, \\ 
\chi ( x, \xi ) = e^{ i \langle x - y_0, \xi - \xi_0 \rangle }
\widehat \psi ( \lambda^{-\frac12} 
( \xi_1 - \xi_{0,1}), \xi' - \xi'_{0} )  \bar \phi ( 
\lambda^{\frac12} ( x_1 - y_{0,1} ) , x' - y_0' ) \,.
\end{gathered} \end{gather}
Writing 
$$ \zeta_1 \stackrel{\rm{def}}{=} \frac{\eta_{0,1} -  \xi_{0,1}}{
\lambda^{\frac12} } \,, \ \ z_1 \stackrel{\rm{def}}{=} \lambda^{\frac12}
( x_{0,1} - y_{0,1} )\,, \ \ \zeta' = \eta'_0 - \xi_0' \,, \ \
z' = x_0' - y'_0 \,, $$
 we therefore get 
\[   {\mathcal F} \left(  
\widetilde \chi \widetilde b_{h,\tilde h}  \right)
 ( \zeta , z ) =   
{\mathcal O} ( 1)  \tilde h^{ (|\alpha'|
+ |\beta '|)/2  +   q} 
\dn  \psi_{x_0, \xi_0} \dn_{( p /2 + k ) }  
\dn \phi_{y_0, \eta_0 }\dn_{ ( - q/2 ) }  
 \,, \]
where 
\[ \widetilde b_{h, \tilde h }  ( X , \Xi ) \stackrel{\rm{def}}{=} 
b_{ h , \tilde h } ( y_0 + \lambda^{-\frac12} X , \xi_0 + 
\lambda^{\frac12} \Xi ) \,, \]
and 
\[  \widetilde \chi ( X , \Xi ) \stackrel{\rm{def}}{=} 
e^{ i \langle X, \Xi \rangle} \widehat \psi ( \Xi ) \bar \phi ( X ) \,.\]
We then have  
\begin{align*}
&  {\mathcal F} \left(  
\partial_{\Xi'}^{\tilde \alpha' } \partial _{\Xi_1}^{\tilde q } 
\partial_{X'} ^{\tilde \beta'} \partial_{X_1}^{\tilde p }  
\widetilde \chi  
\widetilde b_{h,\tilde h} \right) ( \zeta , z ) =  \\
& \ \ \ \ \ 
{\mathcal O}_N ( 1)  \tilde h^{ (|\alpha'|
+ |\beta '|)/2  +   q} 
I_{( p /2 + k + \tilde p/2 ) }  ^N ( \lambda, \xi_{0,1} ) 
I_{ ( - q/2 - \tilde q/2 ) }  ( \lambda, \eta_{0,1} ) 
\lambda^{-\tilde p/2 + \tilde q/2 } 
 \,,\end{align*} 
which putting 
$$ \tilde \alpha = ( \tilde q , \tilde \alpha' ) \,, \ \
 \tilde \beta = ( \tilde p , \tilde \beta ' ) \,, $$ 
we rewrite as 
\[
\begin{split}
&  z^{\tilde \alpha } \zeta^{\tilde \beta}  
{\mathcal F} \left(  
\widetilde \chi   
\widetilde b _{h,\tilde h} \right) ( \zeta , z ) = \\ 
& \ \ \ 
{\mathcal O} ( 1) 
  \tilde h^{ (|\alpha'|
+ |\beta '|)/2  +   q} 
I_{( p /2 + k + \tilde p/2 ) }  ^N ( \lambda, \xi_{0,1} ) 
I_{ ( - q/2 - \tilde q/2 ) }  ( \lambda, \eta_{0,1} ) 
\lambda^{-\tilde p/2 + \tilde q/2 } 
 \end{split}
\,. \] 
We now go back to \eqref{eq:fixed} and choose $ \lambda = \langle
\tilde h \xi_{0,1} \rangle $. Then 
\[ \begin{split}
I^N_{(p/2 + \tilde p/2 + k)} (  \langle \tilde h \xi_{0,1} \rangle, 
\xi_{0,1} ) 
& = 
\left(  \int_{\RR} \langle \tilde h ( \xi_{0,1} + \langle \tilde h\xi_{0,1} 
\rangle^{\frac12} r  ) \rangle^{ p + \tilde p + 2 k } 
\langle r \rangle^{-N} dr \right)^{\frac12} \\
& =
C_N 
\left(  \int_{\RR} \langle R + \langle R\rangle ^{\frac12} \tilde h r  ) 
\rangle^{ p + \tilde p + 2 k } 
\langle r \rangle^{-N} dr \right)^{\frac12} \rest_{ R = \tilde h \xi_{0,1}}  \\
& \sim   \langle R \rangle ^{(p + \tilde p + 2 k )/2 } 
\rest_{ R = \tilde h \xi_{0,1}}  = \langle \tilde h \xi_{0,1} \rangle^{ p/2 
+ \tilde p /2 + k }  \,. 
\end{split} \]
It follows that for any $ \tilde \alpha $ and $ N $ we have
\begin{equation}
\label{eq:follows}
 | z^{\tilde \alpha } {\mathcal F} \left(  
\widetilde \chi 
\widetilde 
b_{h,\tilde h} \right) ( \zeta , z ) | 
 =  {\mathcal O}_N ( 1)
   \tilde h^{ (|\alpha'|
+ |\beta '|)/2  +   q} 
I^N_{ ( - q/2 - \tilde q/2 ) }  
\lambda^{ \tilde q/2  + p/2 + k } 
\langle \zeta \rangle^{-N} 
\,, \ \ \lambda = \langle \tilde h \xi_{0,1} \rangle \,. \end{equation}
We now integrate the left hand side 
in $ \zeta_1 = ( \eta_{0,1} - \xi_{0,1} )/ \lambda^{\frac12}  $:
\[ \begin{split}
& \int_{\RR } I_{-q/2 -\tilde q /2 } ^N  \left\langle \frac{
\eta_{0,1} - \xi_{0,1} }{ \langle \tilde h \xi_{ 0, 1}
 \rangle^{\frac12} } \right\rangle^{-N }
 \langle \tilde h \xi_{ 0, 1} \rangle^{- \frac12} 
 d  \eta_{0,1}  \\
& \ \ \ =  
\int_{\RR }  
\left(  \int_{\RR} \langle \tilde h ( \eta_{0,1} + \langle \tilde h 
\xi_{0,1} 
\rangle^{\frac12} r  ) \rangle^{ - q - \tilde q  } 
\langle r \rangle^{-N} dr \right)^{\frac12} 
\left\langle \frac{
\eta_{0,1} - \xi_{0,1} }{ \langle \tilde h \xi_{ 0, 1} \rangle^{-N}
} \right\rangle^{\frac12} \langle \tilde h \xi_{ 0, 1} \rangle^{- \frac12}
 d  \eta_{0,1}  \\
& \ \ \ \leq 
C'_N \int_\RR \left( \int_\RR \langle \tilde h ( \xi_{ 0, 1 } + 
\langle \tilde h \xi_{ 0,1} \rangle^{\frac12} ( r + \zeta_1 ) )
\rangle^{-q/2 - \tilde q/2 }  \langle r \rangle^{-N}
d r \right)^{\frac12}  \langle \zeta_1 \rangle^{-N} 
d  \zeta_1 \\
& \ \ \ \leq \langle \tilde h \xi_{ 0 , 1} \rangle^{- q/2 - \tilde q/2 } \,.
\end{split} \]
Returning to \eqref{eq:follows} we see that 
\[ \int_{\RR^n} | {\mathcal F} \left(  
\widetilde \chi 
\widetilde 
b_{h,\tilde h} \right) ( \zeta  , z ) | d \zeta 
 =  {\mathcal O} ( 1)
   \tilde h^{ (|\alpha'|
+ |\beta '|)/2  +   q}  \langle \tilde h \xi_{0,1} \rangle^{p/2 + k - q/2 }
\langle z \rangle^{-N} \]
and consequently
\[ \widetilde  b_{h , \tilde h} \widetilde \chi ( 0 , 0 ) = 
  {\mathcal O} ( 1)  \tilde h^{ (|\alpha'| + |\beta '|)/2  +   q}
\langle \tilde h \xi_{0,1}\rangle^{k + p/2 - q/2 } 
\,.\]
Combining this with 
the definition of $ b_{h , \tilde h } $ and $ \widetilde b_{ h , \tilde h }
$ gives 
\[ \partial_{\xi'}^{\alpha'} 
\partial _{\xi_1}^q \partial_{x'}^{\beta'} \partial_{x_1}^p  
a_{h,\tilde h}  ( y_0 , \xi_0 ) =  
{\mathcal O} ( 1) 
   \tilde h^{ (|\alpha'|
+ |\beta '|)/2  +   q} 
\langle \tilde h \xi_{0,1}\rangle^{k + p/2 - q/2 } 
\,, \]
which is \eqref{eq:reqa} for $ k \geq 0 $. When $ k < 0 $ we check that 
the assumption is satisfied for $ \langle {\tilde h} D_{x_1}  
\rangle^{-k} A $ with $ k $ replaced by $ 0$. The composition  
formula then gives the result. 
\end{proof}

\subsection{Modified Sobolev spaces}
\label{msc}
The norm used in Proposition \ref{p:Beals} can be defined 
globally. We generalize it to include standard Sobolev spaces
by adding additional information at a smooth compact hypersurface
$ \Sigma \subset T^* X $.

Let $ Q \in \Psi_h^{0,0} $ be an operator, $ Q = Q^* $, 
with the principal
symbol $ q $ satisfying
\begin{equation}
\label{eq:pq}
\Sigma = \{( x , \xi ) \; : \;  q ( x , \xi ) = 0\} 
 \,, \ \ dq \rest_\Sigma \neq 0\,, \ \ 
| q ( x, \xi )|   \geq 1  \ \text{ for 
$ { |\xi| \geq C } $ } \ 
\,. \end{equation}

For $ s, m \in \RR $ we  define
\begin{equation}
\label{eq:defh}
 H_\Sigma^{s, m } ( X) = \{ u \in 
{\mathcal S}' ( X ) \; : \; \langle ( \tilde h/ h) Q \rangle^m u \in 
H^s_h ( X ) \} \,, \ \ 
\| u \|_{  H_\Sigma^{s, m } ( X) } \stackrel{\rm{def}}{=} 
\| \langle ( \tilde h/h) Q \rangle^m u \|_{ H^s_h  ( X ) } \,.
\end{equation}
Here $ H^s_h ( X ) $ denotes the usual semiclassical Sobolev spaces
defined in \S \ref{rsa}. 
The spaces are (complex) interpolation spaces, in $ m $, and in $ s $.

When $ m \in \ZZ $  the definition is equivalent to 
\[  H_\Sigma^{s, m } ( X) {=} \left\{
\begin{array}{ll}
\{ u \; : \; (\tilde h/h )^k Q^k u \in H^s_h ( X )\,, \ 0 \leq k \leq m  \} &
 m \geq 0 \\
\ & \ \\
\{ u \; : \; u = \sum_{k=0}^{ |m | } ( \tilde h /h )^{ k } Q^k u_k \,, \ u_k \in H^s_h ( X ) \}&  m \leq
 0  \end{array}\right.
\] 
with the norms equivalent to the same natural way as for Sobolev
spaces:
\[ \| u \|_{ H_\Sigma^{s, m } ( X) } {\simeq} \left\{
\begin{array}{ll}
\sum_{k=0}^m \| (\tilde h/h )^k Q^k u \|_{H^s_h ( X )} \,,  &
 m \geq 0 \\
\ & \ \\
\inf\{  \sum_{k=0}^{ |m| } \| u_k \|_{ H^s_h (  X) } \; : \; 
u =  \sum_{k=0}^{ |m| } ( \tilde h /h )^{ k } Q^k u_k  \}\,, &  m \leq
 0  \end{array}\right.
\] 

 This can be seen using a spectral decomposition of $ Q $ which is
assumed to be bounded and self-adjoint.
We use the following simplified notation
\[ H^m_\Sigma ( X ) \stackrel{\rm{def}}{=} H^{0,m}_\Sigma ( X ) 
\,, \ \ m \in \ZZ \,.\]
The spaces $ H_\Sigma^{s, m } ( X ) $
have the following basic invariance property:
\begin{lem}
\label{l:invh}
The definition \eqref{eq:defh} does not depend on the choice
of $ Q \in \Psi_h^{0,0} ( X ) $ satisfying \eqref{eq:pq}. If 
$ A \in \Psi^{0,0} ( X ) $ has the property that $ d ( \WFh ( A ) , \Sigma ) 
> 1 / C $ then for $ M \geq 0 $
\begin{equation}
\label{eq:l.invh}
\begin{split}
&  \| A u \|_{ H_\Sigma^{-M} ( X ) }\leq ( h/\tilde h )^{M} \| u \|_{L^2 ( X)}
 \,, \\  
& \| A u \|_{ L^2 ( X ) } \leq ( h / \tilde h )^{M } \| u \|_{
H^{M}_\Sigma ( X ) } \,, \end{split} \end{equation}
for $ u \in \CIc ( X ) $.

Also, suppose that  $ F $ is a $0$th order $h$-Fourier Integral operator 
 associated to a canonical transformation which 
maps $ \Sigma $ to another hypersurface $ \Sigma' $ satisfying our
hypothesis. Then
\[ F \; : \; H_\Sigma^{s, m }( X )  \ \longrightarrow
H_{\Sigma'}^{s,m} ( X ) \,.\]
\end{lem}
\begin{proof}
Let $ Q' $ be another operator satisfying \eqref{eq:pq}. 
Then 
$$ Q' = A Q + E = Q A + E' \,, $$ 
where $ A \in \Psi^{0,0} ( X) $ is uniformly 
elliptic, and $ E, E' \in \Psi^{-1,-1} ( X) $. Because of the
interpolation property of $ H^{s, m}_\Sigma $ we only need to 
check the independence for $ k \in \ZZ $. 
Then, by induction on $ k $, 
\[ \|   ( \tilde h/ h )^k (Q')^k u \|_{ H^s ( X ) } \leq 
C_k \sum_{ \ell =0}^k \tilde h^{ k - \ell} \| ( \tilde h / h )^\ell
Q^\ell \|_{ H^s ( X ) } \,, \] 
where $ C_k $ depends only on $ A $ and $ E $. This
shows that the definition for $ m \geq 0 $ is independent of
the choice of $ Q $. The case of $ m < 0 $ is similar. To 
see the first inequality in \eqref{eq:l.invh} we recall that
$$ \| A u \|_{H^{-M}_\Sigma ( X ) } = \inf \left\{ 
\sum_{k=0}^M \| u_k \|_{ L^2 ( X ) } \; : \; 
A u = \sum_{k=0}^M ( \tilde h / h )^k Q ^k u_k \right\} \,. $$
Because of the ellipticity of $ Q $ on $\WFh ( A ) $ 
we can find $ B $ such that 
$$ A u = Q^k B A u + v \,, \ \  \| v\|_{L^2 } = {\mathcal O} (
h^\infty ) \| u \|_{ L^2 }\,. $$
Hence we can take $ u_0 = v $ and $ u_M = (  h/ \tilde h)^M B A u $,
$ u_k = 0 $ for all other $ k$'s, so that
$$ \| A u \|_{H^{-M}_\Sigma ( X ) } \leq C ( h / \tilde h )^M 
\| u \|_{L^2 ( X ) } + {\mathcal O} (h^\infty ) \| u \|_{ L^2 ( X ) } $$
proving the first estimate in \eqref{eq:l.invh} ($ h < \tilde h $
everywhere here). Since $ H^M_\Sigma ( X ) $ is the dual of 
$ H^{-M}_\Sigma $ this is equivalent to 
\[ \forall \; v \in \CIc ( X ) \ \ 
\langle A u , v \rangle_{L^2 ( X) } \leq C ( h / \tilde h)^M 
\| u \|_{ L^2 ( X)  }
\| v \|_{ H^M _\Sigma ( X)  } \,, \]
which in turn proves the second estimate with $ A $ replaced by 
$ A^* $.

The Egorov theorem (Proposition \ref{p:egorov} above) shows that
$ Q F $ is equal to $ F Q' + E $ where $ Q'$ satisfies
\eqref{eq:pq} with $ \Sigma $ replaced by $ \Sigma' $. 
The mapping property follows from the boundedness of 
$ F $ on $ H^s ( X ) $ and an argument similar to that 
above.
\end{proof}

\subsection{Globally defined class of operators}
\label{gc}
For the global definition we cannot use the classes 
\[  \Opt\left(\widetilde{\mathcal O}_{\frac12}( \langle \lambda \rangle^m ) \right) \]
due to  the presence of propagation in the $h$-sense.
To see the problem let us consider the operator $ P $
introduced in \S \ref{in} and define
\[  A = h^{-m} \chi \left( \frac{ \tilde h P }{ h } \right) \,, \ \ 
\chi \in \CIc ( \RR ) \,.\]
In the model theory of \S \ref{cmc} this operator is obtained
from the symbol $h^{-m} \chi ( \lambda ) $. However, when  
$ \tilde h $ is fixed, the operator $ A $ is
an $ h$-Fourier integral operator:
\[ A = \frac{h^{-m} 
}{ 2 \pi \tilde h } \int_\RR \widehat \chi ( t / \tilde h 
) \exp ( i t P / h ) dt \,,\]
associated to the relation
\[ \{ ( \rho, \rho' ) \; : \; 
\exists \; t \,, \  \exp ( t H_p ) ( \rho ) = \rho' \,, \ p ( \rho ) = 
p ( \rho' ) = 0 \} \,.\]
Hence the presence of almost closed orbits of the flow 
prevents a definition which would be purely local in the 
$ {\mathcal O} ( h ^\infty ) $ sense. We observe however that
the non-local contributions in $ A $ are of the order $ {\mathcal 
O } ( \tilde h ^\infty )$.

These global concerns suggest the
following definition of the residual class. First we introduce 
a useful cut-off operator. Let $ V_1, V_2 $ be two open neighbourhoods of 
$ \Sigma $, satisfying 
$$ \overline{V_1} \Subset V_2 \,.$$ 
Then we choose
\begin{equation}
\label{eq:ls}
 \gamma_\Sigma \in \Psi^{0,0}( X) \,, \ \ \WFh ( \gamma_\Sigma ) 
\subset V_2 \,, \ \ \WFh ( I - \gamma_\Sigma ) \subset \complement V_1 
\,.\end{equation}
We now define the spaces of  operators.
Let 
\[ {\mathfrak p} = ( m, \tilde m , k_1, k_2 ) \,, \ \
{\mathfrak p}_\infty = ( m, - \infty , - \infty , k_2) \,.\]
As before we start with the definition of a residual class:

\medskip
\noindent
{\bf Definition 1.}
We say that $ A \in \Psi_{\Sigma, 
\frac12}^{{\mathfrak p}_\infty}( X ) $ 
 if 
\[ A : \CIc ( X ) \rightarrow \CI ( X ) \,, \ \ 
 ( I - \gamma_\Sigma ) A\,, \, \,  A  ( I - \gamma_\Sigma )  
 \in \Psi_{\frac12}^{m+k_2,-\infty,- \infty } ( X ) \,, \]
and for 
any $ u \in \CIc ( X ) $, 
any sequence $ \{ b_j \}_{ j=1}^N 
\subset S^0 ( T^* X ) $,
and any $ k $, $ p $, and  $M $,
 we have
\[
\| \ad_{\Op (b_1)} \circ \cdots \ad_{\Op (b_N ) }
\ad^k_Q  A u \|_{ 
H_\Sigma^{ p + N/2 -k_2 }( X )  } \leq C h ^{-m +k} 
\tilde h^M \| u \|_{ H_\Sigma^{p + k/2}  }  \,,
\]
where $ Q $ is the operator in \eqref{eq:defh}.

\medskip
\noindent
{\bf Definition 2.}
 We say that $ A \in \Psi^{ \mathfrak p }_{\Sigma , \frac12}
( X ) $ if 
\begin{gather*}
 A : \CIc ( X ) \rightarrow \CI ( X ) \,, \ \ ( I - \gamma_\Sigma ) 
A \,, \ A  ( I - \gamma_\Sigma )  \in 
\Psi_{\frac12}^{m +k_2,\tilde m - k_2 , k_1} ( X ) \,, \end{gather*}
for some cut-off operator $ \gamma_\Sigma $ satisfying \eqref{eq:ls}.
\begin{itemize}
\item For any $  \chi \in \CIc ( T^* X)  $, 
\[  A = A_\chi^\sharp + A_\chi^\flat \,, \ \ A_\chi^\flat \in \Psi_{\Sigma, 
\frac12}^{{\mathfrak p}_\infty}( X ) \,, \] 
so that for any $ p$,
and any sequences
\[ \{ b_j \}_{ j=1}^N \,,  \ \
\{ a_j  \}_{j=1}^M \subset S^0 ( T^* X )\,, \ 
 H_p a_j \rest_{ {\supp \chi }  } \equiv 0 \,,\]
\begin{gather*}  
\| \ad_{\Op (b_1)} \circ \cdots \ad_{\Op (b_N ) }
\circ \ad_{\Op (a_1)} \circ \cdots \ad_{\Op (a_M ) }
\circ \ad_Q^k \chi^w \gamma_\Sigma A_\chi^\sharp \chi^w u \|_{ 
H_\Sigma^{p+N/2-k_2}( X )  } \\ \ \leq \; C 
h^{ M/2 + k - m  } \tilde h ^{M/2 + N - {\widetilde m} }  \| u \|_{H^{p+k/2}_
\Sigma  ( X) } \,, \ \ \ \chi^w = \Op ( \chi)\,,
\end{gather*}
for any $ u \in \CIc( X ) $.
\item
For any $ \psi_1, \psi_2 \in \CI_{\rm{b}} ( T^* X )  $, $
\supp \psi_1 \cap \supp \psi_2 = \emptyset $,
$$ \Op ( \psi_1 ) A \Op ( \psi_2)  \in \Psi_{\Sigma, \frac12}^{
{\mathfrak p}_\infty} ( X )  \,.$$
\end{itemize}

It is important to record
\begin{lem}
\label{l:import}
Definitions 1 and 2 are independent of the choice of the cut-off
operator satisfying \eqref{eq:ls}.
\end{lem}
\begin{proof}
Suppose that $ \gamma_\Sigma' $ is
another cut-off operator satisfying \eqref{eq:ls}. We need to show that
\begin{equation}
\label{eq:smg}
 ( \gamma_\Sigma - \gamma_\Sigma' ) A \,, \ \ A ( \gamma_\Sigma - 
\gamma_\Sigma' ) \in \Psi_{\frac12}^{m + k_2 , \tilde m - k_2 
, - \infty } ( X ) \,, 
\end{equation}
and for that we will use only the commutators involving $ b_j$'s 
in Definitions 1 and 2.
The first inclusion is, in view of Lemma \ref{l:1hf} and the fact that
$ \gamma_\Sigma - \gamma_\Sigma' $ is smoothing, equivalent to 
\begin{gather}
\label{eq:equiv}
\begin{gathered}
 \| \ad_{\Op (b_1)} \circ \cdots \circ \ad_{\Op (b_N ) }
( \gamma_\Sigma - \gamma_\Sigma' ) A u \|_{L^2( X )  } \\ 
\ \ \ \ \ \ \leq C 
h ^{-m -k_2 + N/2 } 
\tilde h^{ - \tilde m + k_2 + N/2 }  \| u \|_{ L^2 ( X )   }  \,.
\end{gathered} \end{gather}
To see this we first use  
\eqref{eq:l.invh} to see that if $ E \in \Psi^{0,0} ( X ) $ and 
$ d(\WFh ( E ) , \Sigma ) > 1/C 
 $ then for any subsequence $ \{i_j\}_{j=1}^J $ of $ 1, \cdots, N $,
\begin{equation}
\label{eq:usee}
\begin{split}
&  \| E  \ad_{\Op (b_{i_1} )} \circ \cdots \circ \ad_{\Op (b_{i_J} ) } 
A u \|_{ L^2 ( X ) } \\
& \ \ \ \ \  \leq C ( h/ \tilde h)^{ J/2 - k_2 }  \| 
  \ad_{\Op (b_{i_1})} \circ \cdots \circ \ad_{\Op (b_{i_J} ) } 
A u \|_{ H_\Sigma^{J/2 -k_2  } ( X ) } \\
& \ \ \ \ \  \leq 
C h^{-m} ( h / \tilde h)^{ J/2  - k_2 } 
\tilde h^{ J - \tilde m } \| u \|_{L^2 ( X ) } = 
C h^{-m  - k_2 + J/2 } 
\tilde h^{ - \tilde m + k_2 + J/2  } \| u \|_{L^2 ( X ) }
\,.
\end{split} \end{equation}
Using the derivation property of $ \ad_{ \Op ( b_j ) } $ and the fact
that 
\[  \ad_{\Op ( b_j ) } ( \gamma_\Sigma - \gamma_\Sigma' ) = h E_j \,, 
\ \ E_j \in \Psi^{0,0} ( X ) 
\,, \ \ d(\WFh ( E_j ), \Sigma ) > 1/C\,, \ \ \]
we can estimate the left hand side of \eqref{eq:equiv} 
by a linear combination of terms of the form
\[  h^{N-J} \| E  \ad_{\Op (b_{i_1} )} \circ \cdots \circ 
\ad_{\Op (b_{i_J} ) } A u \|_{ L^2 ( X ) } \,, \ \ E \in \Psi^{0,0}( X ) \,,
\ \ d (\WFh ( E ) , \Sigma ) > 1/C \,. \]
Consequently \eqref{eq:equiv} follows from \eqref{eq:usee}.

Since, by the calculus of \S \ref{sma}
$$ B \in \Psi_{\frac12}^{m + k_2  , - \infty , -\infty } ( X ) 
\ \Longleftrightarrow \ 
B^* \in \Psi_{\frac12}^{m + k_2  , - \infty , -\infty } ( X )  \,, $$
and since 
$$  \ad_{\Op (b_1)} \circ \cdots \circ \ad_{\Op (b_N ) } B^* =  
( -1)^N \left( \ad_{\Op (\bar b_1)} \circ \cdots \circ \ad_{\Op (\bar b_N ) } 
B \right)^* \,,$$   
the second inclusion in \eqref{eq:smg} follows from the first one.
\end{proof}

We now have a natural mapping and invariance properties:
\begin{prop}
\label{p:invr}
The operators $ \Psi^{0,0,0,0}_{\Sigma, \frac12} ( X) $ form an algebra
and $ \Psi^{0, - \infty,  - \infty, 0  }_{\Sigma, \frac12}
 ( X ) $ is an ideal in 
that algebra. If $ A 
\in \Psi^{{\mathfrak p}_\infty }_{\Sigma, \frac12} ( X )  $ then 
\[ A = {\mathcal O} ( h^{-m} \tilde h^M ) \; : \; H^{-M,r}_\Sigma ( X ) 
\rightarrow H^{M,r - k_2 }_\Sigma (  X) \,, \]
for any $ m, M $. If $ A \in \Psi_{\Sigma, \frac12}^{{\mathfrak p}  } 
( X) $ then 
\[ A = {\mathcal O} ( h^{-m  } \tilde h^{-{\widetilde m}} ) \; 
: \; H^{s+k_1, p+k_2}_\Sigma ( X ) \longrightarrow H^{s,p} _\Sigma ( X ) \,,
\]
for any $ s \in \RR $ and $ p \in \ZZ $.

Also, suppose that $  F , G $ are 
$0$th order $h$-Fourier Integral operators, $ F $ 
associated to a canonical transformation which 
maps $ \Sigma $ to another hypersurface $ \Sigma' $ satisfying our
hypothesis, and $ G $ to its inverse (both transformation need to be
defined only locally). Then 
\begin{equation}
\label{eq:p.inv1}
 A \in   \Psi_{\Sigma, \frac12} ^{\mathfrak p} ( X )  
\ \Longrightarrow \ F \circ A \circ G \in  
 \Psi_{\Sigma', \frac12}^{\mathfrak p} ( X )  \,.\end{equation}
\end{prop}
\begin{proof}
We first show that $ \Psi_{ \Sigma, \frac12}^{0,-\infty,-\infty,0} ( X) $
is an algebra. In fact, in the notation of Definition 1,
\[ \| \ad_{\Op (b_1)} \circ \cdots \ad_{\Op (b_N ) }
\ad^k_Q  A B u \|_{ 
H_\Sigma^{ p + N/2 }( X )  } \,,  \ \ 
A , B \in  \Psi_{ \Sigma, \frac12}^{0,-\infty,-\infty,0} ( X) \,,\]
is bounded by a linear combination of 
\begin{gather*} 
\| \ad_{\Op (b_{i_1})} \circ \cdots \ad_{\Op (b_{i_J} ) }
\ad^{l}_Q  A 
 \ad_{\Op (b_{k_1})} \circ \cdots \ad_{\Op (b_{i_{N-J}} ) }
\ad^{k-l}_Q B u \|_{ 
H_\Sigma^{ p + N/2 }( X )  } \,,  \\ 
\{ i_j\}_{j=1}^J \cup \{ k_j \}_{j=1}^{N-J} = \{ 1, \cdots , N \} \,.
\end{gather*}
These terms are estimated by 
\begin{gather*}
 C h^{l} \tilde h^{M_1} \| \ad_{\Op (b_{k_1})} \circ \cdots \ad_{\Op (b_{i_{N-J}} ) }
\ad^{k-l}_Q B u \|_{ 
H_\Sigma^{ p + (N-J)/2 + l/2} ( X )  } \\
\ \ \ \ \ \ \ \  \ 
\leq C' h^{k} \tilde h^{M_1 + M_2 } \| u \|_{ p + N/2 + k } \,, 
\end{gather*}
where we used Definition 1 with  $ A $ and with $ B $. Here $ M_1 $ and
$ M_2 $ are any large integers. Checking that 
\[ ( I - \gamma_\Sigma ) A B \,, \ A B ( I - \gamma_\Sigma ) \in 
\Psi_{\frac12}^{0,-\infty,-\infty} ( X ) \,, \]
is the same as in the proof of Lemma \ref{l:import}.

We now check that for $ B \in \Psi^{0,-\infty, -\infty, 0 }_{\Sigma, \frac12}
( X) $ and 
 $ A \in \Psi^{0,0, 0, 0 }_{\Sigma, \frac12} ( X) $ we have 
\[ B A \in \Psi^{0,-\infty, -\infty, 0 }_{\Sigma, \frac12} ( X) \,.\]
In the notation of Definition 2, we can write 
\[ B = \gamma_\Sigma B_1 \gamma_\Sigma + ( 1 - \gamma_\Sigma ) B_1 
 +  \gamma_\Sigma B_1 ( 1 - \gamma_\Sigma ) + B_2 \,, \]
and it suffices to check that 
\begin{equation}
\label{eq:ver}
 A \gamma_\Sigma B_1 \gamma_\Sigma \in 
\Psi_{\frac12}^{0,-\infty,-\infty} ( X ) \,. \end{equation}
Since $ \Sigma $ is compact (and hence we can take 
$ \chi^w = \gamma_\Sigma $) and Definition 2 is independent of the 
choice of $ \gamma_\Sigma $ we have, as its special case,
\[ 
\| \ad_{\Op (b_1)} \circ \cdots \ad_{\Op (b_N ) }
\circ \ad_Q^k \gamma_\Sigma  B_1 \gamma_\Sigma  u \|_{ 
H_\Sigma^{p+N/2 }( X )  } \\ \ \leq \; C 
h^{  k   } \tilde h ^{ N }  \| u \|_{H^{p+k/2}_\Sigma  ( X) } \,,
\]
and the verification of \eqref{eq:ver} follows the 
proof of the algebra property of $ \Psi_{\Sigma, \frac12}^{0,-\infty, 
-\infty, 0 } ( X ) $.

To conclude the algebraic part of the proof we show that 
$ \Psi_{\Sigma, \frac12}^{0,0,0,0} ( X ) $ is closed under composition 
of operators.
Let $ B_1 , B_2 \in  \Psi_{\Sigma, \frac12}^{0,0,0,0} ( X ) $ and 
let $  \chi \in \CIc ( T^* X)  $.
Since we already established that composition with 
elements of $ \Psi_{\Sigma, \frac12}^{0,-\infty,
-\infty, 0} (X ) $ produces an operator in that space we only have to 
show that 
\[ B_1 B_2 = A_1 + A_2 \,, \ \ A_2 \in \Psi^{0,-\infty, - \infty, 0 }
_{ \Sigma , \frac12} ( X) \,,\]
so that 
for any $ p$,
and any sequences
\[ \{ b_j \}_{ j=1}^N \,,  \ \
\{ a_j  \}_{j=1}^M \subset S^0 ( T^* X )\,, \ 
 H_p a_j \rest_{ {\supp \chi }  } \equiv 0 \,,\]
and
\begin{gather*}  
\| \ad_{\Op (b_1)} \circ \cdots \ad_{\Op (b_N ) }
\circ \ad_{\Op (a_1)} \circ \cdots \ad_{\Op (a_M ) }
\circ \ad_Q^k \chi^w \gamma_\Sigma  A_1  \chi^w  u \|_{ 
H_\Sigma^{p+N/2-k_2}( X )  } \\ \ \leq \; C 
h^{ M/2 + k - m  } \tilde h ^{M/2 + N - {\widetilde m} }  \| u \|_{H^{p+k/2}_
\Sigma  ( X) } 
\end{gather*} 
To find the decomposition of $ B_1 B_2 $ we introduce $ \chi_j \in 
\CIc ( T^* X) $ such that $ \chi_1 \equiv 1 $ on $ \supp \chi_0 $,
and $ \chi_0 \equiv 1 $ on $ \supp \chi $. We choose $ \supp \chi_1 $
sufficiently close to the support of $ \chi $ so that the 
functions $ a_j $ satisfy $ H_p a_j \rest_{\supp \chi_1 } \equiv 0 $.
We then put 
\[  A_2 = \chi_0^w (B_1)^\sharp_{\chi_1} 
 ( 1 - \chi_1^w ) (B_2)^\sharp_{\chi_1}
+ B_1 (B_2)^\flat_{\chi_1}  +  (B_1)^\flat_{\chi_1} B_2 \,,\]
which is in $ \Psi_{\Sigma , \frac12}^{0,-\infty,-\infty,0} ( X) $ since
$ \supp \chi_0 \cap \supp ( 1 - \chi_1^w ) = \emptyset $, and 
we can use the second part of Definition 2. 

We have 
\[  A_1 = ( 1 - \chi_0^w )  (B_1)^\sharp_{\chi_1}
 (B_2)^\sharp_{\chi_1}   + \chi_0^w  (B_1)^\sharp_{\chi_1}  \chi_1^w  
 (B_1)^\sharp_{\chi_1}  \,.\]
 Up to negligible, $ {\mathcal O} ( h^\infty ) $, errors 
\[ \chi^w ( 1- \chi_0^w ) \equiv 0 \,, \ \ 
\chi^w \chi_0^w \equiv \chi^w \,.\]
Hence we need to check that 
\begin{gather*}
 \ad_{\Op (b_1)} \circ \cdots \ad_{\Op (b_N ) }
\circ \ad_{\Op (a_1)} \circ \cdots \ad_{\Op (a_M ) }
\circ \ad_Q^k \chi^w (B_1)^\sharp_{\chi_1} \chi_0^w (B_2)^\sharp_{\chi_1}
  \chi^w  \\
 \ \ \ = \; {\mathcal O}  
( h^{ M/2 + k - m  } \tilde h ^{M/2 + N - {\widetilde m} })
\; : \; H_\Sigma^{p+N/2-k_2}( X ) \; \longrightarrow \; {H^{p+k/2}_
\Sigma  ( X) } \,,
\end{gather*}
and this follows from the Leibnitz rule for $ \ad_{b_j} $ and 
assumptions on $ (B_j)_{\chi_1}^\sharp $.

The mapping properties are immediate from the definitions: we apply
them with no commutators, and in the case of Definition 2, with 
$ \chi^w = \gamma_\Sigma$.

Lemma \ref{l:invh} 
shows that the spaces $ H_\Sigma^{\pm M } $ transform correctly under
$ F $ and hence Proposition \ref{p:egorov} (Egorov's Theorem) 
shows that for 
\[  A \in \Psi^{{\mathfrak p}_\infty}_{\Sigma , \frac12 } 
( X ) \ \Longrightarrow \ 
 F \circ A \circ G \in \Psi^{{\mathfrak p}_\infty}  _{\Sigma',\frac12} 
( X) \,, \]  
that is \eqref{eq:p.inv1} holds for the residual class.
Since the conditions on $ a_j $, $ b_j$'s in Definition 2 are 
symplectically invariant, Egorov's theorem (Proposition \ref{p:egorov})
again shows that 
$ F \circ A \circ G \in \Psi^{\mathfrak p }_{\Sigma'} ( X ) $.
\end{proof} 
 
\subsection{The symbol map}
\label{sm}
To define the symbol map we will use the invariance given by Proposition 
\ref{p:invr} and the symbol map in the model case. We start with 
\begin{lem}
\label{l:symb1}
In the notation of Definition 2, suppose that $ X = T^* \RR^n $ and that
\[   \Sigma\cap V = \{ \xi_1 = 0 \} \cap V \,, \ \ V \subset T^* \RR^n \,, 
\ \text{ is open,  and }  \supp \chi \Subset V \,. \] 
Then 
\[ \Op ( \chi) A_1 \Op ( \chi) = \Opt ( a_\chi ) \,, \ \ 
a_\chi =  \widetilde{\mathcal O }_{\frac12} ( h^{-m} \tilde h ^{-\widetilde m}
\langle \lambda \rangle^{k_2 } ) \,. \]
\end{lem}
\begin{proof}
This is a consequence of Proposition \ref{p:Beals} and the properties
of the term $ A_1 $ in Definition 2.
\end{proof}

To construct a symbol map, that is a  homomorphism 
\begin{gather*}
  \sigs \; : \; \Pkh ( X ) \longrightarrow \Ssh ( T^* X )  
/ \Ssgh   ( T^*X )   \,, \\
 {\mathfrak p} = ( m, \tilde m , k_1, k_2 ) \,, \ \
{\mathfrak p}' = ( m, \tilde m-1  ,k_1 - 1 , k_2) \,, 
\end{gather*}
such that the sequence   
\begin{gather*}
0 \longrightarrow \Pkgh ( X ) \longrightarrow \Pkh ( X)  
\stackrel{\sigs}{\longrightarrow}   \Ssh ( T^* X )  
/  \Ssgh   ( T^* X )  \longrightarrow 0  
\end{gather*}
is exact, for an 
arbitrary $ \Sigma $ we will use Lemma \ref{l:symb1}. 
That requires putting $ \Sigma $ locally to the
model hypersurface $ \xi_1 = $ which on the quantum level 
is done using $h$-Fourier integral operators. Hence we need
a local invariance statement given in the next
\begin{prop}
\label{th:3.2} 
Let $ U $ be an  $h$-Fourier Integral Operator,
elliptic in $ V \times V $, where $ V $ is a
neighbourhood  of $  ( 0 , 0 ) \in T^* \RR^n  $, the compact hypersurface
$ \Sigma $ satisfies,
$$ \Sigma \cap V = \{ \xi_1 = 0 \} \cap V \,.$$
Assume also that the canonical transformation associated to $ U $, $\kappa $,  
satisfies: 
\begin{equation} 
\label{eq:th.3} \kappa ( 0, 0 ) = ( 0 , 0 ) \,, \ \  
\kappa ( \{ \xi_1 = 0 \} \cap V ) \subset \{ \xi_1 = 0 \}  \,.
\end{equation} 
Let $ A = \Opt ( a ) $ where  
$  a = {\widetilde{\mathcal O}_{\frac12}}( h^{-m} \tilde h^{-\widetilde m } 
\langle \lambda \rangle^{k} ) \,, \ \
a ( x , \xi , \lambda ) \equiv 0 \ \text{ for $ ( x , \xi ) \notin V $.} $

If $ U^{-1}
$ is the microlocal inverse of $ U $ near $ V \times V  $,
then 
\begin{gather*}
 U^{-1} \circ A \circ U = \Opt ( b ) + E \,, 
\ \  b =  a \circ K  \,, \\ 
E \in \Psi_\Sigma^{{\mathfrak p}'} 
(  \RR^n )  \,, 
\ \ {\mathfrak p}' = ( m ,  \widetilde m -1 , - \infty , k - 1 ) \,, 
\end{gather*}
where $ K $ is the natural 
lift of $ \kappa $ to the $ ( x, \xi, \lambda ) $ variables: 
\begin{equation}
\label{eq:Kh}
 K ( y , \eta , \mu ) = ( x, \xi , \lambda ) 
\ \Longleftrightarrow \
 ( x, \xi  ) = \kappa ( y , \eta ) \,, \ \ 
\lambda =  ( \xi_1 / \eta_1 )\mu \,. \end{equation}
\end{prop} 
\begin{proof} 
We start by observing that the proposition holds in the special cases 
$ \kappa ( x, \xi ) = ( x, \xi ) $ and $ \kappa ( x, \xi) = ( - x,  
-\xi ) $. The first special case concerns conjugation with elliptic 
classical $h$-pseudodifferential operators and it follows from the 
discussion after \eqref{eq:3.5}. The second special case  
follows from the first one and the fact that the proposition is easily  
checked for $ U u ( x  ) = u ( - x ) $. As a consequence we can  
asssume that $ \kappa $ preserves the sign of $ \xi_1$: 
\begin{equation} 
\label{eq:3.9s} 
\kappa ( y , \eta ) = ( x, \xi ) \ \Longrightarrow \ 
\xi_1 \eta_1 \geq 0 \,. 
\end{equation} 
We will prove the proposition by a deformation method inspired 
by the ``Heisenberg picture of quantum mechanics''  
and for that we need the following 
geometric 
\begin{lem} 
\label{l:3.2s} 
Let $ \kappa $ be a smooth canonical transformation satisfying 
\eqref{eq:th.3} and \eqref{eq:3.9s}. Then we can find a  
piecewise smooth family of canonical transformations $ 
[0,1] \ni t \mapsto \kappa_t $  satisfying \eqref{eq:th.3}, 
\eqref{eq:3.9s} and such that $ \kappa_0 = id $ and  
$ \kappa_1 = \kappa $.  
\end{lem} 
\begin{proof} 
Let us denote by $ \Sigma $ the hypersurface given by  
$ \xi_1  = 0 $. We first observe that if $ \kappa $ is a linear
symplectic transformation preserving $ \Sigma $ then we
can find a family of linear symplectic transformations,
$ \kappa_t $, satisfying the conclusions of the lemma: the subgroup
of elements of $ Sp( n, \RR ) $ preserving a half space 
bounded by $ \Sigma $ is connected. 

Hence we can assume that $ d \kappa ( 0 , 0 ) = Id $. Now 
introduce
\[  \kappa_t ( x , \xi ) = t^{-1} \kappa \left( tx , t \xi  
\right) \,,\]
a smooth family of symplectic transformations, preserving a
half space bounded by $ \Sigma $, with $ \kappa_1 = \kappa $ and
$ \kappa_0 = Id $.
\end{proof} 

We now return to the proof of Proposition \ref{th:3.2}. 
For simplicity we can assume that $  m = \tilde m = 0 $.
Let  
$ \kappa_ t $ be a piecewise smooth family with $  \kappa_1 =  
\kappa $, $ \kappa_0 = id $. Using 
\eqref{eq:2.2} we construct  a piecewise  
smooth family of classical elliptic $h$-Fourier integral  
operators,  $ U_t $, defined microlocally near $(0,0)$ and associated to  
$ \kappa_t $. If we demand that 
 $ U_1 = U $ then $ U_0  $ is a pseudodifferential operator 
elliptic at $ ( 0 , 0 ) $. 
 For notational convenience we assume that our deformation 
is smooth in $ t $ -- the piecewise smooth case follows from 
the same argument applied in  several steps and that $ U_0 = Id $ 
(the last condition can be arranged by composing $ U $ with a 
elliptic pseudodifferential operator). 
Thus we have 
\begin{equation} 
\label{eq:3.17s} 
h D_t  U_t  + U_t Q_t = 0 \,, 
\end{equation} 
where $ Q_t $ is a smooth family of classical $ h$-pseudodifferential 
operators of order $0$ with the leading symbol $ q_t$ satisfying 
\begin{equation} 
\label{eq:3.18s} 
\frac{d}{dt} \kappa_t ( x, \xi) = (\kappa_t)_* \left(  
H_{q_t} ( x, \xi ) \right) \,, 
\end{equation} 
in a neighbourhood of $ V $. 
It follows from \eqref{eq:3.18s} that $ H_{q_t } $ is tangent 
to $ \Sigma $ and hence 
\begin{equation} 
\label{eq:3.19s} 
\partial_{x_1} q _t ( x, \xi ) = \xi_1 r_t ( x, \xi )  \,.  
\end{equation} 
We extend $ q_t $ to a globally defined function in 
$ S ( T^* \RR^n,  \langle x , \xi \rangle) $, 
keeping the property \eqref{eq:3.19s}. 
This defines a family of global canonical transformation which coincide
with $ \kappa_t $ near $ V $.

Let $ V_t $ satisfy 
\begin{equation} 
\label{eq:3.20s} 
hD_t V_t = Q_t V_t \,, \ \ V_0 = Id \,. 
\end{equation} 
It follows that $ V_t = U_t^{-1} $, and if we take $ Q_t $ to be 
self-adjoint, $ V_t^* = U_t $.

We will now use \eqref{eq:3.17s} and \eqref{eq:3.18s}
to prove Egorov's theorem \eqref{eq:egorov} for the new class.
Thus we consider $ A_t  = U_t^{-1} A U_t $, so that $ A_1 $ is the  
operator we want to study and $  A_0 = A $  is the given operator. 
From Proposition \ref{p:invr} we already know that $ A_t \in 
\Psi_{\Sigma, \frac12}^{0,0,0,k} ( \RR^n ) $. Using Lemma 
\ref{l:symb1} we can write 
\[  A_t \equiv \Opt ( a_t ) \mod 
\Psi_{\Sigma, \frac12}^{0,-\infty,-\infty,k} ( \RR^n ) \,, \]
microlocally near $ V $. 
From \eqref{eq:3.17s} and \eqref{eq:3.20s}, we get 
\begin{gather} 
\label{eq:3.21s} 
\begin{gathered}
\left\{ \begin{array}{l}  
hD_t  A_t  = [ Q_t , A_t ]  \,, \\ 
A_0 = A \,, \end{array} \right. 
\end{gathered}
\end{gather} 
To compute the commutator on the symbolic level we need:
\begin{lem}
\label{l:com}
Suppose that  $ a \in 
\widetilde {\mathcal O}_{\frac12} ( \langle \lambda \rangle^k ) $  
and 
that $ b \in S^{0,-\infty} ( T^*\RR^n ) $
 (that is, $ b $ is a symbol in the
sense of \S \ref{rsa}), does
not depend on $ \lambda $. In addition let us assume that 
\[  \partial_{x_1} b = \xi_1 r ( x , \xi ) \,. \]
Then
\begin{equation}
\label{eq:fcom}
\frac{i}{h} [ \Opt ( b ) , \Opt ( a) ] = \Opt (c ) \,, \ \ 
c = (H_b  - r \lambda \partial_\lambda ) a + \widetilde {\mathcal O}_{\frac12}
 (  \tilde h   \langle \lambda \rangle^{k}  )
\,. 
\end{equation}
\end{lem}
\begin{proof}
We will use \eqref{eq:tensq} to compute $ a \sharp_{h, \tilde h} b $
(and $ b \sharp_{h , \tilde h } a $), noting that
\[ \Oph a ( x_1, \tilde h D_{x_1} ) \circ
\Oph b ( x_1 , \tilde h D_{x_1} ) = 
\left( \Oph (a ) \sharp_{\tilde h} \Oph (b ) \right) ( x_1, \tilde h D_{x_1} ) 
\,.\]
Since $ b \in S^{0, -\infty } $,  Lemma \ref{l:Sjnew} shows that 
the composition formula in the $ ( x' , \xi' ) $ is given by 
an asymptotic series in $ ( h \tilde h)^{\frac12} $, and a $ {\mathcal O} ( h^\infty ) $ error. The only subtlety lies in the dependence on 
$ ( x_1 , \lambda ) $ and to explain it we suppress the other variables.
We have 
\[  \partial_{x_1}^\ell \partial_\lambda^p b ( x_1, \lambda ) 
= \left\{ \begin{array}{ll} {\mathcal O} ( ( h / \tilde h)^{p} \langle ( h/ \tilde h) 
\lambda \rangle^{-\infty} ) & \ell = 0 \\
{\mathcal O} ( ( h / \tilde h)^{p+1} \lambda \langle ( h/ \tilde h) 
\lambda \rangle^{-\infty} ) & \ell > 0 \,, \end{array}\right. 
\] 
and $ \partial_{x_1}^p \partial_{\lambda}^\ell a = {\mathcal O} ( 
\langle\lambda \rangle^{k-\ell/2 + p/2 } ) $. Consequently the terms in the 
expansions of $ a \sharp_{\tilde h} b $ and $ b \sharp_{\tilde h} a $
 are bounded by 
\[ C_p h \tilde h^{p-1} \langle \lambda \rangle^{k-p/2+1}  \langle
( h/ \tilde h ) \lambda\rangle ^{-\infty} \,, 
\ \ 
 C_p ( h / \tilde h)^{p} 
 \tilde h^{p} \langle \lambda \rangle^{k+p/2}  \langle
( h/ \tilde h ) \lambda\rangle ^{-\infty} \,,  \ \ p> 0 \,, \]
respectively. Hence we have expansions such that for  $ p > 1 $ 
the terms are bounded by $ h \tilde h^{p-1} \langle \lambda \rangle^{k} $, and
the errors are 
 $ {\mathcal O} ( h \tilde h^\infty \langle \lambda \rangle^k ) $.
This argument shows that 
\[
 a \; \sharp_h \; b = ab + h \left( 
\sum_{i=1}^{n}\partial_{\xi_i} a \, D_{x_i }
b  +  r \lambda \partial_\lambda a   \right) +
\widetilde{\mathcal O }_{\frac12} \left(  h \tilde h \langle \lambda 
\rangle^{k } \right) \,, \]
and 
\[ b \; \sharp_h \; a = ab + h \sum_{i=1}^{n}\partial_{\xi_i} b \, D_{x_i }
a +   \widetilde{\mathcal O }_{\frac12} \left( h \tilde h \langle 
\lambda \rangle^{k} \right)
\,, \]
from which  the lemma follows.
\end{proof}
Since $ q_t $ satisfies
\eqref{eq:3.19s}, Lemma \ref{l:com} gives
\[ \frac{i}{h} 
[ Q_t , \Opt ( a_t ) ] = \Opt \left( ( H_{q_t } - r_t \lambda \partial 
_\lambda )a_t  + \widetilde {\mathcal O}_{\frac12} ( 
 h  \tilde h 
\langle \lambda \rangle^{k} ) 
\right) \,. \]
If we write $ a_t = a_t^0 + {\mathcal O}_{\frac12} ( \tilde h 
\langle \lambda \rangle^k ) $ then 
\begin{equation} 
\label{eq:3.23s} 
\frac{\partial }{\partial t } a_t^0 = \left( H_{q_t } -  
r_t ( x, \xi ) \lambda \frac{\partial}{\partial \lambda }
\right) a_t^0 \,, \ \ a_0^0 \equiv a \mod  {\mathcal O}_{\frac12} ( \tilde h 
\langle \lambda \rangle^k ) \,.  
\end{equation} 
We now note that 
\[ \left( H_{q_t } - r_t \lambda \frac{\partial}{\partial \lambda }  
\right) a_t^0 = H_{q_t} ( a_t^0 \rest_{ \lambda = (\tilde h/ h) \xi_1 } ) \,.\]
Hence, if  $ K_t $
is the transformation in $ ( x, \xi , \lambda )$-space 
corresponding to $ \kappa_t $ as in the statement of the proposition, 
it follows that
\[ a_t^0 =  a  \circ K_t    \,, \] 
that the principal symbol of $ A_t $ is  $ a \circ K_t $ and the  
proposition follows. 
\end{proof} 

We can now define the symbol map, 
 \[   \Pkh ( X ) \ni A \longmapsto \sigs ( A ) \in \Ssh ( T^* X )  
/ \Ssgh   ( T^*X )    \,.
\] 
We recall from \S \ref{rsa} and \S \ref{sma} 
that we already have the symbol map
for the $S_{\frac12} $ calculus:
\[ \Psi^{m,\widetilde m , k}_{\frac12} 
 ( X) \ni B \longmapsto \sigma_h ( A ) \in S_{\frac12}
^{m,\tilde m , k} ( T^* X) 
/ S^{m,  \widetilde m - 1 , k-1}_{\frac12} ( T^*X ) \,. \]
The definition of the symbol classes \eqref{eq:9.0} shows that
\[\begin{split}
&  ( \Ssh ( T^* X )  \cap \CI ( T^* X \setminus U_\Sigma ) 
)  / \Ssgh   ( T^*X )  =  \\
& \ \ \ \ ( S^{m+k_2, \widetilde m - k_2 ,k_1 } ( T^* X) \cap 
\CI ( T^* X \setminus U_\Sigma ) ) / S^{m + k_2, \widetilde m - k_2 
- 1 , k-1} ( T^*X ) \,,
\end{split} \]
for any open ($h$-independent) neighbourhood of $ \Sigma$, $ U_\Sigma $.
Hence we define
\[ \sigs (( I - \gamma_\Sigma) A ) \stackrel{\rm{def}}{=} \sigma_h (
( I -  \gamma_\Sigma ) A ) \,, \]
and as in the proof of Lemma \ref{l:import} we see that this 
definition is independent of the choice of $ \gamma_\Sigma $.

To define $ \sigs ( \gamma_\Sigma A ) $ we use Lemma \ref{l:symb1}.
We choose a partition of the cut-off operator, $ \gamma_\Sigma $,
\[  \sum_{j=1}^J \gamma_j^2 = \sigma_h (\gamma_\sigma) 
 \,, \ \ \gamma_j \in \CIc ( T^* X ) \,,\]
such that for each $ j $, $ \supp \gamma_j \cap \Sigma $ can be 
put into the normal form $ \Omega \cap \{ \xi_1 = 0 \} $ by a local 
canonical transformation,
\begin{gather}
\label{eq:kaj}
\begin{gathered}
 \kappa_j : \Omega \longrightarrow \Omega_j \,,  \ \ 
( 0 , 0 ) \in \Omega \subset T^* \RR^n \,, \\
 \kappa_j ( \{ \xi_1 = 0 \} \cap \Omega ) = \Sigma \cap \Omega_j \,. 
\end{gathered}
\end{gather}
We then {\em choose} elliptic $h$-Fourier Integral Operators, $U_j $,
microlocally defined in  neighbourhoods of $ \Omega \times \Omega_j $ and 
associated to $ \kappa_j$'s. By Proposition \ref{p:invr} 
\[  U_j^{-1} \Op ( \chi_j ) A \Op ( \chi_j ) U_j = \Op ( \tilde \chi_j )
\widetilde A \Op ( \tilde \chi_j )
\in \Psi_{\Sigma', \frac12 }^{\mathfrak p} ( \RR^n ) \,, \ \ \Sigma' \cap \Omega 
= \{ \xi_1 = 0 \} \,.\]
In the notation of Lemma \ref{l:symb1} we then define
\begin{equation}
\label{eq:symb2}
\sigs ( \gamma_\Sigma A ) = \sum_{ j=1}^J (\kappa_j^{-1})^* a_{ \tilde \chi_j }
\gamma_j^2 \,.\end{equation}
Proposition \ref{th:3.2} shows that this definition is independent of
the choices made here.

\subsection{Global quantization map}

From the local quantization given in \S \ref{cmc}
we can define a global map $ 
\Ops $.  Thus let $ a \in S_\Sigma^{\mathfrak p} (T^* X) $ be
a symbol in the class defined in \eqref{eq:9.0}. 

Let $ \gamma_\Sigma $ (where we will use the same letter for the symbol
and the operator) and $ V_j $'s be as in \eqref{eq:ls} 
By shrinking $ V_2 $ if necessary we can find 
a finite open cover
\[ V_2 \subset \bigcup_{j=1}^J \Omega_j \]
such that for each $ j$ there exists a canonical transformation
$ \kappa_j $ satisfying \eqref{eq:kaj}. 

Let $ \phi_j $ be a partition of unity on $ V_2 $ 
subordinate to the cover by $ \Omega_j$'s.  
Let $ a_j $ be the unique symbol of the form 
$$ a_j = a_j ( x, \xi_2 , \cdots, \xi_n, \lambda ; h ) $$
such that
\[ \left( a_j \right)_{ \lambda = (\tilde h / h )  \xi_1 }  = 
\left( \gamma_\Sigma \phi_j a \right) \circ \kappa_j  \,. \]
Using the $h$-Fourier integral operators defined after \eqref{eq:kaj}
we put 
\begin{equation}
\label{eq:3.op} 
\Ops ( a )  \stackrel{\text{def}}{=} \Op \left( ( 1 - \gamma_\Sigma
) a \right) + 
\sum_j U_j \Opt \left( 
a_j \right) U_j^{-1} \,.
\end{equation}
In view of Proposition \ref{p:invr}
we have  $ \Ops (a) \in \Psi^{\mathfrak p}_{\Sigma, \frac12} ( X) $. 

The construction of the symbol map in \S \ref{sm} shows that 
\[ \sigma_{\Sigma, h , \tilde h} ( \Ops ( a ) ) \equiv a \; \text{mod} \; 
S^{\mathfrak p' } _{\Sigma }   ( X ) 
 \,, \ \ \mathfrak p' = ( m  , {\widetilde m}-1, k_1-1, k_2 ) \,. \]
We recall from \eqref{eq:9.0} that {\em away from $\Sigma$},
\[  S^{\mathfrak p}_{\Sigma, \frac12}  \ \text{ becomes } \ 
 S_{\frac12} ^{m + k_2 ,  \widetilde m -k_2 , k_1 } \,, \]  and 
\[  S^{\mathfrak p'}_{\Sigma, \frac12}  \ \text{ becomes } \
  S_{\frac12} ^{m + k_2   ,  \widetilde m - k_2 - 1 , k_1 -1 } \,.\]
If $ \delta < 1/2 $ in \eqref{eq:9.0} then 
we have the usual filtration of the $h$-pseudodifferential
calculus: near $ \Sigma $ we only gain in $ \tilde h$ and away from 
$ \Sigma $, in $ h $. In the case of $ \delta = 1/2 $ considered here
in detail we only gain in $ \tilde h $, near and away $ \Sigma $.

This completes the proof of Theorem \ref{th:3.1} and provides  
an explicit quantization $ \Ops $.

\subsection{Approximation by finite rank operators}
\label{apr}

To estimate the number of resonances 
we will need to use approximation by finite
rank operators.

For $  a \in S^{m , {\widetilde m} , - \infty, k_2 } _{\Sigma, 
\frac12}  ( T^* X ) $ we 
need a notion of essential support. Unlike the essential support defined
in \S \ref{rsa} it now has to depend on $ h, \tilde h $. As in \cite{SjZw99},
rather than introduce an equivalence
class of families of sets,  we will say that for an 
$ (h, \tilde h) $-dependent 
family of sets $ W_{h, \tilde h} \subset T^* X  $
\begin{gather*}
  \esssupp a \subset W_{h, \tilde h}  \ \; \Longleftrightarrow   \\
 \; \ \exists \; a' \in
S^{m,\tilde m ,-\infty,
k_2 }_{\Sigma, \frac12} ( T^* X ) \,, \ \supp a' \subset W_{ h,\tilde h}
\,, \ \ 
a - a'  \in S_{\Sigma, \frac12}^{m,-{\infty}, - \infty, k_2 } 
( T^* X )   \,. \end{gather*}
We notice that 
\[ \esssupp a \subset V_{h, \tilde h}^j \, , \ j = 1,  \cdots , N 
\ 
 \; \Longrightarrow \; \ \esssupp a \subset V_{h, \tilde h} ^1 \cap \cdots 
\cap V_{h, \tilde h}  ^N   \,,  \]
so that this {\em formal notion} of
essential support behaves correctly under finite products and sums.
We can now state
\begin{prop}
\label{p:finite} 
Suppose that 
 $ a \in S_{\Sigma,\frac12}^{0 , 0, - \infty , - \infty } ( T^* X) $  
and 
$$ \esssupp a \subset  W_{h, \tilde h}  \,, $$
 where $ W_{h, \tilde h} $ satisfies
\begin{gather*}
 W_{h, \tilde h} \subset \{ ( x , \xi ) \; : \; d ( ( x , \xi ) , \Sigma\cap W_{h, \tilde h}  ) 
\leq C_1 h/ \tilde h \} \,,
\\ W_{h, \tilde h} \subset \bigcup_{k=1}^{K(h)}
\exp( [-1, 1] H_q ) ( B_k) \,, \ \ {\rm{diam}}\; B_k \leq C_1 (h/\tilde h )^{\frac12} \,,  \,. \end{gather*}
Then 
for $ 0 < h < h_0  $,
there exists 
a finite rank operator $ R(h) $ such that for 
\[ \Ops ( a ) - R (h) \in  \Psi^{0,  - \infty, - \infty, - \infty  }_{\Sigma, \frac12}
 ( X) \,,  
\ \ \text{\em rank}\; R (h) = C_2 \tilde h ^{  - n } 
K( h) 
  \,.   \] 
\end{prop} 
\begin{proof} 
We  take an open covering of $ W_{h, \tilde h} $, 
\begin{gather*}
 W_{h, \tilde h} \subset \bigcup_{ k=1}^{ K' ( h ) } U_k \,, \ \ K' ( h ) \leq C' K(h ) 
\,, \ \  U_k = \exp( [ - 1/C , 1/C  ] H_q ) V_k \\ 
 {\rm{diam}} \,( V_k ) \leq C (h/\tilde h)^{\frac12} \,, 
\ \ \ {\rm{sup}}_{(x,\xi) \in U_k} d(( x , \xi) 
, U_k \cap \Sigma ) \leq C h/\tilde
h \,,  \end{gather*}
with a partition of unity on $ W_{h, \tilde h} $,
\[ \sum_{ k=1}^{K' ( h ) } \chi_k = 1 \ \ \text{ on $ W_{h, \tilde h} $,} \ \
\supp \chi_k \subset U_k \,, \ \ \chi_k \in S_{\Sigma, \frac12}^{0,0,-\infty,
- \infty } ( T^* X ) \,.\]
If $ \psi = 1 - \sum_k \chi_k \in  S_{\Sigma, \frac12}^{0,0,-\infty,
- \infty } $ then the condition on the support of $ a $ shows that
\[ \forall \; \alpha, \beta \in \NN^{n}\,, \  \ 
 \partial^\alpha a \, \partial^\beta \psi  \equiv 0 \,.\]
Consequently the calculus of Theorem \ref{th:3.1} 
gives
\[  \Ops ( \psi )  A \in \Psi_{\Sigma, \frac12}^{0, - \infty , 
- \infty , - \infty } ( X ) \,. \]
Hence it suffices to show that for each $ k $ there exists an operator $ R_k $ 
such that 
\[  \Ops ( \chi_k ) A - R_k \in \Psi_{\Sigma, \frac12}^{0, - \infty , 
- \infty , - \infty 
} ( X ) \,, \ \ \rank( R_k ) \leq C   \tilde h^{- n} \,,
\] 
with $ C $ independent of $ k $.
By taking a finer partition (with a 
number of elements $ K'' ( h ) \leq C'' K ( h ) $) we can assume that
$$ \Ops ( \chi_k ) A =  U_k A_k V_k $$
where $ U_k, V_k $ are $h$-semiclassical Fourier Integral Operators 
of the form used in the construction of $ \Ops $, and 
\begin{gather*} A_k = \Opt ( a_k ) \\  
  \supp a_k \subset \{ ( x , \xi', \lambda) \; : \;
  |\lambda | \leq C  \,,  \ \ |x'| + | \xi' | \leq  
C ( h / \tilde h )^{\frac12} \,, |x_1 | < 1/ C \}   \end{gather*}
Consider commuting 
operators 
\begin{gather*}
Q = \left( \tilde h D_{x_1} \right)^2 
+ x_1^2 + ( hD_{x_2} )^2 + x_2^2  
+ \cdots ( hD_{x_n} )^2 + x_n^2 \,, \\  
Q = \Opt ( q ) \,, \ \  q = \lambda^2 + 
x_1^2 + \cdots \xi_n^2 +  
x_n^2 \,, \\ Q' =  ( hD_{x_2} )^2 + x_2^2  
+ \cdots ( hD_{x_n} )^2 + x_n^2 \,.
\end{gather*} 
 If $ \chi \in \CI_{\rm{c}} ( \RR ) $, 
$ \chi ( t) = 1 $ for $ t \leq  \widetilde C $, $ \chi ( t)  
= 0 $ for $ t > \widetilde 2 C $,
then 
\[  \chi ( \tilde h Q'/ h  ) \in \Psi_{\Sigma_0 , \frac12}^{0,0,0,0} 
( \RR^n ) \,, \]  
\[  \chi ( \tilde h Q' / h ) \chi ( Q ) A_k - A_k \in \Psi_{ \Sigma_0 , 
\frac12} ^{0,-\infty, - \infty , - \infty } \,.\]
The standard analysis of the spectrum of harmonic oscillators 
shows that $   \chi ( \tilde h Q' / h ) \chi ( Q )  $ is a finite
rank operator and its rank is bounded by 
$   C' \tilde h^{-n} $.
Hence we can take $ R_k = \chi ( Q) \chi ( \tilde h Q' / h ) A_k $. 
\end{proof}

\section{General upper bounds in regions of size $ h$}
\label{gub}

\subsection{Bound for the number of eigenvalues of a self-adjoint
operator}
\label{pe}

With the calculus developed in \S \ref{smc} we can follow the 
standard procedure of modifying the operator near the energy
surface, now at the limiting scale. For that we introduce the
second small parameter $ \tilde h $ which eventually will 
be fixed, as $ h \rightarrow 0 $.

Let $ \chi \in \CIc ( \RR ; [0,1] ) $ be equal to one near $ 0 $.
We then define 
\begin{equation}
\label{eq:a}
 a ( x, \xi , h ) \stackrel{\rm{def}}{=} \chi \left( \frac{ \tilde h 
p ( x, \xi )}{
 h } \right)   \,.
\end{equation}
Then, in terms of Definition \ref{eq:9.0}, $ a \in S^{0, 0, -\infty , 
- \infty }_{ \Sigma, 0 } ( T^* X) $, $ \Sigma = p^{-1}( 0 ) $.
Although this class of symbols corresponds to a class 
$ \Psi_{\Sigma, 0} ^{\star} $ we will use the larger class $ 
\Psi_{\Sigma, \frac12}^\star $ since it was presented in detail in 
\S \ref{smc}. We stress that this is done for convenience only
and an examination of \S \ref{smc} shows how the simpler calculus
is constructed without the $ S_\frac12 $ complications.

The operator 
\[ \widetilde P \stackrel{\rm{def}}{=} P + i ( h/ \tilde h) A - z \,, \ \
A  \stackrel{\rm{def}}{=} \Opt  ( a ) \,, \]
is elliptic in $ \Psi_{\Sigma}^{0,0,2,1} ( X ) $,
in the
sense that
\[ | \sigs ( \widetilde P ) | > C \left(  d ( \bullet , \Sigma ) 
+ h /\tilde h \right)
\,.\]
This is most clearly seen locally when $ p = \xi_1 $ and $ 
\lambda = (h/\tilde h)^{-1} \xi_1 $:
\[ \tilde p = \frac{ h}{\tilde h} \left( \lambda + i  \chi \left( 
{ \lambda }\right) \right) \,, 
\ \ | \tilde p | > C \frac{ h}{\tilde h} \langle \lambda \rangle \,. \]
Using Theorem \ref{th:3.1} we can construct a parametrix for 
$ \widetilde P - z $, $ | z | \leq C h $, 
uniformly in $ h $, which for $ \tilde h $ small 
enough (keeping $ \tilde h M $ large and constant) gives an exact inverse:
\[  ( \widetilde P - z )^{-1} = {\mathcal O} \left(\frac{ \tilde h}{
h } \right) \; : \; L^2 ( X ) \ \longrightarrow L^2 ( X ) \,, \ \ 
| z | \leq C h \,. \]

Proposition \ref{p:finite} gives a finite 
rank operator $ R $ such that 
\[ A = R + E \,, \ \ \rank( R ) \leq C \tilde h^{-n}  h^{-n+1} \,, \ \ 
E \in \Psi_{\Sigma, \frac12} ^{0, - \infty, - \infty, -\infty} \,, \] 
so that in particular, by Proposition \ref{p:invr},
\[ 
\| E \|_{L^2 \rightarrow L^2} = {\mathcal O} (\tilde h^\infty ) \,. \]
Now we write 
\[ 
\begin{split} P - z & = ( \widetilde P - z ) ( I - i  ( 
h / \tilde h )  (  \widetilde P - z )^{-1} 
( R + E ) ) \\
& = 
( \widetilde P - z ) ( I - i (h/\tilde h) (  \widetilde P - z )^{-1} 
 E  ) ( I + K ( z ) ) \,,\end{split} 
 \]
where 
\begin{gather*}
 K ( z ) = - i (h/\tilde h) ( I + i (h/\tilde h) (  \widetilde P - z )^{-1}  E  )^{-1} 
 (  \widetilde P - z )^{-1} R \,, \\ 
\rank( K ( z ) ) \leq M h^{-n+1} \,, \ \ 
K( z) = {\mathcal O} ( 1 ) \; : \; L^2 ( X ) \; \longrightarrow \; 
L^2 ( X ) \,. \end{gather*}
This implies that 
\[ h ( z ) \stackrel{\rm{def}}{=} \det ( I + K ( z ) ) = {\mathcal O}( 
\exp ( C M h^{-n+1 } ) )\,, \ \ | z | \leq C_0 h \,,\]
and that the zeros of $ h ( z ) $ are the eigenvalues of $ P $ in 
$ |z | \leq Ch $.

\renewcommand\thefootnote{\ddag}%

The bound on the number of eigenvalues will follow standard 
estimates
\footnote{The estimate we need is this: if $ \log | h ( z ) |
\leq K $ in $ R_1 $, where $ R_1 $ is a rectangle, and $ R_2 $ is
another rectangle strictly inside $ R_2 $, $ \log | h ( z_0 )|  > - K $,
for some $ z_0 \in R_2 $, then the number of zeros of $ h( z ) $ in $ R_2 $
is bounded by $ C ( R_1 , R_2 ) K $, with a dilation invariant constant,
$ C( t R_1 , tR_2 ) = C ( R_1 , R_2 ) $.}
once we show that 
\[ \text{ $ |h ( z_0 )| > \exp ( - C M h ^{-n+1} ) $ \ 
at some \ $ z_0 $, \ $ |z_0 | \leq  C_1 h $, \ $ C_1 < C_0 $. } \]
For that we take $ z_0 $ with $ |\Im z_0| > C_1 h $ so that, 
by self-adjointness, 
$$ ( P - z_0 ) ^{-1} = {\mathcal O} ( (C_1h)^{-1} ) \; : \; 
L^2 ( X ) \ \longrightarrow \ L^2 ( X ) \,.$$ 
 We then see that
\[ ( I + K ( z_0 ) )^{-1} =  I + L( z_0  ) \,,\]
where 
\[ I + L ( z ) = ( I + i(h/\tilde h) ( P -z )^{-1} ( R + E ) ) ( I + i (h/\tilde h) ( \widetilde P - z )^{-1} E )  \,, \ \ \  |\Im z| > C_1 h \,. \]
The operator $ L ( z_0 ) $ is of trace class and using the rank of $ R $ 
we have the estimate 
\[  h(z_0) ^{-1} =
\det ( I + L ( z_0 ) ) = {\mathcal O} ( \exp ( C M h^{-n+1 } ) ) \,,\]
which shows that 
\begin{equation}
\label{eq:eig}
|{\rm Spec}\; P \cap D( 0 , h )  | = {\mathcal O} (h^{-n+1} ) \,.
\end{equation}

\subsection{Proof of Theorem \ref{t:2}}
\label{pt2}

To establish the estimate \eqref{eq:th2} we proceed as in the case of
eigenvalues 
but using the scaled operator $ P_\theta $ instead -- see 
\S \ref{rcs}.  
In this section we take
\[ \theta = C_0 h \]
for $ C_0 $, a large and fixed constant. We recall from \S \ref{rsa} that 
\[ X = X_0 \sqcup (\RR^n \setminus B ( 0, R_0 ) ) \,, \]
We can also have more neighbourhoods of infinity (see \S \ref{rsa}) 
but for simplicity of notation we restrict ourselves to the case above.

The operator
 $ P_\theta $ can be written as 
$ P_\theta = P_1 + i P_2 $, where $ P_j $ are formally self-adjoint
on $ L^2 ( X ) $. We consider the {\em Weyl symbols}, $ p_1 $ and
$ p_2 $ (defined, we
recall,  modulo $ {\mathcal O}(h^2 ) $) of these two operators. In 
view of \eqref{eq:gpr} we have 
\begin{gather}
\label{eq:grr}
\begin{gathered}
p_1 \rest_{T^*_{ X_0 \cup \RR^n \setminus B( 0, 2R_0 )  } X} = p \rest_{
T^*_{ X_0 \cup \RR^n \setminus B( 0, 2R_0 )  } X} \,, \ \ \ | p_1 - p | \leq c_1 h 
\langle \xi \rangle^2\,, \\
p_2 \rest_{ X_0} = 0 \,, \ \
- 3 C_1 h |\xi|^2 
\leq p_2  \rest_{  T^*_{\RR^n \setminus B( 0, 2R_0 ) } X } 
 \leq - C_1 h |\xi|^2 \,,
\ \  | p_2 | \leq c_1 h \langle \xi \rangle^2\,. 
\end{gathered}
\end{gather}
We also note that our assumptions give
\[ | p_1 | \leq \delta \ \Longrightarrow \ \langle \xi \rangle \sim |\xi| \,,\]
with some small fixed $ \delta $.
Now let $ \chi \in \CIc ( \RR , [0, 1] )$, satisfy
\[ \chi ( t) = \left\{ \begin{array}{ll} 1 & |t| \leq 1 \\
0 & |t| \geq 2 \end{array} \right. \]
With this $ \chi $ we define
\[  a ( x , \xi , h ) \stackrel{\rm{def}}{=} \chi \left( 
\frac{ \tilde h p ( x , \xi ) }
{   h } \right) \chi \left ( \frac{|x|}{ 3R_0 } \right) \,,\]
where $\tilde h$ is small, 
\[  \tilde h^{-1}  \sim C_1 \sim c_1 \sim C_0 \,,\]
in \eqref{eq:grr} but eventually fixed.

We now choose a {\em compact } $ \CI $ hypersurface $ \Sigma $ so that
$$ \Sigma \cap T^*_{X_0 \sqcup B ( 0 , 6R_0 ) \setminus B( 0, R_0 )  } X 
 = \{ ( x, \xi ) \; : \; 
p ( x , \xi ) = 0 \} \,. $$ 
It follows that 
\[ a ( x , \xi , h ) \in S^{0,0,-\infty,0}_\Sigma  ( T^* X ) \,, \]
and continuing in the same spirit as in \S \ref{pe}, we put
\[ A \stackrel{\rm{def}}{=} \Opt ( a ) \,.\]
On the level of symbols we have 
\begin{equation}
\label{eq:p1}
 | p_1 | \leq (h/\tilde h)  \
\Longrightarrow \ 
 - p_2 ( x , \xi ) + (h/\tilde h) a ( x ,\xi , h )  \geq  (h/\tilde h)/C 
 \,. \end{equation}
Hence if we put
\[ \widetilde P = P_\theta - i (h/\tilde h) A \,, \]
then for 
$ | \Re z | \leq Ch $,  $ \Im z \geq - Ch $,
\[  (\tilde h /h) ( \widetilde P - z ) \in 
\Psi_{\Sigma}^{0,0,2,1} (  X) \,.\]
We claim that $ \widetilde P - z $ is invertible for $ z \in 
D ( 0 , C h ) $ and 
\begin{equation}
\label{eq:wpi}
 \| ( \widetilde P - z )^{-1} \| = {\mathcal O} ( \tilde h / h ) \,. 
\end{equation}
In fact, let $ W \subset T^* X $ be set in which 
$  ( \tilde h / h ) ( \tilde p - z ) $ is elliptic in 
$ S_{\Sigma}^{0,0,2,1} ( T^* X ) $. The estimate \eqref{eq:p1} shows that
$ \Im ( \tilde p -z ) \leq - C ( h /\tilde h ) $ on a neighbourhood 
of $ \complement W \cap 
\{ | \Re(  \tilde p -z ) | < \delta \}  $. Now, 
let $ \Psi_1 $ and $ \Psi_2 $ be as in Lemma \ref{l:r2}, with 
$ \esssupp \Psi_1 \subset W $. Then proceeding as in 
\S \ref{cop} we obtain 
\[ \| ( \widetilde P - z ) \Psi_2 u \| \geq ( h /\tilde h ) 
\| \Psi_2 u \|/C  - {\mathcal O} ( h^\infty ) \| u \| \,, \ \ 
u \in \CIc ( X ) \,.\]
Ellipticity of $ ( \tilde h / h ) ( \widetilde P - z ) $ in 
$ W $ shows that 
\[ \| ( \tilde h / h ) ( \widetilde P - z ) \Psi_1 u \| \geq 
 \| \Psi_1 u \|/C  - {\mathcal O} ( \tilde h^\infty ) \,, \]
that is 
\[ \| ( \widetilde P - z ) \Psi_1 u \| \geq 
( h / \tilde h )  \| \Psi_1 u \|/C  - 
{\mathcal O} (h  \tilde h^\infty ) \,. 
\]
Lemma \ref{l:r3} then gives \eqref{eq:wpi}.

We can now proceed  as in \S \ref{pe}
and obtain the bound on the number of resonances \eqref{eq:th2}. 
The complex analytic argument outlined in the footnote is used in a
rectangle 
$$ [-C h , Ch] + i[ C_2 h , -C h ] $$ 
where $ C_2 $ is large
enough to guarantee a lower bound for $ - \Im ( P_\theta - z ) $ 
when $ \Im z \sim C_2 h $. We can still take $ C_2 $ proportional
to the other large constants.
\stopthm

\section{The escape function for hyperbolic flows and its $h$
dependent regularizations}
\label{efh}

In this section we modify \cite[Sect.5]{SjDuke} and 
construct a regularized escape function depending on a small
parameter, essentially $ h/ \tilde h$. We recall that we 
assume that $ p \in \CI ( T^* X; \RR ) $ satisfies
\begin{gather}
\label{eq:asp}
\begin{gathered} p = 0 \ \Longrightarrow \ dp \neq 0 \ \\
|x | \geq R \,, \ \ | p ( x , \xi ) | < 2 \delta  \ \Longrightarrow \
\exp t H_p ( x , \xi ) \rightarrow \infty \ \text{ for 
either $ t \rightarrow \infty $ or $ t \rightarrow - \infty$.}
\end{gathered}
\end{gather}
We also recall the result of \cite[Appendix]{GeSj}:
\begin{prop}
\label{l:gsa}
Suppose that \eqref{eq:asp} holds and that $ \widehat K $ is the trapped set,
\begin{equation} 
\label{eq:wdK}
\widehat K\stackrel{\rm{def}}{=} \{ \rho \in T^* X \; : \; 
\exp ( t H_p ) ( \rho ) \not \rightarrow \infty \,, \ t \rightarrow
\pm \infty \,, |p ( \rho ) |\leq \delta \} \Subset T^*X \,. \end{equation}
Then for any two neighbourhoods, $ U , V $, 
 of $ \widehat K $,  $ \overline U \subset 
V $ there exists $ G_0 \in \CI ( T^*X ) $ such that 
\begin{gather}
\label{eq:gsa}
\begin{gathered}
\supp G_0 \subset T^*X \setminus U \,, \  \ H_p G_0 \geq 0 \,, \ \ 
H_p G\rest_{ p^{-1}( [ 2 \delta, 2 \delta ]) } \leq C \,, \\
H_p G_0 \rest_{ p^{-1}( [ - \delta , \delta ] ) \setminus V } \geq 1 \,. 
\end{gathered}
\end{gather}
\end{prop}

\subsection{Dynamical assumptions}
\label{esf}

We start with the hyperbolicity assumptions \cite[\S 5]{SjDuke} 
weaker than the more standard assumptions in \S \ref{in}.
Let 
$ \widehat K $ be the compact trapped set near zero energy given by 
\eqref{eq:wdK}.
The trapped set at zero energy is given by 
$ K = \widehat K \cap p^{-1} ( 0 ) $. 
We also have 
$ \widehat K = \widehat \Gamma_+ \cap \widehat \Gamma_- $,
where 
\begin{equation}
\label{eq:gaha}
 \widehat \Gamma_\pm  \stackrel{\rm{def}}{=} \{  ( x , \xi ) \in T^* X \; : 
\; 
| p ( x , \xi )  | \leq \delta \,, \ \ \exp ( t H_p ) ( x , \xi ) 
\not \rightarrow \infty \,, \ t \rightarrow \mp \infty \} \,,
\end{equation}
and the sets $ \widehat K $, $ \widehat \Gamma_\pm $ are clearly 
invariant under the flow, 
\begin{equation}
\label{eq:invf}
 \exp ( t H_p ) ( \widehat K ) \subset \widehat K \,, \ \
\exp ( t H_p ) ( \widehat \Gamma_\pm ) \subset \widehat \Gamma_\pm \,.
\end{equation}
We can now state the dynamical hypothesis.
\begin{itemize}
\item In a neighbourhood, $ \Omega_{\rho_0 } $ of any $ \rho_0 \in K $,
\[    \widehat \Gamma_\pm = \bigcup_{\rho \in \Omega_{\rho_0 }
\cap \widehat \Gamma_{\pm}}  \widehat 
\Gamma_{\pm,
\rho}  \,, \ \  \rho \in \widehat \Gamma_{\pm, \rho} \,, \]
\[ \widehat \Gamma_{\pm,\rho} \cap \widehat \Gamma_{\pm,\rho'} = \emptyset\,,
\ \text{ or } \ \widehat \Gamma_{\pm,\rho} = \widehat \Gamma_{\pm,\rho'}
 \,. \] 
\item Each $ \widehat \Gamma_{\pm,\rho} $ is a closed 
$ {\mathcal C}^1 $ manifold of 
dimension $ n + d $, with $ d \geq 0$ fixed, and the dependence 
\[   \Omega_{\rho_0} \cap \widehat \Gamma_\pm  \ni \rho \longmapsto
 T_\rho \widehat \Gamma_{\pm,\rho}  \]
{is continuous.} 
 \item If  $ E_\rho^\pm \stackrel{\rm{def}}{=} 
 T_\rho   \widehat \Gamma_{\pm, \rho}  $, 
then $ E_\rho^+ + E_\rho^- = T_\rho p^{-1} ( p ( \rho ) ) \subset 
T_\rho (T^* X )  $,  $ \RR H_p ( \rho ) 
\in E_\rho^\pm $, and 
\begin{equation}
\label{eq:hyp8}
  \| d (\exp t H_p )_\rho ( X ) \| \leq C e^{ \pm \lambda t } \| X \| \,, 
\ \  \rho \in K \,, \ \ \text{ for all $ X \in T_\rho ( T^*X ) / E^\mp _ \rho  $, $ \mp t \geq 0 $.} \end{equation}
\end{itemize}
The above definition makes sense 
since by \eqref{eq:invf} $ d ( \exp  t H_p  )_\rho ( E_\rho^\pm ) = 
E_{ \exp t H_p ( \rho) } $, $ \rho \in \widehat \Gamma_\pm $, 
we have 
$$ d (\exp t H_p )_\rho \; T_\rho ( T^* X) /  E^\mp _ \rho  
\longrightarrow T_{\exp t H_p(\rho)} ( T^* X) /  E^\mp _ {\exp t H_p(\rho)}  
\,, \ \ \rho \in K \,, 
$$ 
and we choose continuously dependent norms in the last estimate in 
\eqref{eq:hyp8}. We also note that $ X \in T_\rho ( T^* X )/ E^\mp_\rho $
implies that $ X $ can be identified with a vector 
tangent to $ p^{-1} ( p( \rho ) ) $.

In \cite[\S 5]{SjDuke} it is shown that there exist two 
functions, $  \varphi_\pm 
 \in {\mathcal C}^{1,1} ( T^* X ) $, 
$ \varphi_\pm \geq 0 $, 
$ H_p^k \varphi_\pm \in  {\mathcal C}^{1,1} ( T^* X ) \,,  \ k \in \NN $, 
such that for $ \rho  $ in a small neighbourhood of $ K $, 
\[ \begin{split}
\mp  H_p \varphi_\pm ( \rho) 
& \sim  \varphi_\pm ( \rho ) \,, \ \ H_p^k \varphi_\pm (\rho ) = {\mathcal O} ( \varphi_\pm ( \rho ) )\,, \   \ k \in \NN 
 \,,
\,, \\ \varphi_\pm ( \rho ) & \sim d ( \rho, \widehat \Gamma_\pm)\,, \ \ 
\varphi_+ ( \rho ) + \varphi_- ( \rho )
 \sim d( \rho   , \widehat 
K)^2 \,, 
\end{split} \]
and where $ d( \bullet , \Gamma ) $ is the distance to a closed set $ \Gamma $.
The notation $ f \sim g $, 
 means that there exists a constant $ C > 0 $ such 
$$ 0 \leq  g /C  \leq f \leq C g \,. $$
The simple model \eqref{eq:simo} is given in \S \ref{oof}. Here we
modify the construction to obtain suitably regularized functions $ \hph_\pm $

\subsection{Regularization of $ \varphi_\pm $.}
We start with two general lemmas:
\begin{lem}
\label{l:whitney?}
Suppose $ \Gamma \subset \RR^m $ is a closed set. For any $ \epsilon > 0 $
there exists $ \varphi_\epsilon \in \CI ( \RR^m ) $ such that 
\[ \varphi_\epsilon \geq \epsilon \,, \ \ 
 \varphi_\epsilon \sim d ( \bullet, \Gamma)^2 + \epsilon \,, \ \
\partial^\alpha \varphi_\epsilon = {\mathcal O} ( \varphi_\epsilon ^{1- 
|\alpha|/2 } ) \,,\]
uniformly on compact sets.
\end{lem}
\begin{proof}
We can find a sequence $ x_j \in \RR^m $ such that 
\begin{gather*}
 \bigcup_j B ( x_j , d( x_j , \Gamma ) / 8 ) = \RR^m \setminus \Gamma 
\,, \\ 
\text{every $  x \in Q \setminus \Gamma $, $ Q \Subset 
\RR^m $, is in at most $ N_0 = N_0 ( Q ) $ 
balls $ B ( x_j , d( x_j , \Gamma ) / 2 ) $.}
\end{gather*}
Let $ \chi \in \CIc ( \RR^m ; [0,1]) $ be supported in $ B (0 , 1/4 ) $,
and be identically one in $ B ( 0 , 1/8)$. We define
\[ \varphi_\epsilon ( x ) \stackrel{\rm{def}}{=} \epsilon 
+ \sum_{ d( x_j , \Gamma ) > \sqrt \epsilon } d ( x_j , \Gamma )^2 
\chi \left( \frac{  x - x_j  } { d ( x_j , \Gamma ) + \sqrt \epsilon }
\right) \]
We first note that the number non-zero terms in the sum is uniformly bounded
by $ N_0 $. In fact, $ d ( x_j , \Gamma) + \sqrt \epsilon < 2 d ( x_j , 
\Gamma) $, and hence if $ \chi( ( x - x_j ) / ( d ( x_j , \Gamma ) + 
\sqrt \epsilon ) ) \neq 0 $ then  
$$ 1/4 \geq |x-x_j| / ( d ( x_j , \Gamma ) + \sqrt \epsilon ) 
\geq (1/2) |x - x_j | / d( x_j , \Gamma )  \,,$$
and $ x \in B ( x_j , d ( x_j , d ( x_j , \Gamma) )  / 2 ) $.
This shows that $ \varphi_\epsilon ( x ) 
\leq 2 N_0 (\epsilon +  d ( x, \Gamma )^2 ) $, and 
\[ \partial^\alpha \varphi_\epsilon ( x) = 
{\mathcal O} ( ( d ( x, \Gamma )^2 + \epsilon )^{1 - |\alpha|/2 } ) \,,\]
uniformly on compact sets.

To see the lower bound on $ \varphi_\epsilon $ 
we first consider the case when $ d ( x , \Gamma ) \leq C \sqrt \epsilon $.
\[ \varphi_\epsilon ( x ) \geq \epsilon \geq ( \epsilon + d ( x , \Gamma)^2)
 / C'  \,.\]
If $ d ( x , \Gamma ) > C \sqrt \epsilon $ then for at least one 
$ j $, $ \chi ( ( x - x_j ) / (d ( x_j , \Gamma ) + \sqrt \epsilon ) )
= 1 $ (since the balls $ B( x_j , d ( x_j , \Gamma ) /8 ) $ 
cover the complement of $ \Gamma $, and $ \chi ( t ) = 1 $ if $ |t| \leq 1/8 
$). Thus 
$$ \varphi_\epsilon ( x ) \geq \epsilon + d ( x_j , \Gamma ) ^2 
\geq ( \epsilon + d ( x , \Gamma)^2 )/ C \,, $$
which concludes the proof. \end{proof}
For future use we also record the following
\begin{lem}
\label{l:ef}
Suppose $ \varphi \in {\mathcal C}^{1,1} ( \RR^m )  $, $ \varphi \geq 0 $,
and for a vectorfield 
$ V \in \CI( \RR^m ; \RR^m ) $, $ V^k \varphi = {\mathcal O} ( \varphi)  $, 
$ V^k \phi \in {\mathcal C}^{1,1} ( \RR^m ) $, $ k \in \NN $. 
Then, uniformly on compact sets,  
\[ d V ^k \varphi = {\mathcal O} ( \varphi^{\frac12}) \,, \ \ 
k \in \NN \,. \]
\end{lem}
\begin{proof}
For some $ C > 0 $ the $ {\mathcal C}^{1,1} $ function $ C \varphi -
V^k  \varphi $
is non-negative. Hence using the standard estimate based on Taylor's 
formula, 
\[ | d \varphi |^2 = {\mathcal O} ( \varphi ) \,, 
\ \ 
| d ( C \varphi - V^k \varphi ) |^2 = {\mathcal O} ( C \varphi - V^k \varphi ) 
= {\mathcal O} ( \varphi ) \,. \]
The lemma follows.
\end{proof}

We now have 
\begin{prop}
\label{p:es1}
Let $ \widehat \Gamma_\pm $ be given by \eqref{eq:gaha}. For any 
small $ \epsilon > 0 $ there 
exist functions $ \hph_\pm \in \CI ( T^* X ; [0, \infty) ) $ such that
in a neighbourhood of $ \widehat K $, 
\begin{equation}
\label{eq:es1'}
\begin{split}
\hph_\pm  ( \rho ) & \sim d( \rho  , \widehat \Gamma_\pm )^2 + C  \epsilon 
 \,, \\
\mp  H_p \hph_\pm ( \rho ) + C \epsilon 
& \sim  \hph_\pm ( \rho ) 
\,, \\
\partial^\alpha 
H_p^k \hph_\pm (\rho ) & = {\mathcal O} ( \hph_\pm ( \rho )^{ 1 
- |\alpha|/2} )\,, \   \ k \in \NN \,, \\
 \hph_+ ( \rho ) + \hph_- ( \rho )
& \sim d( \rho , \widehat 
K)^2 +  C \epsilon \,.
\end{split}
\end{equation}
\end{prop}
\begin{proof}
We modify the arguments of \cite[\S 5]{SjDuke}, roughly speaking,
adding an $ {\mathcal O} ( \epsilon ) $ error to all the estimates.
Let $ \varphi_\pm $ be the functions obtained using Lemma \ref{l:whitney?}
with $ \Gamma = \Gamma_\pm $. 
We now put
\begin{gather*}
 \hph_\pm ( \rho ) 
\stackrel{\rm{def}}{=} \int_\RR g_T ( t ) \varphi
_\pm ( \exp t H_p ( \rho ) ) dt \,, \\ g_T \in \CIc ((-1,T+1) ) \,, \ \ 
\supp g'_T \subset [-1,1]\cup[T-1,T+1]\,, \\ g'_T \rest_{[-1,1]} \geq 0 \,, \
 \ g'_T \rest_{[ T-1, T+1] } \leq 0 \,, \ \ g'_T ( 0 ) = 1 \,, \ \ 
g_T' ( T ) = - 1 \,.
 \end{gather*}
To check \eqref{eq:es1'} we note that, by definition, 
 $ \varphi_\pm ( \rho ) \sim 
d ( \rho , \widehat \Gamma_\pm )^2 +  C \epsilon  $.
The assumptions \eqref{eq:hyp8} imply (see \cite[Lemma 5.2]{SjDuke}) 
that 
\[
  \exists \; C\,, \ \forall \; T \geq 0 \,, \ \exists \; \Omega_T 
\supset K \,, \text{an open set, } \ \ 
d ( \exp (\pm T H_p) ( \rho ) , \widehat \Gamma_\pm ) \leq C 
e^{ - T / C } d( \rho, \widehat \Gamma_\pm ) \,. \]
Hence, with constants depending on $ T $, 
\begin{gather*}
 \widehat \varphi_+ ( \rho ) \sim \varphi_+ ( \exp ( T H_p ) ( \rho )) 
\sim  \varphi_+ ( \rho) \sim d( \rho, \Gamma_+)^2 + C \epsilon \,, \\
\widehat \varphi_- ( \rho ) \sim \varphi_- (  \rho )
 \sim d( \rho, \Gamma_-)^2 + C \epsilon \,. \end{gather*}
This shows the first statement in \eqref{eq:es1'}.

The assumptions on $ g_T $ also show that 
\[ H_p \hph_\pm ( \rho ) \sim \varphi_\pm ( \exp T H_p (\rho )) 
- \varphi_\pm ( \rho ) \sim d (  \exp T H_p (\rho ) , \widehat 
\Gamma_\pm )^2 - 
d ( \rho, \widehat \Gamma_\pm )^2 + {\mathcal O} ( \epsilon ) \,. \] 
so that for $ T $ large enough
and for $ \rho $ in a small neighbourhood of  $K $, 
(again with $ T $ depenendent constants)
\[ \mp H_p \hph_\pm ( \rho ) + C \epsilon 
\sim  d ( \rho, \widehat \Gamma_\pm )^2 + 
C' \epsilon 
\sim \hph_\pm ( \rho ) \,.\]
This proves the second part of \eqref{eq:es1'}.
The third part is proved using Lemma \ref{l:ef} for $ |\alpha | =1 $
and the estimates on $ \varphi_\pm $ in general.

To prove the last statement in \eqref{eq:es1'} we first see that the
transversality, $ E_{\rho_0} ^+ 
+ E_{\rho_0}^- = T_{\rho_0 } ( T^* X ) $, and  the continuity,
$ \rho \mapsto E_{\rho}^\pm $, assumed  
 in \eqref{eq:hyp8} imply that for $ \rho\,, \rho_1\,, \rho_2\,, $ 
near a point
$ \rho_0 \in K $, 
\[ d ( \rho , \widehat \Gamma_{+ , \rho_1 } \cap \widehat \Gamma_{-, \rho_2} ) 
\sim d ( \rho , \widehat  \Gamma_{+ , \rho_1 } )  + 
d ( \rho , \widehat \Gamma_{-, \rho_2} )  \,. \]
Hence
\[ \begin{split} \hph_+ ( \rho) + \hph_- ( \rho ) + {\mathcal O}
(\epsilon )  & \sim 
 d ( \rho , \widehat  \Gamma_{+ } )^2  + 
d ( \rho , \widehat \Gamma_{-} )^2  + C\epsilon  \\
& \leq 
 d ( \rho , \widehat  \Gamma_{+ , \rho' } )^2  + 
d ( \rho , \widehat \Gamma_{-, \rho'} )^2  + C \epsilon \\
&  \sim 
 d ( \rho , \widehat \Gamma_{+ , \rho' } \cap \widehat \Gamma_{-, \rho'} )^2  + C \epsilon \,.
\end{split}\]
If we choose $ \rho' \in K $ so that $ d ( \rho , \widehat K ) = 
d ( \rho, \rho' ) $
then 
$$  d ( \rho , \widehat \Gamma_{+ , \rho' } \cap \widehat \Gamma_{-, \rho'} )^2 \leq d ( \rho , \rho' )^2 = d ( \rho,  \widehat K) ^2\,,  $$
proving that 
$$ \hph_+ ( \rho ) + \hph_- ( \rho )  
\leq d ( \rho,  \widehat  K)^2 + {\mathcal O} ( \epsilon ) \,.$$
The opposite inequality is obtained by choosing $ \rho_\pm 
\in  \widehat  \Gamma_\pm $ such that
$ d ( \rho, \rho_\pm ) = d ( \rho, \widehat \Gamma_\pm ) $. Then using 
the transversality of $ \widehat \Gamma_+ $, $ \widehat \Gamma_-$
\[  \begin{split} d ( \rho , \widehat K )^2 & \leq 
 d ( \rho , \widehat \Gamma_{+ , \rho_+} \cap \widehat \Gamma_{-, \rho_-} )^2
\sim
 d ( \rho ,  \widehat  \Gamma_{+, \rho_+}  )^2  +
  d ( \rho ,  \widehat  \Gamma_{- , \rho_-}  )^2 \\ & 
\leq  d ( \rho, \rho_+ )^2 + d ( \rho, \rho_- )^2 = 
d ( \rho, \widehat \Gamma_+ )^2+ d ( \rho, \widehat \Gamma_-)^2 \\
& \leq  \hph_+ ( \rho ) + \hph_- ( \rho ) + {\mathcal O} ( \epsilon ) \,.
\end{split} \]
\end{proof}

\subsection{Regularized escape function}
\label{res}

We now use the functions constructed in Proposition \ref{p:es1} to 
obtain an escape function near $ K$. We first need the following
\begin{lem}
\label{l:fes}
Then for 
$ |\alpha | + k \geq 1 $ we have 
\begin{equation*}
\partial_\rho^\alpha H_p^k \log ( \widehat \varphi_\pm ) =  {\mathcal O} ( \widehat \varphi_\pm ^{-\frac{|\alpha|}2} ) 
\end{equation*}
\end{lem}
\begin{proof}
Let $ f ( t ) = \log ( t ) $. Then 
\[ f^{(k)} ( \hph_\pm ) = {\mathcal O} \left( \frac{1}{ 
\hph_\pm ^k } \right) \,, \ \  k \geq 1 \,, \]
and for $ |\alpha | + k \geq 1 $, $ \partial_\rho^\alpha H_p ^k f ( 
\hph_\pm ) $ is a finite linear combination of terms
\[ f^{(l)} ( \hph_\pm ) \left(  \partial_\rho^{\alpha_1} H_p ^{k_1} \hph_\pm \right)
\cdots  \left(  \partial_\rho^{\alpha_\ell} H_p ^{k_\ell} \hph_\pm \right) 
= {\mathcal O} ( 1) \prod_{ j=1}^\ell \frac{\partial_\rho^{\alpha_j} H_p ^{k_j}\hph_\pm  }
{ \hph_\pm } \,, \]
with
\[ | \alpha_j | + k_j \geq 1 \,, \ \ \alpha_1 + \cdots + \alpha_\ell = 
\alpha \,, \ \ k_1 + \cdots + k_\ell = k \,. \]
The estimates in \eqref{eq:es1'} show that 
${\partial_\rho^{\alpha_j} H_p ^{k_j} \hph_\pm  }/
{ \hph_\pm }  = {\mathcal O}(   \hph_\pm^{-|\alpha_j|/2}  )
$, 
and hence 
\[ \partial_\rho^\alpha H_p ^k f ( \hph_\pm ) = {\mathcal O} (  
\hph_\pm^{-\frac{|\alpha|}2 } ) \,,\]
proving the lemma.
\end{proof}

We are now ready for the main results of this section. 
\begin{lem}
\label{l:es2}
Let $ \hph_\pm $ be given in Proposition \ref{p:es1} and 
\begin{equation}
\label{eq:es2}
\widehat G \stackrel{\rm{def}}{=} \left( 
\log (  M \epsilon + \hph_- ) - \log (  M \epsilon + \hph_+ ) \right) \,.
\end{equation}
Then in a neighbourhood of $ K $ we have 
\begin{gather}
\label{eq:es4}
\begin{gathered}
\partial_\rho^\alpha H_p^k \widehat G = 
 {\mathcal O}_M (   \min ( 
\widehat \varphi_+ , \hph_- )^{-\frac{|\alpha|}2
} )  = {\mathcal O}_M (  \epsilon^{ -\frac{|\alpha|}2
} ) \,, \ \ |\alpha | + k \geq 1 \,, \\
  d ( \rho, \widehat K)^2 \geq C \epsilon
\Longrightarrow \ H_p \widehat G \geq 1 / C  \,, 
\end{gathered}
\end{gather}
where, for the second 
estimate, $ M $ has to be chosen 
large enough, independently of $ \epsilon $, and $ C $ is a large
constant.
\end{lem}
\begin{proof}
We observe that, with constants depending on $ M $, $ \hph_\pm + M \epsilon $
has the same properties as $ \hph_\pm $. Hence
the estimates on $ \partial_\rho^\alpha H_p^k \widehat G $
follow directly from the definition
\eqref{eq:es2} and from Lemma \ref{l:fes}. To check 
the second part of \eqref{eq:es4}
we compute, using  Proposition \ref{p:es1},
\[ H_p \widehat G = \left( \frac{ H_p \hph_- }{ \hph_- 
+  M \epsilon } - \frac{ H_p \hph_+ }{ \hph_+ +  M \epsilon } \right)
\geq \frac{1}{C_1} \left( \frac{ \hph_- -C_2 \epsilon }{\hph_- + 
 M \epsilon } + \frac{ \hph_+ -C_2 \epsilon }{\hph_+ +  M \epsilon } 
\right) \,.
\]
From \eqref{eq:es1'} we also have 
\[ d ( \rho , \widehat K )^2 \geq C \epsilon \ \Longrightarrow \
\max ( \hph_+ , \hph_- ) \geq  (C/2 - {\mathcal O} (1 ) ) \epsilon
>  C_3 \epsilon \,,\]
where $ C_3 $ can be as large as we like depending on the choice of 
$ C $.
Hence, since $ x \mapsto ( x - C_2 )/( x + M ) $ is increasing,
\[ H_p \widehat G \geq 
 \frac{1}{C_1} \left( \frac{ C_3  -C_2 }{ C_3 + 
 M  } - \frac{ C_2 }{ M  } 
\right) \geq \frac{1}{C} \,,\]
if we choose $ C_3 \gg M \gg C_2 $.
\end{proof}

We now modify $ \widehat G $ using $ G_0 $ given in Proposition \ref{l:gsa}:
\begin{prop}
\label{p:es2}
Let us fix $ \delta_0 > 0 $.
Then there exist $ \widehat \chi, \chi_0 \in \CIc ( T^* X ) $, $ C_0 > 0 $, 
and a
neighbourhoood $ V $ of $ K $, such that
\[ G \stackrel{\rm{def}}{=} \widehat \chi \widehat G + 
C_0 \left( \log \frac{1}{\epsilon } \right) \chi_0 G_0 \,, \]
satisfies
\begin{gather}
\label{eq:es3}
\begin{gathered}
  \partial^\alpha H^k_p G  = \left\{ 
\begin{array}{ll} {\mathcal O} ( 
\log ( 1 / \epsilon ) ) & \alpha = 0 \\
{\mathcal O}( \epsilon^{-|\alpha|/2} ) & \text{otherwise} \end{array}
\right. \,, \\
  d ( \rho, \widehat K)^2 \geq C \epsilon \,, \ \rho \in V \ 
\;  \Longrightarrow \;  H_p  G ( \rho ) \geq 1 / C  \,,  \\  
 \rho \in p^{-1}([-\delta, \delta] )
\setminus V \,, \ \ | x ( \rho ) | \leq 3 R_0 \
\;  \Longrightarrow \;  H_p  G ( \rho ) \geq \log ( 1/ \epsilon )  \,, \\
 H_p G ( \rho )   \geq - \delta_0 \log( 1/ \epsilon ) \,, \ \ \rho
\in T^* X \,. 
\end{gathered}
\end{gather}
In addition we have 
\begin{equation}
\label{eq:orderf}
 \frac{\exp G ( \rho ) }{ \exp G( \mu ) } \leq C_0 \left \langle
\frac{ \rho - \mu }{\sqrt \epsilon } \right \rangle^{N_0} \,,\end{equation}
for some constants $ C_0 $ and $ N_0 $.
\end{prop}
\begin{proof}
We obtain $ G_0 $ from Proposition \ref{l:gsa} taking for $ V $ a
neighbourhood of $ \widehat K $ in which 
the estimates of Lemma \ref{l:es2} hold. 
We have 
$  \partial^\alpha H^k_p G_0 = {\mathcal O}_{k, |\alpha|} ( 1 ) $, and  
consequently for any $ \chi_0 \in \CIc ( T^* X) $, 
\[ \partial^\alpha H^k_p \left( \log ( 1/\epsilon ) \chi_0 G_0 \right)
=  {\mathcal O}_{k, |\alpha|} ( \log( 1/ \epsilon  ))
=  \left\{ 
\begin{array}{ll} {\mathcal O} ( 
\log ( 1 / \epsilon ) ) & \alpha = 0 \\
{\mathcal O}( \epsilon^{-|\alpha|/2} ) & \text{otherwise} \end{array}
\right. \,. \]
From Lemma \ref{l:es2} we obtain, again for any $ \widehat \chi \in 
\CIc ( T^* X ) $, 
\[ \partial^\alpha H^k_p ( \widehat \chi \widehat G ) 
=  \left\{ 
\begin{array}{ll} {\mathcal O} ( \log( 1/\epsilon ) )
 & \alpha = 0 \\
{\mathcal O}( \epsilon^{-|\alpha|/2} ) & \text{otherwise} \end{array}
\right. \,. \]
The loss compared to 
\eqref{eq:es4} is due to the presence of the cut-off function.

We take $ \chi_0 \in \CIc ( T^* X ; [0,1] ) $ to be identically equal to 
$ 1 $ in 
$$ p^{-1} ( [ - \delta , \delta ] ) \cap \{ ( x , \xi) \; : \;
|x | \leq 3 R_0 \} \,. $$
 For $ \widehat \chi \in \CIc ( T^* X ) $ we 
take a function which is supported in a neighbhourhood of $ \widehat K $
where \eqref{eq:es4} holds, and identically $ 1 $ in $ V $. Hence
for $ \rho \in  p^{-1}([-\delta, \delta] ) 
\setminus V $, $ | x ( \rho ) | \leq 3 R_0 $, 
\[ H_p G ( \rho ) = C_0  \log( 1/ \epsilon ) H_p G_0 ( \rho ) 
+ H_p ( \widehat \chi 
\widehat G ) ( \rho ) \geq C_0 \log ( 1/ \epsilon ) - {\mathcal O} ( 1 ) 
\log( 1/ \epsilon ) \geq \log ( 1/ \epsilon ) \,, \]
if $ C_0 $ is taken large enough. For $ \rho \in V $, $ \widehat \chi ( \rho ) 
= 1$, and 
\[ H_p G ( \rho ) = 
 C_0  \log( 1/ \epsilon ) H_p G_0 ( \rho ) 
+ H_p  \widehat G ( \rho ) \geq H_p \widehat G ( \rho ) \,, \]
and if $ d ( \rho, \widehat K ) \geq C \epsilon $, $ H_p G( \rho ) 
\geq 1/C $.
To complete the proof of \eqref{eq:es3} we need to define $ \chi_0 $
for $ |x| \geq R_0 $. That is essentially 
done as in \S \ref{cef} where it was based on Lemma \ref{l:chi}. 
Let $ T $ and $ R $ be large positive constants to be fixed later, 
$ \chi( t) $ be given by Lemma \ref{l:chi}, and let $ \psi \in 
\CIc ( \RR ; [0,1]) $ be equal to $ 1 $ for $ |t| \leq 1 $, and to 
$ 0 $ for $ |t| \geq 2 $. We define
\[ \chi_0 ( \rho ) \stackrel{\rm{def}}{=} \frac{\chi ( G_0  ( \rho ) ) }
{G_0 ( \rho ) } 
\psi \left( \frac{ p( \rho ) }{\delta} \right) \psi \left( \frac{ 
| x( \rho ) | }{ R } \right) \,.\]
Then
\[ \begin{split}
 H_p ( \chi_0 G_0 ) ( \rho ) & = \chi' ( G_0 ( \rho ) ) H_p G_0 ( \rho) 
\psi \left( \frac{ p( \rho ) }{\delta} \right) \psi \left( \frac{ 
| x( \rho ) | }{ R } \right) \\ 
& \ \ \ \ \ + \; 
\frac{1}{R} \chi ( G_0 ( \rho ) ) 
\psi \left( \frac{ p( \rho ) }{\delta} \right) \psi' \left( \frac{ 
| x( \rho ) | }{ R } \right) H_p ( |x| ) ( \rho ) \,,
\end{split} \]
and 
\[  H_p ( \chi_0 G_0 ) ( \rho ) \geq - C_1 \left ( \alpha 
+  \frac{T}{R} \right)\,,\]
where $ C_1 $ is independent of $ T $ and $ R $: we note that 
\eqref{eq:gsa} guarantees the boundedness of $ H_p G_0 $, and the
assumptions on $ p $ imply that $ H_p ( |x|) $ is uniformly bounded
for $ |p | \leq 2 \delta $.
For any $ \alpha > 0 $ we can choose $ T = T( \alpha ) $ such that
$ |G_0 ( \rho )| \leq \alpha T $ for $ |x( \rho ) | \leq 3 R_0 $, 
$ | p ( \rho ) | \leq 2 \delta $. We then choose $ \alpha $ and $ R $
so that 
\[ C_0 C_1 ( \alpha + T ( \alpha ) / R  ) < \delta_0\,. \]
Hence for $ |x ( \rho ) | \geq R_0 $
\[ H_p G = C_0 \log( 1 / \epsilon) H_p ( \chi_0 G_0 ) \geq - \delta_0 
\log ( 1 / \epsilon ) \,, \]
which is the last statement in \eqref{eq:es3}.

It remains to show \eqref{eq:orderf} and for 
simplicity of presentation we replace $ T^* X $ with $ \RR^{2n}$.
We
first prove that
\begin{equation}
\label{eq:ordph} \frac{ \hph_{\pm } ( \rho) + M\epsilon } {\hph_\pm ( \mu)
+ M\epsilon } \leq 
C_1 \left \langle
\frac{ \rho - \mu }{\sqrt \epsilon } \right \rangle^2 \,, \ \ M \geq 0 \,, 
\end{equation}
with constants depending on $ M $. We can replace $ \hph_\pm + M \epsilon $
with $ \hph_\pm $, as $  \hph_\pm + M \epsilon \sim_M \hph_\pm $. 
Thus we claim that, 
\[ 
 \frac{ \hph_{\pm } (\rho ) } {\hph_\pm ( \mu )
 } \leq 
C_1 \left \langle
\frac{ \rho - \mu }{\sqrt \epsilon } \right \rangle^2 \,. \]
Since $ \hph_\pm \sim d ( \bullet , \Gamma_\pm )^2 + \epsilon $, 
$ \hph_\pm \geq \epsilon $, we have 
\[ \begin{split} 
\hph_\pm ( \rho )  & \leq C ( d ( \rho , \Gamma_\pm )^2 + \epsilon ) 
\leq C ( d ( \mu, \Gamma_\pm)^2  + | \mu - \rho|^2 + 
\epsilon ) \\
& \leq C' ( \hph_\pm ( \mu )  + | \mu - \rho|^2 
) = C' ( \hph_\pm ( \mu ) 
 + \epsilon \langle ( \rho - \mu )/ \sqrt \epsilon \rangle^2 ) \\
& \leq 2 C'   \hph_\pm ( \mu )   
  \langle ( \rho - \mu )/ \sqrt \epsilon \rangle^2 \,.
\end{split}
\]
In the notation of Lemma \ref{l:es2}, \eqref{eq:ordph} gives
\[ | \widehat G ( \rho ) - \widehat G ( \mu ) | \leq C + 2  \log 
\langle ( \rho - \mu ) / \sqrt \epsilon \rangle \,,\]
and with $ \widehat \chi \in \CIc $, 
\[ | \widehat \chi ( \rho ) \widehat G ( \rho ) - 
\widehat \chi ( \mu ) \widehat G ( \mu ) | \leq 
C | \rho - \mu | \log ( 1 / \epsilon ) + C \log \langle ( \rho - \mu ) /
\sqrt \epsilon \rangle \,. \]
Clearly,
\[ | \chi_0 ( \rho ) G_0 ( \rho ) - \chi_0 ( \mu )  G_0 ( \mu ) | 
\leq C | \rho - \mu | \log ( 1/ \epsilon ) \,, \]
and hence to obtain \eqref{eq:orderf} 
we need 
\[ | \rho - \mu | \log ( 1 / \epsilon ) \leq C  \log \langle ( \rho - \mu ) /
\sqrt \epsilon \rangle + C  \,, \  \ \rho, \mu \in Q \Subset \RR^{2n} \,.\]
If we put $ \delta = \sqrt \epsilon $, 
$ t = | \rho - \mu |/ ( C\delta) $ this becomes
\[  \delta \log \frac{1}{\delta } 
\leq \frac{\log \langle t \rangle + 1  }t   \,,  \ \ 0 \leq t \leq 
\frac1{\delta} \,, \]
and that is clear as $ t \mapsto ( \log  \langle t \rangle  + 1) / t $
is decreasing.
\end{proof}

\section{Proof of the main result}
\label{pfo}

Let $ G $ be the escape function given in 
Proposition \ref{p:es2}, $ \epsilon = h/ \tilde h $. and 
let $ G^w $ be its Weyl quantization, 
\[  G^w = {\mathcal O} ( \log ( \tilde h/ h )) \; : \; L^2 ( X ) 
\rightarrow L^2 ( X ) \,.\]
 We define a family of conjugated operators:
\begin{equation}
\label{eq:ptt}  P_{\theta, t } \stackrel{\rm{def}}{=} e^{ -t G^w } P_\theta 
e^{ t G^w } \,, \ \ \theta = C_0 h \log ( 1/h ) \,. \end{equation}
It is easy to see that, in the notation of \S \ref{sma}, 
\begin{equation}
\label{eq:triv}
 \exp ( t G^w ) \in \Psi^{|t|C, 0 , 0 }_{\frac12} ( X) \,, 
\end{equation}
that is $ \exp ( t G^w ) = B_t^w $, $ \partial^\alpha B_t = 
{\mathcal O} ( h^{-|t|C - |\alpha|/2 } \tilde h^{|\alpha|/2} ) $.
Finer estimates are however possible thanks to the results of
Bony-Chemin \cite{BoCh}. The first of these is given in 
\begin{lem}
\label{l:cG}
Suppose that $ Q \in  \Psi_{\frac12}^{0,0,0} ( \RR^n ) $. Then
\begin{equation}
\label{eq:ficl}
 \exp( {-t G^w }) Q \exp( t G^w ) 
\in \Psi_{\frac12}^{0,0,0} ( \RR^n ) \,. \end{equation}
\end{lem}
\begin{proof}
We follow \S \ref{sma} and change to the
variables
\begin{gather*} 
 ( \tilde x , \tilde \xi ) = ( \tilde h / h )^{\frac12} ( x , 
\xi ) \,, \\
\widetilde G ( \tilde x , \tilde \xi ) =  G ( x , \xi ) \,, \ \ 
\widetilde Q_t (  \tilde x, \tilde \xi) = Q_t ( x , \xi ) \,, \\
U^{-1} G^w ( x , h D) U = \widetilde G^w ( \tilde x , 
\tilde h D_{\tilde x } ) \,, \ \ 
U^{-1} Q_t^w ( x , h D) U = \widetilde Q_t^w ( \tilde x , 
\tilde h D_{\tilde x } ) \,, \\
U v ( \tilde x ) = ( \tilde h / h )^{\frac{n}4} v ( ( h/ \tilde h)^{\frac12}
\tilde x ) \,.
\end{gather*}
We also note that 
\[  R \in \Psi_{\frac12}^{0,0,0} ( \RR^n ) \ \Longleftrightarrow \ 
U^{-1} R U \in \Psi^{0,0} ( \RR^n ) \,, \]
where on the right, $ \tilde h $ is the small parameter -- see the 
proof of Lemma \ref{l:1hf}.
The estimate \eqref{eq:orderf} shows that, in $ ( \tilde x, \tilde \xi) $
coordinates, $ \widetilde G $ satisfies the hypothesis of 
Proposition \ref{p:bbc} and that proves \eqref{eq:ficl}. 
\end{proof}

The basic properties of $ P_{t, \theta } $ are given in 
\begin{prop}
\label{p:propptt}
Let $ P_{\theta, t} $ be given by \eqref{eq:ptt} and let $ \Sigma 
\Subset T^* X $ be a compact surface coinciding with $ p^{-1} ( 0 ) $
in a neighbourhood of the support of $ G $. 
Then for $ |t | \leq C $, 
\begin{gather}
\label{eq:p.ptt} 
\begin{gathered}
P_{\theta , t }  = 
P_\theta - i t  h \Op (H_p G) + E_t \,,  
\ \ 
P_\theta - i t  h \Op (H_p G) \in 
 \Psi_{ \Sigma, \frac12}^{0,0,2,0} ( X) \cap 
\Psi_{ \frac12} ^{0,0,2} ( X )\,, 
\\ 
E_t \in  \Psi^{-1,-1,0}_{ \frac12} ( X ) \,, \ \ 
E_t = {\mathcal O} ( h \tilde h ) \; : \; L^2 ( X) \; \longrightarrow \;
L^2 ( X ) \,,
\end{gathered}
\end{gather}
uniformly in $ h $ and $ \tilde h $.
\end{prop}
\begin{proof}
Let $ V_1 , V_2 $ be open neighbourhoods of $ \supp G $, 
\[ \supp G \subset V_1 \Subset \overline{V}_2 \Subset T^*X \,.\]
We first observe that if $ \Psi \in \Psi^{0,-\infty} ( X ) $ satisfies
\[ \WFh ( \Psi ) \subset V_2 \,, \ \ \WFh ( I - \Psi ) \subset 
\complement V_1 \,, \]
then 
\begin{equation}
\label{eq:8.cut}
 [\exp ( t G^w ),  \Psi] \in \Psi^{-\infty , -\infty } ( X ) \,, \ \ 
( I - \Psi) ( \exp ( t G^w ) - I ) \in \Psi^{-\infty , -\infty } ( X ) \,,
 \ \ | t| \leq 1 \,. \end{equation}
In fact, using the calculus in \S \ref{sma} we see that
$ [ G^w , \Psi ] \in  \Psi^{-\infty , -\infty } ( X ) $, Hence, 
using \eqref{eq:triv}
\[ \begin{split} 
\frac{d}{dt}  [\exp ( t G^w ),  \Psi] & =  G^w  [\exp ( t G^w ),  \Psi] 
+ [ G^w , \Psi ] \exp ( t G^w )  \\
& = G^w  [\exp ( t G^w ),  \Psi]  + A_t \,, \ \  
A_t \in  \Psi^{-\infty , -\infty } ( X )  \,. \end{split} \]
Thus
\[ [ \exp ( t G^w) , \Psi ] = \int_0^t \exp( ( t -s ) G^w ) A_s ds 
\in 
  \Psi^{-\infty , -\infty }  ( X ) \,, \] 
which is the first statement in \eqref{eq:8.cut}. 
We also compute
\[ \frac{d}{dt} ( I - \Psi) ( \exp ( t G^w ) - I ) = 
( I - \Psi ) G^w \exp ( t G^w ) \in 
\Psi^{-\infty , -\infty } ( X )  \,, \]
and the second statement in \eqref{eq:8.cut} follows. Treating 
the equivalence of 
$ ( I - \Psi) P_\theta e^{t G^w } $ and $ ( I - \Psi ) P_\theta $
similarly we conclude that
\[ P_{\theta,t} -  e^{ -t G^w } \Psi P_\theta e^{ t G^w } -
( I - \Psi ) P_\theta \in \Psi^{-\infty, - \infty } ( X ) \, . \]
We now put
\[ Q_\theta \stackrel{\rm{def}}{=}  \Psi P_\theta \in 
\Psi^{0,0} (X ) \,, \ \ Q_{\theta, t} 
 \stackrel{\rm{def}}{=} e^{-t G^w} Q_\theta e^{t G^w } \,,
 \]
and we only need to prove \eqref{eq:p.ptt} with $ P_\bullet $ 
replaced by $ Q_\bullet $.
By a localization argument similar
to the one used to construct $ Q_\bullet $ we can assume that
$ X = \RR^n $ when applying Lemma \ref{l:cG} and that shows that
\[ Q_{\theta , t} \in \Psi_{\frac12}^{0,0,0} ( X ) \,. \]

We now establish the expansion in \eqref{eq:p.ptt}. Lemma 
\ref{l:Sjnew} implies that 
\[  [ Q_\theta, G^w ] = ( h/ i )  \Op ( H_{p_\theta}
G ) + R \,, \]
where $ R \in \Psi^{-3/2,-3/2,0}_{ \frac12} ( X ) 
\subset \Psi^{-1,-1,0}_{\frac12} ( X )  $. It also shows that
\[ [[ Q_\theta, G^w ], G^w]  = ( h/ i )  [ \Op ( H_{p_\theta}
G ) , G^w] + [ R , G^w ] \in \Psi^{-1,-1,0}_{\frac12} ( X ) \,.\]
Here we used the special structure of $ G$, 
\[ G = \widehat \chi \widehat G + C_0 \log ( 1/ h ) \chi_0 G_0 \,, \]
where $ \widehat \chi , \chi_0 $ and $ G_0 $ are uniformly smooth.
When derivatives fall on these 
terms in error estimates \eqref{eq:Sjnew1} 
the gain in $ h$ compensates for the logarithmic growth, while 
for $ |\alpha | > 0 $, $ \partial^{\alpha} \widehat G \in 
S^{|\alpha|/2, -|\alpha|/2}_{\frac12} $.

This gives,
\[
 \frac{d}{dt} E_t = [ Q_{\theta, t}, G^w] - ( h/i ) \Op ( H_{p_\theta}
G ) + ( h / i ) \Op ( H_{p_\theta - p } G ) = [ Q_{\theta, t } - Q_\theta, G^w ] + R_t   \,, \]
with $$ E_0 =  
 ( h / i ) \Op ( H_{p_\theta - p } G ) \in ( h \log( 1/ h ) )^2
\Psi^{0,0,0}_{\frac12} ( X )  \subset  \Psi^{-1,-1,0}_{\frac12} ( X ) 
\,, $$  
and $ R_t \in \Psi^{-1,-1,0}_{\frac12} ( X )  $.
We also have
\[  \frac{d}{dt} [ ( Q_{\theta, t } - Q_\theta ) , G^w] = 
e^{-tG^w }[ [ Q_{\theta} , G^w ] , G^w] e^{ t G^w } \in 
\Psi^{-1,-1,0}_{\frac12} 
( X) \,, \ \ Q_{ \theta, 0} - Q_\theta = 0 \,. \]
Hence $ [ Q_{\theta, t } - Q_\theta , G^w ] 
\in  \Psi^{-1,-1,0}_{\frac12}  
( X) $, and consequently $ E_t \in \Psi^{-1,-1,0}_{\frac12} $. 

To show that 
\[  Q^0_{ \theta, t }
\stackrel{\rm{def}}{=} Q_\theta - t h i  \Op ( H_{p} G ) 
 \in \Psi_{\Sigma, \frac12}^{0,0,0,0} ( X) \,, \]
it suffices to show, in view of Definition 2, that
for any sequence
\[ \{ a_j  \}_{j=1}^M \subset S^0 ( T^* X )\,, \ 
\] 
we have 
\[ \| \ad_{\Op (a_1)} \circ \cdots \ad_{\Op (a_M ) }
\circ \ad_P^k  Q^0_{\theta , t}   u \|_{ 
 L^2 ( X )  } \\ \ \leq \; C 
h^{ M/2 + k   } \tilde h ^{M/2  }  \| u \|_{L^2  ( X) } \,.
\] 
This will follow if we show that 
\[ \ad_P^k Q^0_{\theta, t } 
\in \Psi^{-k,0,0}_{\frac12} ( X )\,,
\]
and since that is clear for $ Q_\theta$
we only need to check that
\begin{equation}
\label{eq:adPQ} h \ad_P^k \Op ( H_{p} G )
\in \Psi^{-k,0,0}_{\frac12} ( X )\,.
\end{equation}

This follows from the following stronger result
\begin{lem}
Let $ G $ and $ P $ be as above. Then for any $ \epsilon > 0 $
\[  \ad_P^\ell ( H_p G)^w \in \Psi^{-\ell+\epsilon, 0, 0}_{\frac12} ( X ) \,, 
\ \ \ell \in {\mathbb N } \,, \ \ \ell \geq 1 \,. \]
\end{lem}
\begin{proof}
We start by 
proving that
\begin{equation}
\label{eq:Sjnew2}  \ad_P^\ell ( H_p G_1) ^w  \in 
\Psi^{-\ell+ \epsilon,0,0}_{\frac12} ( X ) \,, \ \ \ell \geq 1 \,, 
\end{equation}
where $ G_1 = \widehat G \widehat \chi $ is given in Proposition \ref{p:es2}.
We claim that 
$  \ad_P^\ell ( H_p G_1)^w = E_\ell^w $, where 
\[ E_\ell \sim \left( \frac{h}{i} \right)^\ell 
\left(  H_p ^{\ell+1} G_1 + 
\sum_{s=0}^{\ell-1} \sum_{r = 2(\ell - s)   }^\infty h^{r } V^{\ell}
_{rs} ( H_p^{s+1} G_1 ) \right)  \,, \]
with $ V_{rs}^\ell $ a differential operator of order less than or equal to 
$ r + \ell - s $, and 
where in view of good symbolic properties of $ P$ the 
error is $ {\mathcal O} ( h^\infty ) $.  In fact, Lemma \ref{l:Sjnew} gives for any 
$ a \in S^{0,0,0}_{\frac12} $, $ \ad_P a^w = a_1^w $,
\[ a_1 \sim (h/i) \left( H_p a + \sum_{r=2}^\infty h^{r} V_r ( a ) \right) 
\,, \]
where $ V_{r} $ is a differential operator of order $ r +1 $. 
The expansion for $ E_\ell $ comes from iterating this and observing that
for any differential operator $ B $ of order $ q $ 
\begin{equation}
\label{eq:obs}
 H_p^m B  H_p^k G_1 
= \sum_{s=0}^m B_s H_p^{k+s} G_1 \,, 
\end{equation}
where each $ B_s $ is a differential operator of order $ q $. 
Using \eqref{eq:es4} we now see
that for $ r \geq 2 ( \ell - s) $,
\[ h^\ell h^{r  } V^\ell_{rs} ( H_p^{s+1} G_1) \in 
h^{\ell +r} S_{\frac12}^{ ( r  + \ell - s  )/2 + \epsilon , -( r + \ell -s)/2  ,0} \subset 
S_{\frac12}^{-\ell - r/2 + \epsilon + ( \ell - s )/2   , 0 , 0 } 
\subset  S_{\frac12}^{-\ell - 1/2 + \epsilon  , 0 , 0 } \,, \]
where the $ \epsilon > 0 $ correction in the order comes from the term $ 
( H_p^{s+1} \widehat \chi ) \widehat G \in S^{\epsilon, 0 , 0 }_{\frac12} $.
Hence
\[  E_\ell  = ( h / i )^\ell (H_p^{\ell+1} G_1)^w + {\mathcal O} ( h^{\ell + 
\frac12} ) \Psi_{\frac12}^{\epsilon ,0,0} ( X ) 
\in \Psi^{-\ell + \epsilon ,0,0}_{\frac12} ( X ) \,, \]
which gives \eqref{eq:Sjnew2}. 
On the other hand, again in the notation of Proposition \ref{p:es2},
\[ \ad_P^\ell (H_p G_0) ^w \in \Psi^{-\ell,0} ( X ) \,, \]
and consequently 
\[ \ad_P^\ell ( \log ( 1/h ) H_p G_0 )^w  \in 
\log( 1/ h ) \Psi_{\frac12}^{-\ell, 0 , 0 } ( X ) 
\subset \Psi_{\frac12}^{-\ell + \epsilon  ,0 , 0} ( X ) \,. \]
Since $ G = G_1 + \log ( \tilde h / h ) \chi_0 G_0 $ the lemma follows.
\end{proof}

\renewcommand\thefootnote{\P}%
The lemma immediately gives \eqref{eq:adPQ} completing the proof of 
Proposition \ref{p:propptt}.
\end{proof}
The next lemma follows from Proposition \ref{p:es2}:
\begin{lem}
\label{l:chHpG}
Let $ \widehat 
G $ be given in Lemma  \ref{l:es2}, $ \chi \in \CIc ( 
\RR )$, and let $ \psi \in \CIc ( T^* X ) $ be one in a fixed 
small 
neighbourhood of $ K $ and zero outside of another sufficiently
 small neighbourhood
of $ K $. Then 
\[ \psi ( x, \xi) \chi( H_p \widehat G ( x, \xi ) ) 
\in S^{0,0,-\infty,0}_{ \Sigma , \frac12} ( T^*X ) \,.\]
\end{lem}
\begin{proof}
Lemma \ref{l:es2} gives  a stronger condition 
\[  H_p^\ell \partial^\alpha_\rho H_p^k ( \psi ( \rho) \chi ( H_p G (\rho ) )
= {\mathcal O} ( ( \tilde h / h )^{|\alpha|/2} ) \,,\]
as can be verified using \eqref{eq:obs}
\end{proof}

As in \S \ref{pt2} we modify our operator to obtain global invertibility.
Thus we define $ a \in S_{\Sigma, \frac12}^{0,0,-\infty,-\infty} ( T^* X ) 
$ as follows
\begin{gather*}
 a ( x , \xi ) \stackrel{\rm{def}}{=}
 \chi \left( \frac{\tilde h}{ h } p ( x , \xi ) \right) 
\chi ( H_p G ( x, \xi )  ) \psi  (x, \xi)   \,, \\ 
\chi \in \CIc ( \RR ; [0,1] ) \,, \ \ \chi ( t ) \equiv 1 \,,  \ \ 
|t | \leq 1 \,, 
\end{gather*}
and $ \psi $ is as in Lemma \ref{l:chHpG}. In particular by taking its
support close to $ K $ we can replace $ G $ by $ \widehat G $ in the
definition of $ a $. 

We then put
\begin{equation}
\label{eq:ptht}
 \widetilde P_{\theta, t} =  P_{\theta, t}  - i (h/\tilde h) \Opt( a ) 
\in \Psi^{0,0,2,1}_{\Sigma , \frac12} ( X ) \,,\end{equation}
and first treat the region away from the trapped set:
\begin{lem}
\label{l:e1}
Suppose that $ P_{\theta, t } $ is given by \eqref{eq:ptht} and
$ \Psi_0 \in \Psi^{0,0} ( T^* X) $ satisfies
$$ \WFh ( \Psi_0 ) \cap K  = \emptyset \,. $$
Then for $ u \in \CIc ( X ) $,  $ z \in D ( 0 , Ch ) $  we have 
\begin{gather*}
 \| ( \widetilde P_{\theta,t} - z ) \Psi_0 u \|_{L^2 } \geq  t h \| \Psi_0
u \|_{L^2 ( X ) }/ C - 
{\mathcal O} ( h^\infty ) \| u \|_{L^2 (  X) }  \,, \\ 0 < h \leq h_0 ( \tilde h ) \,, \ \ 0 <  \tilde h
\leq \tilde h_0 ( t ) \,. 
\end{gather*}
\end{lem}
\begin{proof}
Let us assume that $ \| u \|=1 $.
Once $ h $ is small enough, 
$ a \equiv 0 $ in a neighbourhood of $ \WFh ( \Psi_0 ) $
and Theorem \ref{th:3.1} gives
\[ \| ( \widetilde P_{\theta,t} - z ) \Psi_0 u \|_{L^2 } =
\| ( P_{\theta, t} - z ) \Psi_0 u \|_{L^2 } + {\mathcal O} (h  \tilde 
h^\infty ) 
\,. \]
We then observe that $ P_{\theta, t} \in \Psi^{0,0,2}_{1/2} 
( X ) $, and consequently we can use the simpler calculus of 
\S \ref{sma}. Microlocally 
near $ \WFh ( \Psi_0 ) $,  for $ z \in D ( 0 , C h ) $, and for
$ t $ sufficiently 
large,  Proposition \ref{p:es2} and the choice of the angle of 
scaling give,
\begin{gather*}
P_{\theta , t } - z = \Op ( \Re p_\theta - \Re z ) + i 
 \Op ( \Im p_\theta  - i h t H_{p} G - \Im z ) 
+ {\mathcal O}_t ( h \tilde h + h^2 \log ( 1/ h ) )\,, \\
| \Re p_\theta - \Re z | < \delta \ \Longrightarrow \ 
-\Im p_\theta + h t H_p G + \Im z \geq t h / C \,.
\end{gather*}
(This is the analogue of \eqref{12} in the non-trapping case of 
\S \ref{rfr}.)
Lemma \ref{l:r3} applied with $ \Psi_j $'s such that  
$ | \Re p_\theta - \Re z | > \delta $ on $\WFh ( \Psi_1 ) $
(with $ \Psi_j$'s constructed using Lemma \ref{l:r2}) 
completes the proof.
\end{proof}

Near the trapped set we  use the second microlocal calculus to 
obtain
\begin{lem}
\label{l:e2}
Suppose that $ P_{\theta, t } $ is given by \eqref{eq:ptht} and
let $  z \in D ( 0, Ch ) $.
For $ u \in \CIc ( X ) $, $ \|u \| = 1$, with $ \WFh ( u ) $ in a fixed small
neighbourhood of $ K $ we have 
\begin{equation}
\label{eq:l.e2}
\begin{split}
& \|  ( \widetilde P_{\theta,t} -z ) u \|_{L^2 ( X ) } 
\geq t h \| u \|_{L^2 ( X ) } / C 
\,, \ \ 0 < h \leq h_0 ( \tilde h ) \,, 0 <  \tilde h
\leq \tilde h_0 ( t ) \,. 
\end{split}
\end{equation}
provided that $ t $ is large enough.
\end{lem}
\begin{proof}
In a small neighbourhood of $ K $ the operator is microlocally equal to 
\[   P^\flat_t \stackrel{\rm{def}}{=}  P
 - i t h \Opt ( H_p \widehat G ) - i ( h/ \tilde h  ) \Opt ( a ) 
+ {\mathcal O}_{ L^2 \rightarrow L^2 } ( h \tilde h) \,,\]
that is,
\[  \| (\widetilde P_{\theta, t } - z ) u \|_{ L^2 ( X) } 
 = \| ( P^\flat_t - z ) u \|_{ L^2 ( X ) } + {\mathcal O} 
( h ^\infty ) \,, \ \ \| u \|_{L^2 ( X ) } = 1 \,, \]
for $ u $ with $ \WFh ( u ) $ near $ K $.
For $ z \in D( 0 , C h ) $,
\[ P^\sharp_t - Z
 \stackrel{\rm{def}}{=}
( \tilde h / h ) ( P^\flat_t - z ) \in \Psi^{0,0,2,1}_{\Sigma , 
\frac12} (  X) \,, \ \  Z \stackrel{\rm{def}}{=} ( \tilde h / h ) z \,, \]
has the symbol given by 
\[ p^\sharp - Z= \lambda - Z - i t \tilde h  H_p G  - i \chi ( \lambda) 
\chi ( H_p G ) + {\mathcal O} ( \tilde h^2 ) \,, \ \  
Z \in D ( 0 , C \tilde h ) \,, \ \ \lambda = (\tilde h / h) p \,.\]
Now let $ \psi_0, \psi_1 \in \CI_{\rm{b}} ( \RR ) $ satisfy
\[ \psi_0 ( t)^2 + \psi_1^2 ( t ) = 1 \,, \ \  \supp \psi_0 \subset 
\{ t \; : \; \chi ( t)  = 1 \}\,, \ \ \psi_1 ( t ) \equiv 0 \,, 
\ \ |t| \leq 1/2 \,. \]
As in Lemma \ref{l:r2} we can now find two operators 
$ \Psi_j^\sharp \in \Psi^{0,0,0,0}_{\Sigma, \frac12} ( T^* X ) $ such 
that 
\[ \sigma_{h, \tilde h } ( \Psi_j^\sharp ) = \psi_j ( H_p G ) \,, 
\ \ (\Psi_0^\sharp)^2 + (\Psi_1^\sharp)^2 = Id + {\mathcal O}_{L^2 
\rightarrow L^2 } ( \tilde h^\infty ) \,. \]
In a neighbourhoood of the support of $ \psi_1 ( H_p  G) $ the 
operator $ P_\theta^\sharp - Z $ is elliptic in 
$ \Psi^{0,0,2,1}_{\Sigma , \frac12} $: 
\[ | \lambda - Z - i \tilde h H_p G - i \chi ( \lambda ) \chi ( H_p G ) |
\geq ( | \lambda - Z | + \chi ( \lambda ) \chi ( H_p G ) )/2 \geq 
\langle \lambda \rangle/ C \,, \]
when $ Z \in D ( 0 , C \tilde h ) $ and $ \chi ( H_p G ) > 1/2 $, say.
This implies (see Lemma \ref{l:r1}) that for $ u $ with $ \WFh ( u ) $
near $ K $, $ \| u \|_{ L^2 ( X ) } = 1 $,
\[  \| ( P^\sharp_t - Z ) \Psi_0^\sharp  u \|_{L^2 ( X) } \geq \| \Psi_0^\sharp 
u \| _{L^2 ( X ) } /C - {\mathcal O} ( \tilde h^\infty )  \,, \ 0 < \tilde h \leq \tilde h_0 \,. \]
To estimate $ \| ( P^\sharp_t - Z ) \Psi_1^\sharp  u \| $ from 
below we proceed as in \S \ref{cop}. Let 
\[ B^\sharp_t = \frac1{2i} \left( P^\sharp_t - ( P^\sharp_t)^* \right) \,, \]
so that the Weyl symbol (in the sense of 
$ \Psi_{\Sigma, \frac12}$) of $ B^\sharp_t $ is equal to 
\[ - \chi ( \lambda ) \chi ( H_p G ) + {\mathcal O} ( \tilde h ) 
- t \tilde h H_p G  + {\mathcal O}_t (  \tilde h^2 ) \,,\]
where we indicated the dependence on $ t $ in the second bound. 
Since
$ H_p G \geq 1/C $ in a neighbourhood of 
the support of $ \psi_1 ( H_p ) $ we see  that for $ u $ with $ \WFh ( u ) $
near $ K $, $ \| u \|_{ L^2 ( X) } $,
\[ - \langle B_t^\sharp \Psi^\sharp_1 u ,  \Psi^\sharp_1 u \rangle
\geq ( t \tilde h  - {\mathcal O}( \tilde h ) 
- {\mathcal O}_t ( \tilde h^2 ) )  \| \Psi^\sharp_1 u  \|^2 
-  {\mathcal O} ( \tilde h^\infty )  \,.\]
We now first take $ t $ large enough to dominate the first error
term and then $ \tilde h $ small enough to dominate the second one.
Hence,
\begin{equation*}
\begin{split}
\| ( P^\sharp_t - Z) \Psi_1^\sharp u \| \| \Psi_1^\sharp u \| 
& \geq | \langle ( P^\sharp_t - Z )  \Psi_1^\sharp u , 
\Psi_1^\sharp u \rangle | 
\geq | \Im  \langle ( P^\sharp_t - Z )  \Psi_1^\sharp u , 
\Psi_1^\sharp u \rangle |\\
& = - \langle  ( B^\sharp_t - \Im Z ) 
\Psi_1^\sharp u , \Psi_1^\sharp u \rangle  \geq 
t {\tilde h} \| \Psi_1^\sharp u \|^2 /2  -  {\mathcal O} ( \tilde h^\infty )  \,,
\end{split}
\end{equation*}
provided that $ t $ was large enough, and then $ \tilde h $ 
small enough. Lemma \ref{l:r3} (or rather its proof) gives
\[ \| ( P^\sharp_t - Z) u \| \geq t \tilde h \| u \|/C \,, 
t \geq t_0 \gg 1 \,, \ \ 0 < h< h_0 ( t) \,, \]
for $ u $ with $ \WFh( u ) $ near $ K $.

We complete the proof by writing
\[ \begin{split}
 \| (\widetilde P_{\theta, t } - z ) u \|_{ L^2 ( X) } 
& = \| ( P^\flat_t - z ) u \|_{ L^2 ( X ) } + {\mathcal O} 
( h ^\infty ) \| u \|_{L^2 ( X ) } \\
& = ( h / \tilde h ) \| ( P^\sharp _t - Z ) u \|_{ L^2 ( X ) } 
 + {\mathcal O} 
( h ^\infty ) \| u \|_{L^2 ( X ) } \\ & \geq t h \| u \|_{ L^2 ( X ) }/ C \,.
\end{split}\]
\end{proof}
The two lemmas are now combined using Lemma \ref{l:r3} 
which gives for large $ t $, $ 0 <  \tilde h 
\leq \tilde h_0 ( t ) $, and $ 0 < h < h_0 ( t ,\tilde h ) $,
the invertibility of $ \widetilde P_{\theta, t } - z $, $ z 
\in D( 0 , C h ) $:
\[ ( \widetilde P_{\theta , t } - z )^{-1} = {\mathcal O} ( 
1/h ) \; : \; L^2 ( X ) \; \longrightarrow \; L^2 (X ) \,. \]

As in \S \ref{pt2}, Theorem \ref{t:3} 
is a consequence of writing
\begin{equation}
\label{eq:fit}
  \Opt( a )  = R + E \,, \ \ \rank ( R ) = {\mathcal O} ( h^{-\nu} )\,,\ \ 
E = {\mathcal O} ( \tilde h^\infty ) \; : \; L^2 ( X ) \rightarrow 
L^2 ( X ) \,, \end{equation}
$ \nu > \nu ( E ) $, where $ m( E ) = 2 \nu ( E ) + 1 $ is the dimension
of the trapped set at energy $ E $, allowing $ \nu = \nu( E )$ if
the trapped set is of pure dimension.

The decomposition \eqref{eq:fit}
 follows from Proposition \ref{p:finite} and the definition of the
Minkowski dimension:
\[ m_0 = 2 n - 1  - \sup\{ d \; : \; \limsup_{ \epsilon \rightarrow 0 }
\epsilon^{-d} \vol ( \{ \rho \in p^{-1} ( 0 ) \; : \; d ( \rho , K ) < 
\epsilon \} ) < \infty \} \,,\] 
with the set being of pure dimension if
\[  \limsup_{ \epsilon \rightarrow 0 }
\epsilon^{-2n + 1 + m_0} 
\vol ( \{ \rho \in p^{-1} ( 0 ) \; : \; d ( \rho , K ) < 
\epsilon \} ) < \infty \,. \]
In other words, for $ \epsilon $ small
\[  \vol ( \{ \rho \in p^{-1} ( 0 ) \; : \; d ( \rho , K ) < 
\epsilon \} ) \leq C \epsilon^{ 2n - 1 - m } \,, \ \ m > m_0 \,, \]
and $ m $ replaceable by $ m _0 $ when $ K $ is of pure dimension.
In particular, 
\[ \vol ( \supp a \cap p^{-1} (0) ) \leq C_{\tilde h} h^{ ( 2n - 1 - m )/2 } 
= C_{\tilde h } h^{ n - \nu - 1 } \,,  \ \ m = 2 \nu + 1 > m_0\,, \]
with equality if $ K $ is of pure dimension.
Since 
\[ \supp a \cap p^{-1}(0) \subset \bigcup_{ \rho \in K } B_\Sigma ( \rho, 
M ( h/ \tilde h )^{\frac12}) \,,  \]
where $ B_\Sigma $ are balls in $ \Sigma $ with respect to some fixed
smooth metric, and since $ K $ is invariant under the flow,
the standard covering arguments (see \cite[Lemma 3.3]{SjDuke}) 
show that the hypothesis
of Proposition \ref{p:finite} are satisfied with 
$$ K (h ) \leq C_{\tilde h} h^{-\nu}\,,  $$
which completes the proof of Theorem \ref{t:3}.

\vspace{0.5cm}
\noindent
{\bf Appendix}
\vspace{0.4cm}
\renewcommand{\theequation}{A.\arabic{equation}}
\refstepcounter{section}
\renewcommand{\thesection}{A}
\setcounter{equation}{0}

\renewcommand{\Op}{{\operatorname{Op}^{{w}}}}
\newcommand{\Opp}{{\operatorname{Op}^{{w}}_{\Phi_0}}}

We present a direct proof of Proposition \ref{p:bbc}. 
The hypotheses on $ G$ in \eqref{eq:bc1} are equivalent to 
the statement that
$  \exp ( t G ) \in S ( m^t ) $, for all $ t \in \RR $.
We start with 
\begin{lem}
\label{l:a1}
Let $ U ( t) \stackrel{\rm{def}}{=} ( \exp t G )^w ( x , D ) 
\; : \; {\mathcal S} ( \RR^n ) \; \longrightarrow \; 
 {\mathcal S} ( \RR^n ) $. For $ |t| < \epsilon_0 ( G) $, 
the operator $ U ( t ) $ is invertible, and
\[ U(t)^{-1} = B_t ^w ( x, D ) \,, \ \ B_t \in S( m^{-t} ) \,. \]
\end{lem}
\begin{proof}
We apply the composition formula \eqref{eq:usual} to obtain 
\begin{gather*}
 U ( - t) U ( t ) = Id + E_t^w ( x , D )  \,, \ \ E_t \in S ( 1) \,.
\end{gather*}
More explicitely we write (see \cite[Proposition 7.7]{DiSj} and 
Lemma \ref{l:Sjnew} here)
\[ \begin{split} E_t ( {x_1} , \xi ) & = 
\int_0^s e^{  s A ( D)}
A( D )  ( e^{- t G ( {x_1} , \xi_1) + 
 t G ( {x_2} , \xi_2 ) }  ) \rest_{ {x_2} = {x_1} = x 
 , \xi_2 = \xi_1 = \xi } ds \\
& =  \int_0^s (i t /2 )  e^{  s A ( D ) }
( D_{\xi_1} G D_{{x_2} } 
G  - D_{x_1} G D_{\xi_2} G  )
 e^{- t G ( {x_1} , \xi_1) + 
 t G ( {x_2} , \xi_2 ) }
\rest_{ {x_2} = {x_1} = x , \xi_2 = \xi_1 = \xi }  ds \,, 
\end{split}
\]
where $ A ( D ) = i \sigma ( D_{x_1}, D_{\xi_1}; D_{x_2} , D_{\xi_2} ) /2  $.

Hence $ E_t = t \widetilde E_t $ where $ \widetilde E_t \in 
S ( 1 ) $ uniformly, and thus
\[  E_t^w ( x, D ) = {\mathcal O} ( t) : L^2 ( \RR^n ) \rightarrow
L^2 ( \RR^n ) \,. \]
This shows that
 for $ |t| $ small enough $ Id + E_t^w ( x, D ) $ is invertible, 
and
Beals's lemma (see for instance \cite[Proposition 8.3]{DiSj})
gives 
\[  ( Id + E_t ^w ( x ,D) )^{-1} = C_t^w ( x, D ) \,, \ \ C_t 
\in S ( 1 ) \,. \]
 Hence $ B_t = C_t \# \exp ( - t G (  x, \xi ) ) \in 
S ( m^{-t } ) $.
\end{proof}
We now observe that 
\begin{gather}
\label{eq:a1}
\begin{gathered}
  \frac{d}{dt} \left( U (- t ) \exp ( t G^w ( x , D ) ) 
\right)  = V ( t) \exp ( t G^w  ( x , 
D ) ) 
\,, \\   V ( t )  = A_t^w ( x , D) \,, \ \ A_t \in S ( m^{-t} ) \,.
\end{gathered} \end{gather}
In fact, we see that 
\[ \frac{d}{dt} U ( - t ) = - ( G \exp ( - t G ) ) ^w ( x , D ) \,, \ \ 
 U ( - t ) G^w( x , D) = (  \exp ( t G ) \# G )  ^w ( x , D ) \,.\]
As before, the composition formula \eqref{eq:usual} gives 
\begin{gather*} 
 \exp ( -t G ) \# G -  G \exp ( - t G )  = \\
\int_0^1  \exp (s A ( D ) )  A( D ) 
\exp ( - t G ( x^1 , \xi^1 )  G  ( x^2 , \xi^2 )  \rest_{ x^1 = x^2 = x ,
\xi^1 = \xi^2 = \xi } \,, \\ 
 A ( D )  =  i \sigma ( D_{x^1}, D_{\xi^1} ; D_{x^2 } , D_{\xi^2} )/2 \,.
\end{gather*}
The hypothesis on $ G $ shows that $ A ( D) \exp ( t 
G ( x^1 , \xi^1 ) )  G  ( x^2 , \xi^2 ) 
 $ is a sum of terms of the form $ a ( x^1 , \xi^1 ) b ( 
x^2 , \xi^2 ) $ where $ a \in S ( m^{-t}  ) $ and $ b \in S ( 1 ) $. 
The continuity of $ \exp ( A ( D) ) $ 
on the spaces of symbols (see \cite[Proposition 7.6]{DiSj}) 
gives \eqref{eq:a1}. 

If we put
\[  C( t) \stackrel{\rm{def}}{=} - V( t) U ( -t)^{-1} \,,\]
then by Lemma \ref{l:a1}, $ C( t ) = c_t^w $ where $ c_t \in S( 1 ) $.
Symbolic calculus shows that $ c_t $ depends smoothly on $ t $ and
\[ (  \partial_t  + C ( t ) ) (U ( - t ) \exp ( t G^w ( x , D) )) = 0 \,.\]
The proof of Proposition \ref{p:bbc} is now reduced to showing 
\begin{lem}
\label{l:a2} 
Suppose that $  C( t) = c_t^w ( x, D)  $, where  $ c_t \in S ( 1 ) $, depends 
continuously  on $ t \in (-\epsilon_0, \epsilon_0) $. Then the solution 
of 
\begin{equation}
\label{eq:a2}
( \partial_t + C ( t) ) Q ( t) = 0 \,, \ \ Q ( 0 ) = q^w ( x, D) 
\,, \ \ q \in S( 1 ) \,,
\end{equation}
is given by $ Q ( t ) = q_t ( x, D ) $,  where $ q_t \in S( 1 ) 
$ depends 
continuously  on $ t \in (-\epsilon_0, \epsilon_0) $.
\end{lem}
\begin{proof} 
The Picard existence theorem for ODEs shows that $ Q ( t) $ is
bounded on $ L^2 $. If $ \ell_j  ( x , \xi )  $ are linear 
functions on $ T^* \RR^n $ 
then 
\begin{gather*}
 \frac{d}{dt} \ad_{ \ell_1 ( x, D ) } \circ \cdots
\circ \ad_{ \ell_N( x , D ) } Q ( t) +
\ad_{ \ell_1 ( x, D ) } \circ \cdots
\circ \ad_{ \ell_N ( x , D ) } ( C ( t) Q ( t) ) = 0 \,, 
\\ 
 \ad_{ \ell_1 ( x, D ) } \circ \cdots
\circ \ad_{ \ell_N( x , D ) } Q ( 0) \; : \; L^2 ( \RR^n ) 
\longrightarrow L^2 ( \RR^n ) \,. \end{gather*}
If we show that for any choice of $ \ell_j's $ and any $ N $
\begin{equation}
\label{eq:bea}  \ad_{ \ell_1 ( x, D ) } \circ \cdots
\circ \ad_{ \ell_N( x , D ) } Q ( t) 
 \; : \; L^2 ( \RR^n ) 
\longrightarrow L^2 ( \RR^n ) \,, \end{equation}
then Beals's lemma (see \cite[Chapter 8]{DiSj}) concludes the proof.
We proceed by induction on $ N $:
\[ \ad_{ \ell_1 ( x, D ) } \circ \cdots
\circ \ad_{ \ell_N ( x , D ) } ( C ( t) Q ( t) ) = 
 C( t) 
\ad_{ \ell_1 ( x, D ) } \circ \cdots
\circ \ad_{ \ell_N ( x , D ) } Q ( t) + R ( t) \,, \]
where $ R ( t ) $ is the sum of terms of the form 
\[  A_k  ( t)  \ad_{\ell_1 ( x, D ) } 
\circ \ad_{\ell_k  ( x , D)}  Q ( t) \,, \ \ 
k < N \,, \ \ A_k( t) = a_k( t)^w \,, \]
where $ a_k ( t) \in  S( 1 ) $ depend continuously on $ t $ (this
statement can also be proved by induction using the derivation property 
of $ ad_\ell $: $ \ad_{\ell } ( C D ) = ( \ad_\ell C) D + C ( \ad_\ell D ) $).
Hence by the induction hypothesis $ R ( t) $ is bounded on $ L^2 $,
and depends continuously on $ t $.
Thus
\[ 
\left(  \frac{d}{dt} + C ( t) \right) \ad_{ \ell_1 ( x, D ) } \circ \cdots
\circ \ad_{ \ell_N( x , D ) } Q ( t) = R ( t) \; : \; L^2 ( \RR^n )
\; \longrightarrow \; L^2 ( \RR^n ) \,. \]
Since  \eqref{eq:bea}  is valid at $ t = 0 $ we obtain it for
all $ t \in ( - \epsilon_0, \epsilon_0 ) $.
\end{proof}

\end{document}